\documentclass[11pt]{amsart}
\usepackage{
  amsmath,
  amsfonts,
  amssymb,
  graphicx, 
  times,
  color,
  hyphenat,
  pinlabel
}
\input xy
\xyoption{all}
\newcommand{\FKh}{\mathcal{F}_{sl(2),n}}
\newcommand{\Fmv}{\mathcal{F}_{sl(3),n}}

\newcommand{\n}{\noindent}

\newcommand{\cmw}{\mathbf{Foam}_{2}}
\newcommand{\foamt}{\mathbf{Foam}_{3}}

\newcommand{\figins}[3] 
{\raisebox{#1pt}{\includegraphics[height=#2 in]{figs/#3}}}
\newcommand{\figwhins}[4] 
{\raisebox{#1pt}{\includegraphics[height=#2 in, width=#3 in]{figs/#4}}}
\newtheorem{thm}{Theorem}[section]

\newtheorem{cor}[thm]{Corollary}
\newtheorem{prop}[thm]{Proposition}

\theoremstyle{definition}
\newtheorem{defn}[thm]{Definition}

\newcommand{\bZ}{\mathbb{Z}}

\newcommand{\bR}{\mathbb{R}}
\newcommand{\bC}{\mathbb{C}}

\catcode`\@=11
\long\def\@makecaption#1#2{%
    \vskip 10pt
    \setbox\@tempboxa\hbox{%
\small{#1: }\ignorespaces #2}%
    \ifdim \wd\@tempboxa >\captionwidth {%
        \rightskip=\@captionmargin\leftskip=\@captionmargin
        \unhbox\@tempboxa\par}%
      \else
        \hbox to\hsize{\hfil\box\@tempboxa\hfil}%
    \fi}
\newdimen\@captionmargin\@captionmargin=2\parindent
\newdimen\captionwidth\captionwidth=\hsize
\catcode`\@=12
\title{The diagrammatic Soergel category and $sl(2)$ and $sl(3)$ foams}
\author{Pedro Vaz}
\address{Institut de Math\'ematiques de Jussieu \\ Universit\'e Paris 7 \\   
175 Rue du Chevaleret, F-75013 Paris, France 
and  
CAMGSD\\Instituto Superior T\'{e}cnico\\ Avenida Rovisco Pais\\ 1049-001 Lisboa\\ Portugal}
\email{pfortevaz@ualg.pt}
\begin{document}
%
%
% this must be defined locally (to center captions)
\newdimen\captionwidth\captionwidth=\hsize
%
%%%%%%%%%%%%%%%%%%%%%%%%%%%%%%%%%%%%%%%%
%%%                                  %%%
%%%            Abstract              %%%
%%%                                  %%%
%%%%%%%%%%%%%%%%%%%%%%%%%%%%%%%%%%%%%%%%
%
\begin{abstract} 
We define two functors from Elias and Khovanov's diagrammatic Soergel category, one targeting Clark-Morrison-Walker's 
category of disoriented $sl(2)$ cobordisms and the other the category of (universal) $sl(3)$ foams.
\end{abstract}
\maketitle
%
%
%%%%%%%%%%%%%%%%%%%%%%%%%%%%%%%%%%%%%%%%
%%%                                  %%%
%%%        Introduction              %%%
%%%                                  %%%
%%%%%%%%%%%%%%%%%%%%%%%%%%%%%%%%%%%%%%%%
\section{Introduction}\label{sec:intro}

In this paper we define functors between the Elias-Khovanov diagrammatic version of 
the Soergel category $\mathcal{SC}$ defined in~\cite{EKh} and the categories of 
universal $sl(2)$ and $sl(3)$-foams defined in~\cite{CMW} and~\cite{MV1}.

The Soergel category provides a categorification of the Hecke algebra and was used by 
Khovanov in~\cite{Kh-bim} to construct a triply graded link homology categorifying the 
HOMFLYPT polynomial. 
Elias and Khovanov constructed in~\cite{EKh} a category defined diagrammatically by 
generators and relations and showed it to be equivalent to $\mathcal{SC}$.

The $sl(2)$ and $sl(3)$-foams were introduced in~\cite{BN, CMW} and in~\cite{Kh-sl3, MV1} 
respectively, to give topolo\-gical constructions of the $sl(2)$ and $sl(3)$ link homologies.

This paper can be seen as a first step towards the construction of a family of functors 
between $\mathcal{SC}$ and the categories of $sl(N)$-foams for all $N\in\bZ_+$, 
to be completed in a subsequent paper~\cite{MV2}. 
The functors $\FKh$ and $\Fmv$ are not faithful. In~\cite{MV2} we will 
extend the construction of these functors to all $N$. 
The whole family of functors is faithful in the following sense: 
if for a morphism $f$ in $\mathcal{SC}_1$ we have $\mathcal{F}_{sl(N),n}(f)=0$ for all $N$, 
then $f=0$. 
With these functors one can try to give a graphical interpretation of 
Rasmussen's~\cite{rasmussen} spectral sequences from the HOMFLYPT link homology to the 
$sl(N)$-link homologies.

The plan of the paper is as follows. In Section~\ref{sec:soergel} we give a brief 
description of Elias and Khovanov's diagrammatic Soergel category. 
In Section~\ref{sec:sl2} we describe the category $\cmw$ of $sl(2)$-foams and construct a 
functor from $\mathcal{SC}$ to $\cmw$. Finally in Section~\ref{sec:sl3} we give the 
analogue of these results for the case of $sl(3)$-foams.

We have tried to keep this paper reasonably self-contained. Although not mandatory, 
some acquaintance with~\cite{CMW, EKh,EKr,MV1} is desirable.

%%%%%%%%%%%%%%%%%%%%%%%%%%%%%%%%%%%%%%%%
%%%                                  %%%
%%%           Soergel                %%%
%%%                                  %%%
%%%%%%%%%%%%%%%%%%%%%%%%%%%%%%%%%%%%%%%%
\section{The diagrammatic Soergel category revisited}
\label{sec:soergel}

This section is a reminder of the diagrammatics for Soergel categories introduced by 
Elias and Khovanov in~\cite{EKh}. 
Actually we give the version which they explained in Section 4.5 of~\cite{EKh} and which 
can be found in detail in ~\cite{EKr}. 
%In a nutshell a diagrammatic category consists of
%planar diagrams which are morphisms between their bottom and top components, which in turn 
%are the objects of the category. Composition corresponds to concatenation of diagrams and 
%this category can be endowed with a monoidal structure by the disjoint union of diagrams.

Fix a positive integer $n$. The category $\mathcal{SC}_1$ is the category whose objects 
are finite length sequences of points on the real line, where each point is colored by an 
integer between $1$ and $n$. 
We read sequences of points from left to right. 
Two colors $i$ and $j$ are called adjacent if $\vert i-j\vert=1$ and distant if 
$\vert i-j\vert >1$.  
The morphisms of $\mathcal{SC}_1$ are given by generators modulo relations. 
A morphism of $\mathcal{SC}_1$ is a $\bC$-linear combination of planar diagrams 
constructed by horizontal and vertical gluings of the following generators 
(by convention no label means a generic color $j$):
\begin{itemize}
\item Generators involving only one color:
\begin{equation*}
\xymatrix@R=1.0mm{
\figins{-15}{0.5}{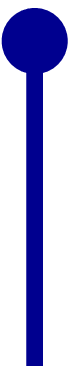}
&
\figins{-9}{0.5}{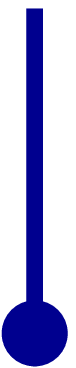}
&
\figins{-15}{0.55}{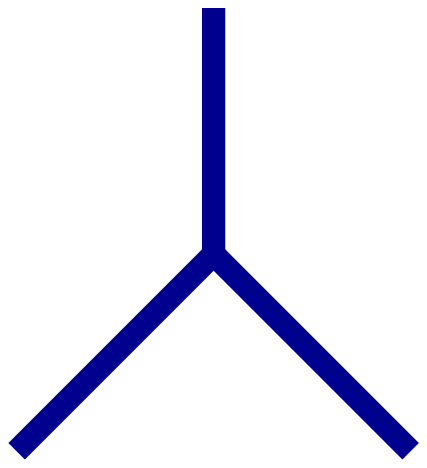}
&
\figins{-15}{0.55}{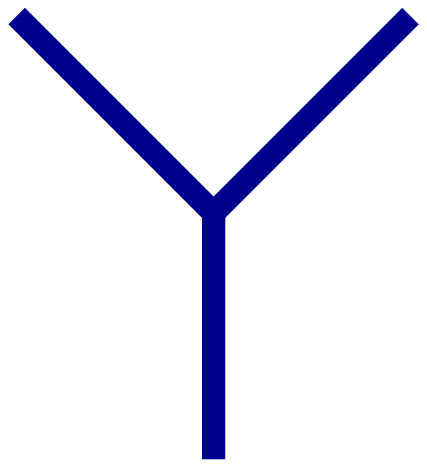}
\\
\text{EndDot} & \text{StartDot} & \text{Merge} & \text{Split}
}
\end{equation*}

\medskip

It is useful to define the cap and cup as
\begin{equation*}
\figins{-17}{0.55}{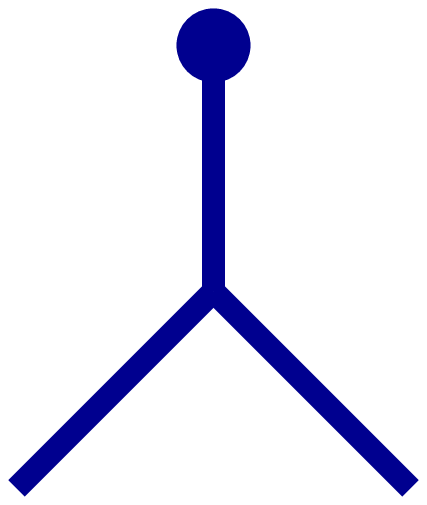}\
\equiv\
\figins{-17}{0.55}{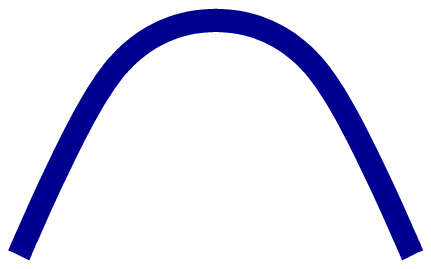}\
\mspace{50mu}
\figins{-17}{0.55}{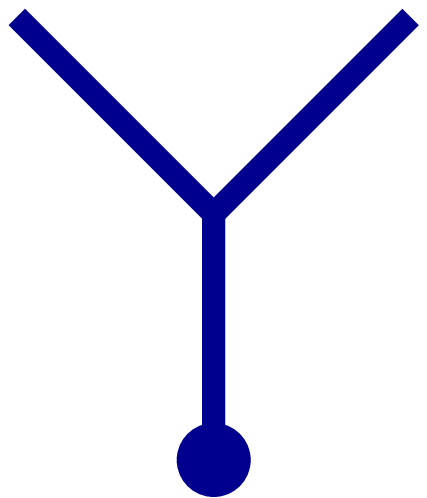}\
\equiv\
\figins{-17}{0.55}{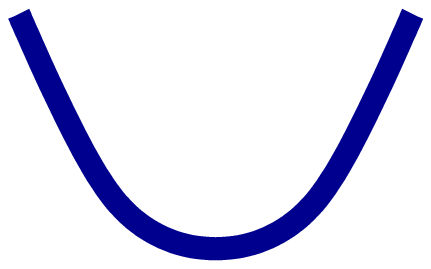}
\end{equation*}

\medskip

\item Generators involving two colors:
\begin{itemize}
\item The 4-valent vertex, with distant colors,
\begin{equation*}
\labellist
\tiny\hair 2pt
\pinlabel $i$ at  -4 -10
\pinlabel $j$ at 134 -10
\endlabellist
\figins{-15}{0.55}{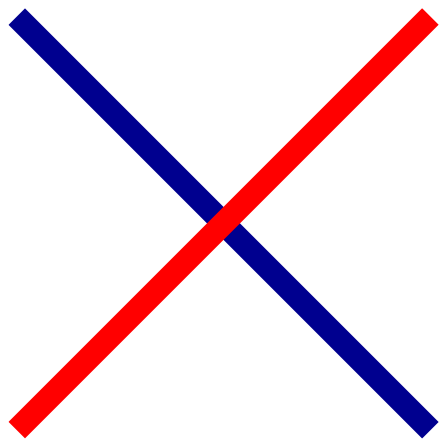}\vspace{1.5ex}
\end{equation*}
\item and the 6-valent vertex, with adjacent colors $i$ and $j$
\begin{equation*}
\labellist
\tiny\hair 2pt
\pinlabel $i$   at  -4 -10 
\pinlabel $j$ at  66 -12
\pinlabel $j$ at 232 -12
\pinlabel $i$   at 298 -10
\endlabellist
\figins{-17}{0.55}{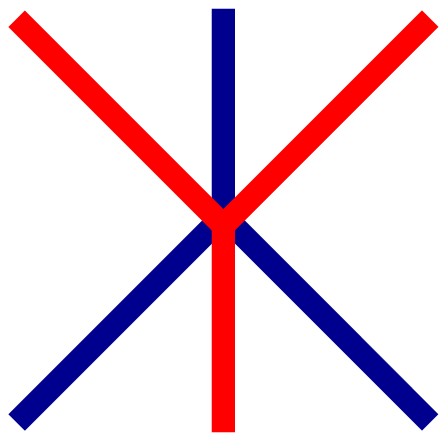}
\mspace{55mu}
\figins{-17}{0.55}{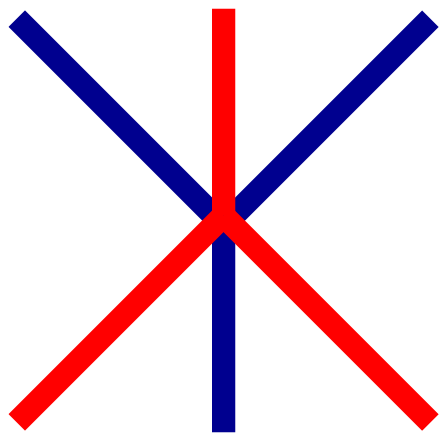}\ .
\vspace{1.5ex}
\end{equation*}
\end{itemize}
\end{itemize}

\n read from bottom to top. In this setting a diagram represents a morphism from the 
bottom bounda\-ry to the top.
We can add a new colored point to a sequence and this endows $\mathcal{SC}_1$ with a 
monoidal structure on objects, which is extended to morphisms in the obvious way. 
Composition of morphisms consists of stacking one diagram on top of the other.

We consider our diagrams modulo the following relations.

%%%%%%%%%%%%%%%%%%%%%%%%%%%
\n\emph{''Isotopy'' relations:}
\begin{equation}\label{eq:adj}
\figins{-17}{0.55}{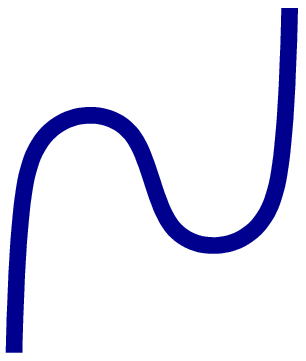}\
=\
\figins{-17}{0.55}{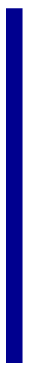}\
=\
\figins{-17}{0.55}{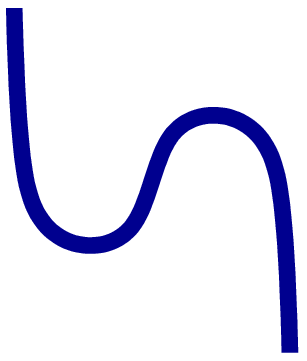}
\end{equation}

\begin{equation}\label{eq:curldot}
\figins{-17}{0.55}{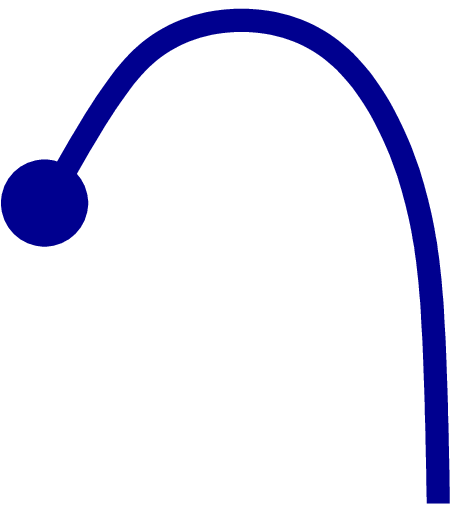}\
=\
\figins{-17}{0.55}{enddot}\
=\
\figins{-17}{0.55}{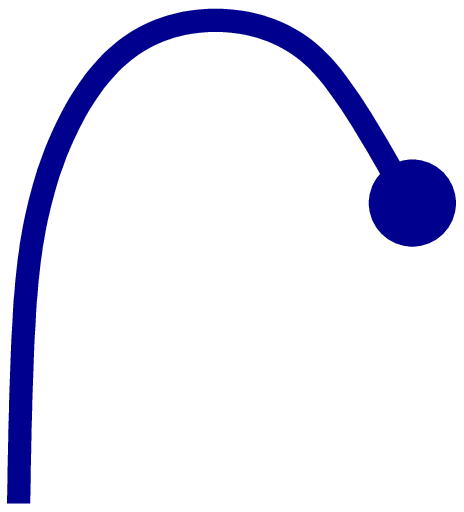}
\end{equation}

\begin{equation}\label{eq:v3rot}
\figins{-17}{0.55}{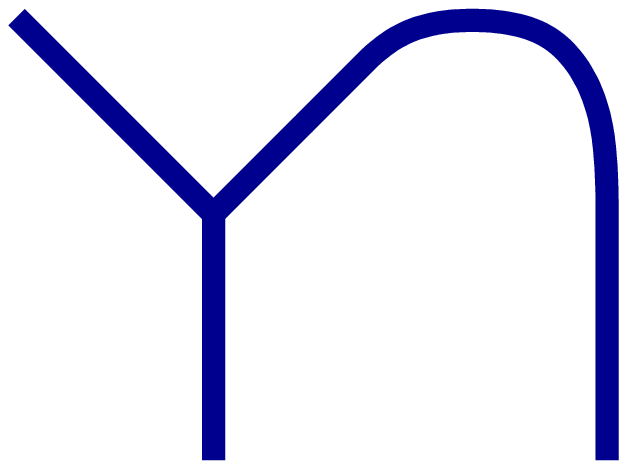}\
=\
\figins{-17}{0.55}{merge}\
=\
\figins{-17}{0.55}{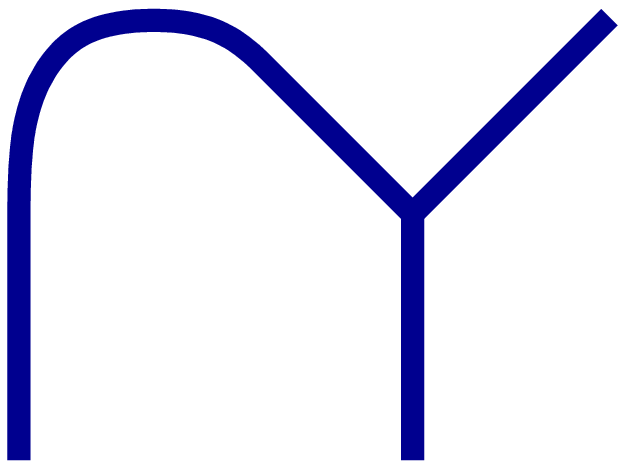}
\end{equation}

\begin{equation}\label{eq:v4rot}
\figins{-17}{0.55}{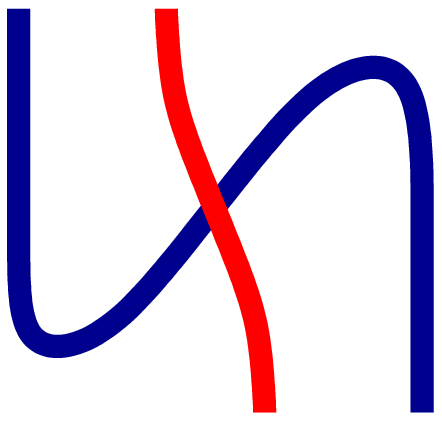}\
=\
\figins{-17}{0.55}{4vert}\
=\
\figins{-17}{0.55}{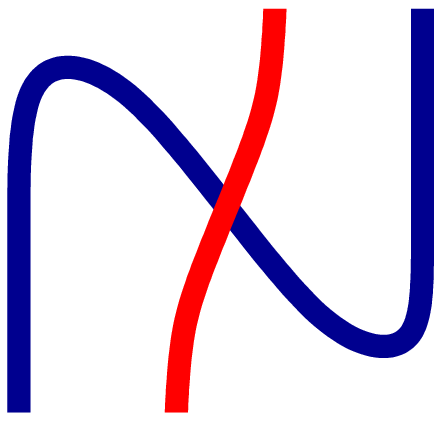}
\end{equation}

\begin{equation}\label{eq:v6rot}
\figins{-17}{0.55}{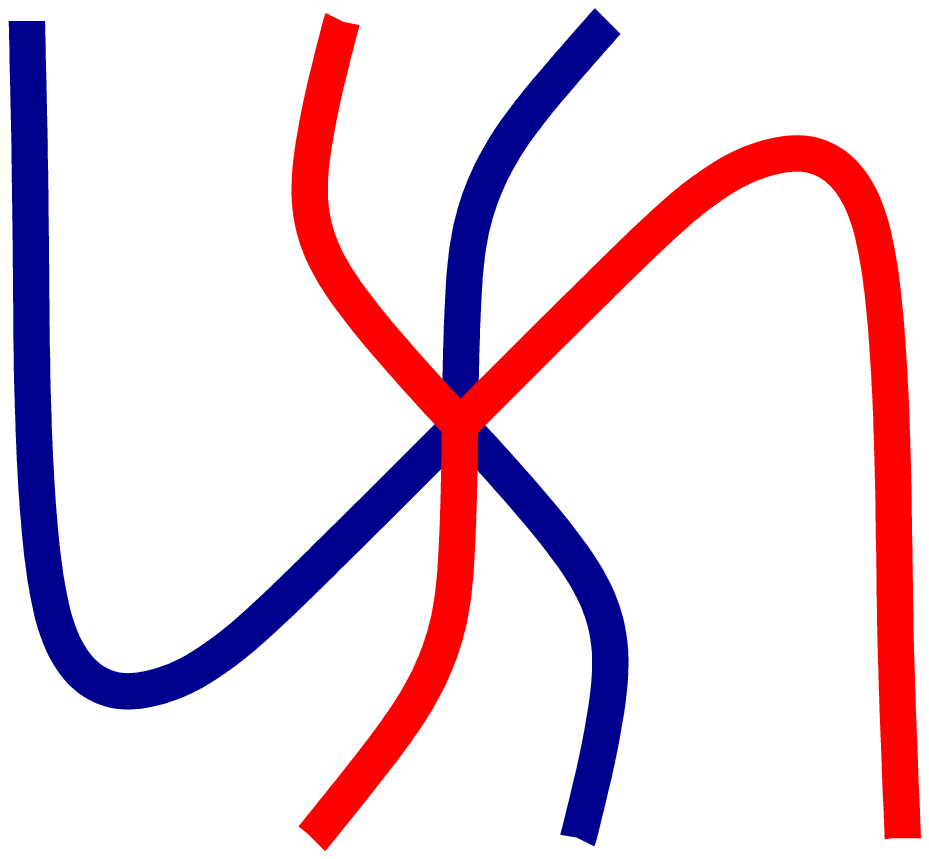}\
=\
\figins{-17}{0.55}{6vertu}\
=\
\figins{-17}{0.55}{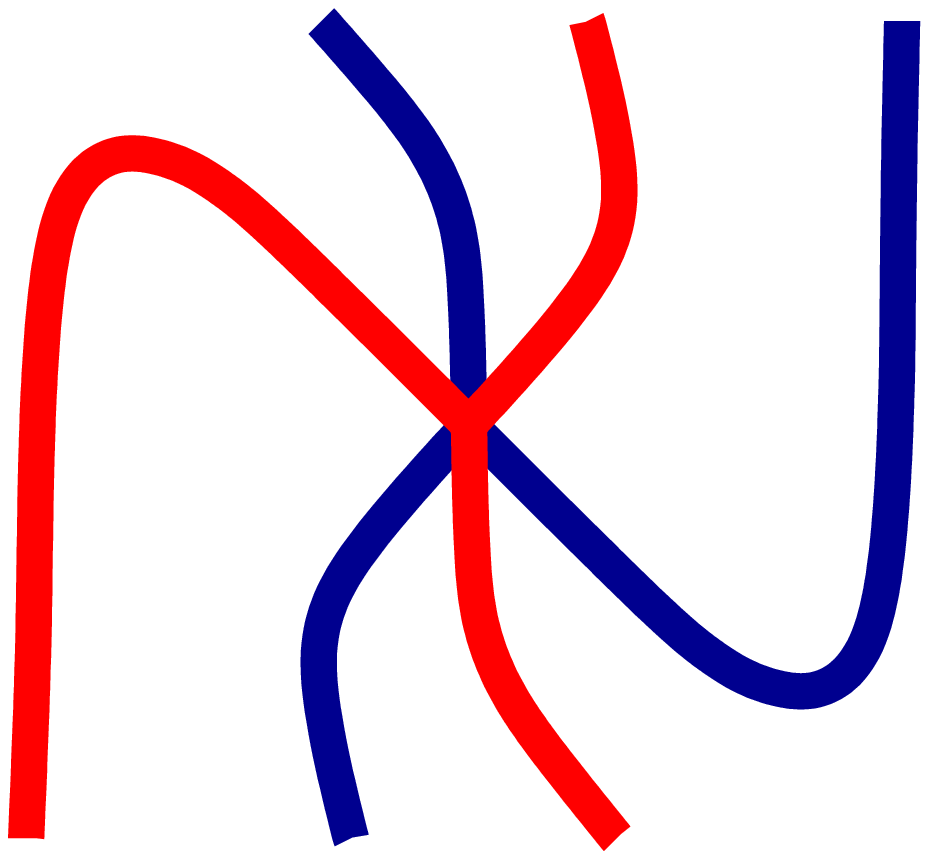}
\end{equation}

The relations are presented in terms of diagrams with generic colorings. 
Because of isotopy invariance, one may draw a diagram with a boundary on the side,
and view it as a morphism in $\mathcal{SC}_1$ by either bending the line up or down.
By the same reasoning, a horizontal line corresponds to a sequence of cups and caps.

%%%%%%%%%%%%%%%%%%%%%%%%%
\medskip
\n\emph{One color relations:}

\begin{equation}\label{eq:dumbrot}
\figins{-16}{0.5}{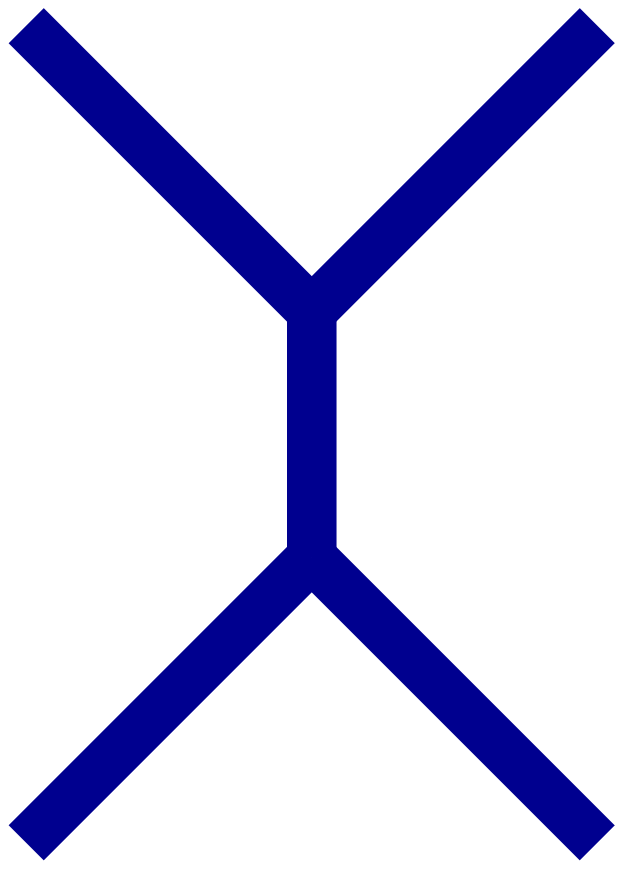}\
=\
\figins{-14}{0.45}{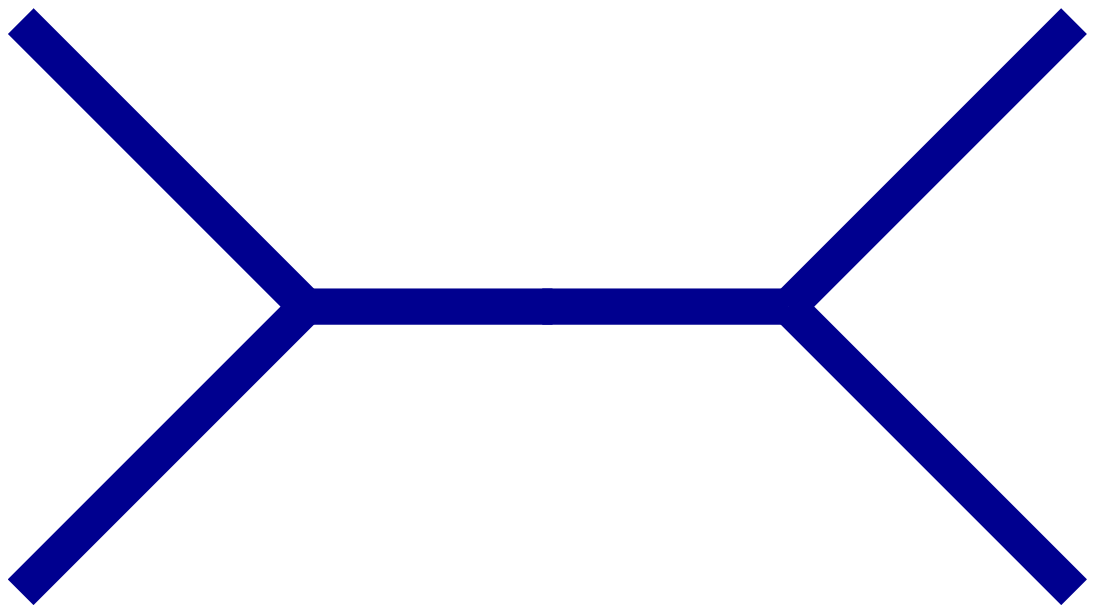}
\end{equation}

\begin{equation}\label{eq:lollipop}
\figins{-17}{0.55}{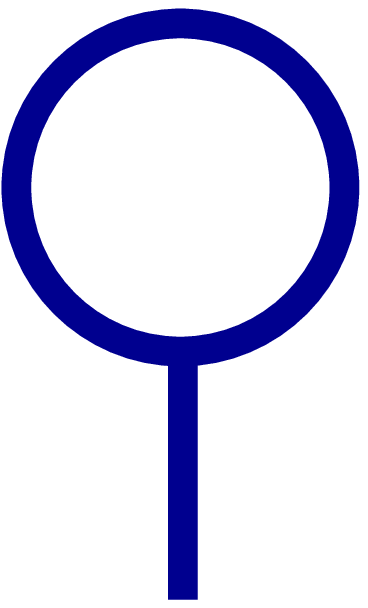}\
=\
0
\end{equation}

\begin{equation}\label{eq:deltam}
\figins{-17}{0.55}{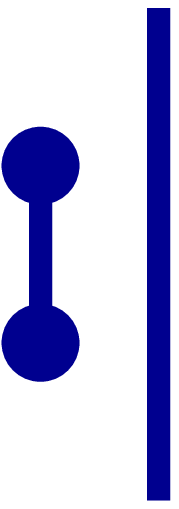}\
+\
\figins{-17}{0.55}{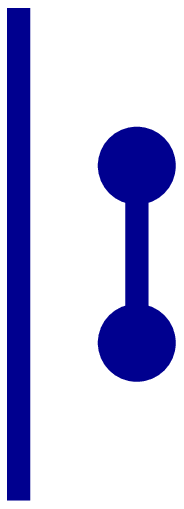}\
=\ 2\ 
\figins{-17}{0.55}{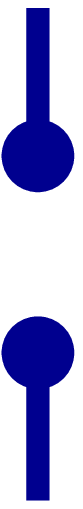}
\end{equation}

%%%%%%%%%%%%%%%%%%%%%%%%
\medskip
\n\emph{Two distant colors:}
\begin{equation}\label{eq:reid2dist}
\figins{-32}{0.9}{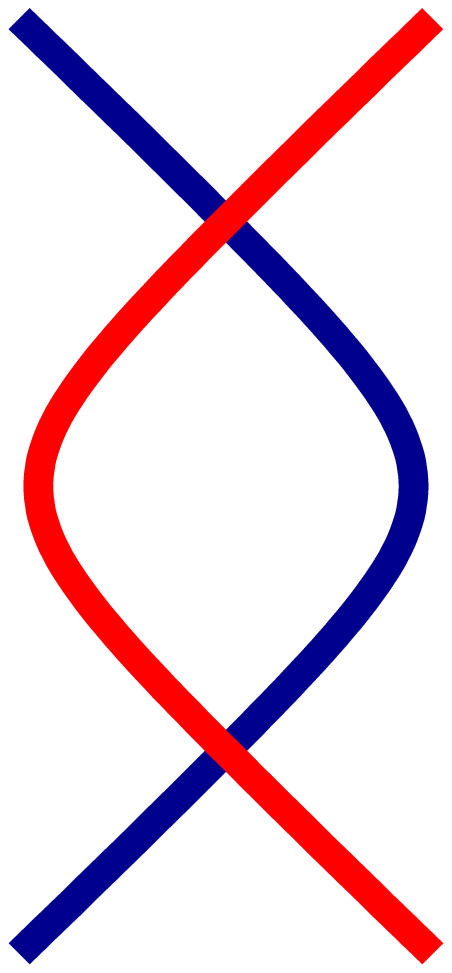}\
=\
\figins{-32}{0.9}{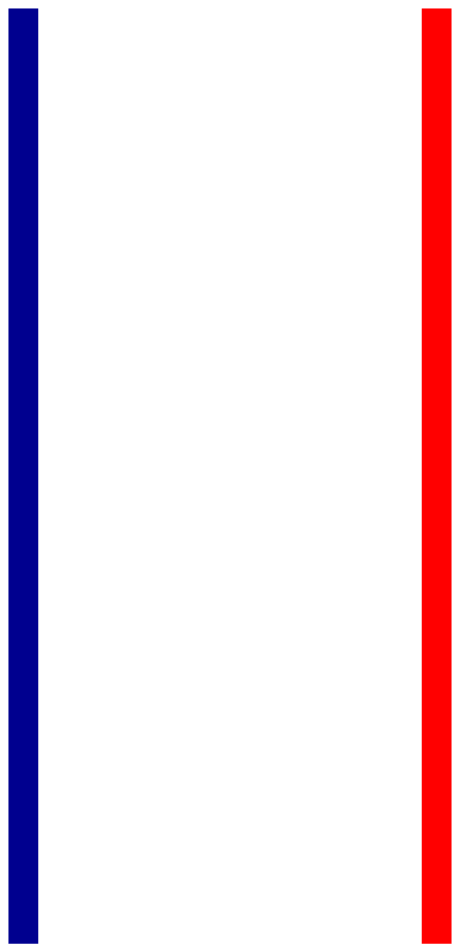}
\end{equation}

\begin{equation}\label{eq:slidedotdist}
\figins{-16}{0.5}{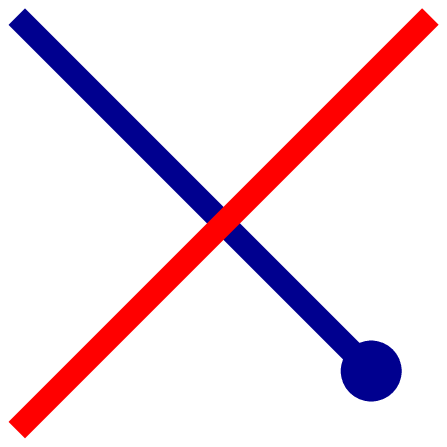}\
=\
\figins{-16}{0.5}{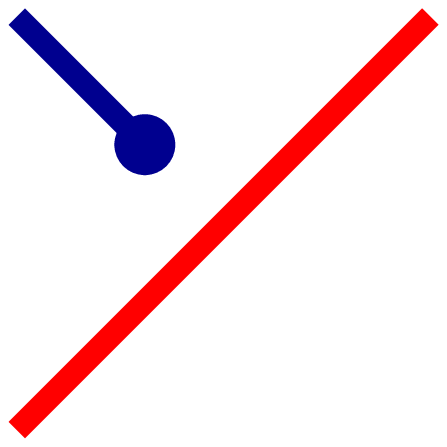}
\end{equation}

\begin{equation}\label{eq:slide3v}
\figins{-17}{0.55}{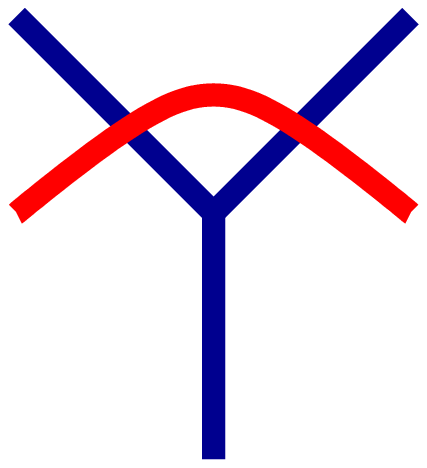}\
=\
\figins{-17}{0.55}{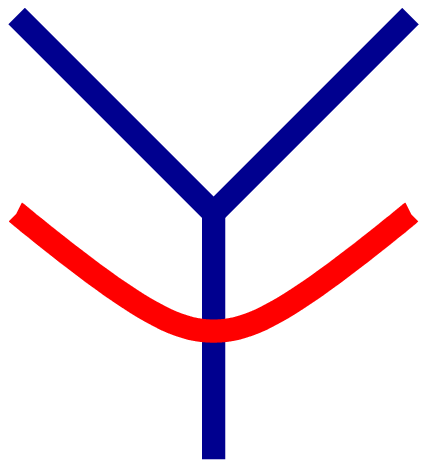}
\end{equation}

%%%%%%%%%%%%%%%%%%%%%
\medskip
\n\emph{Two adjacent colors:}
\begin{equation}\label{eq:dot6v}
\figins{-16}{0.5}{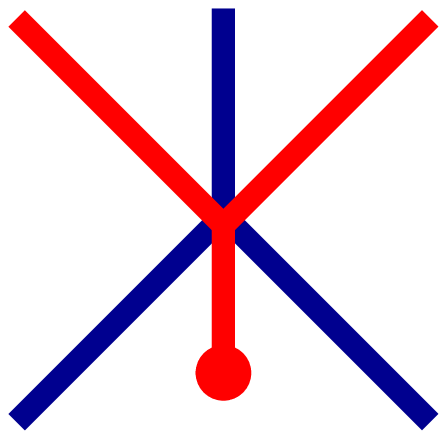}\
=\
\figins{-16}{0.5}{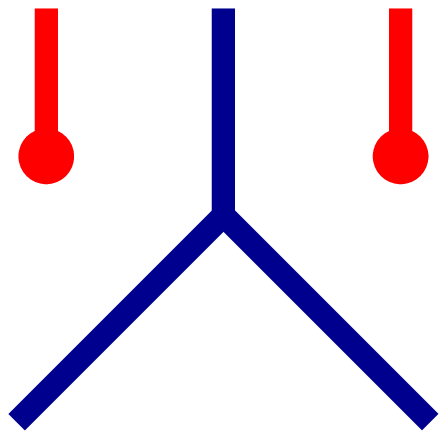}\
+\
\figins{-16}{0.5}{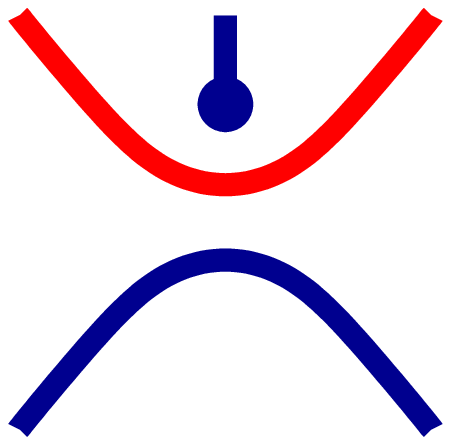}
\end{equation}

\begin{equation}\label{eq:reid3}
\figins{-30}{0.85}{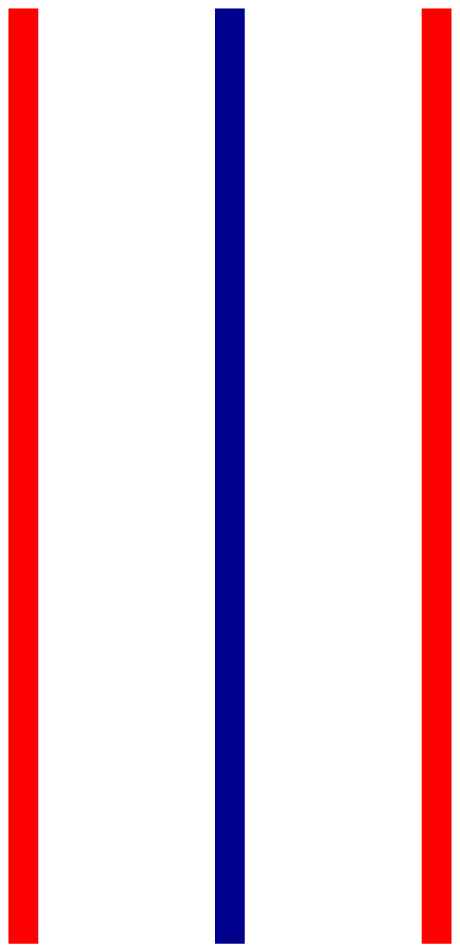}\
=\
\figins{-30}{0.85}{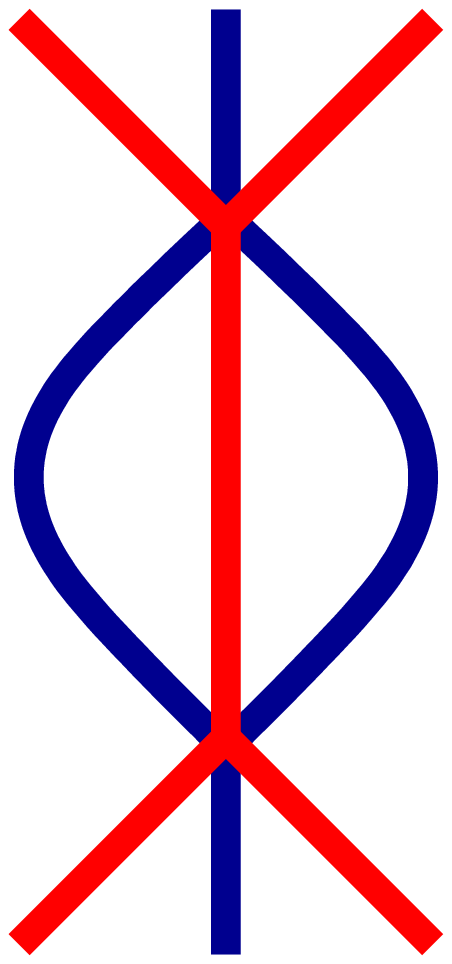}\
-\
\figins{-30}{0.85}{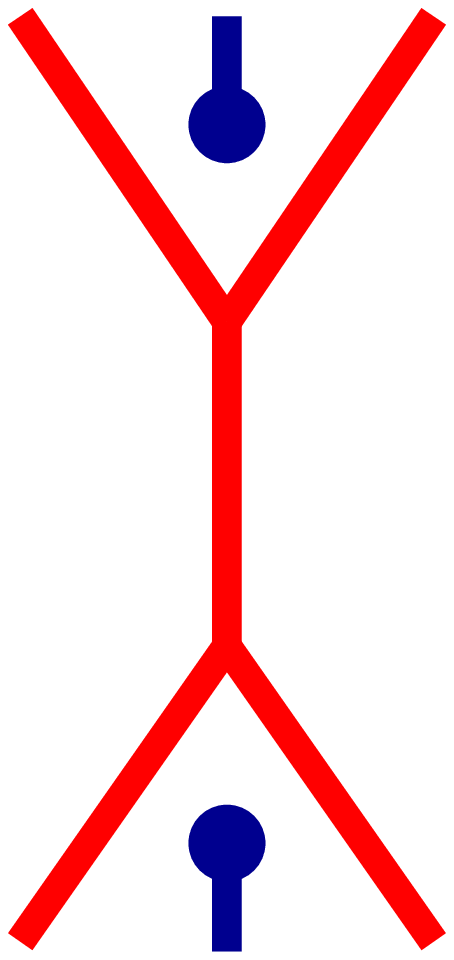}
\end{equation}

\begin{equation}\label{eq:dumbsq}
\figins{-30}{0.85}{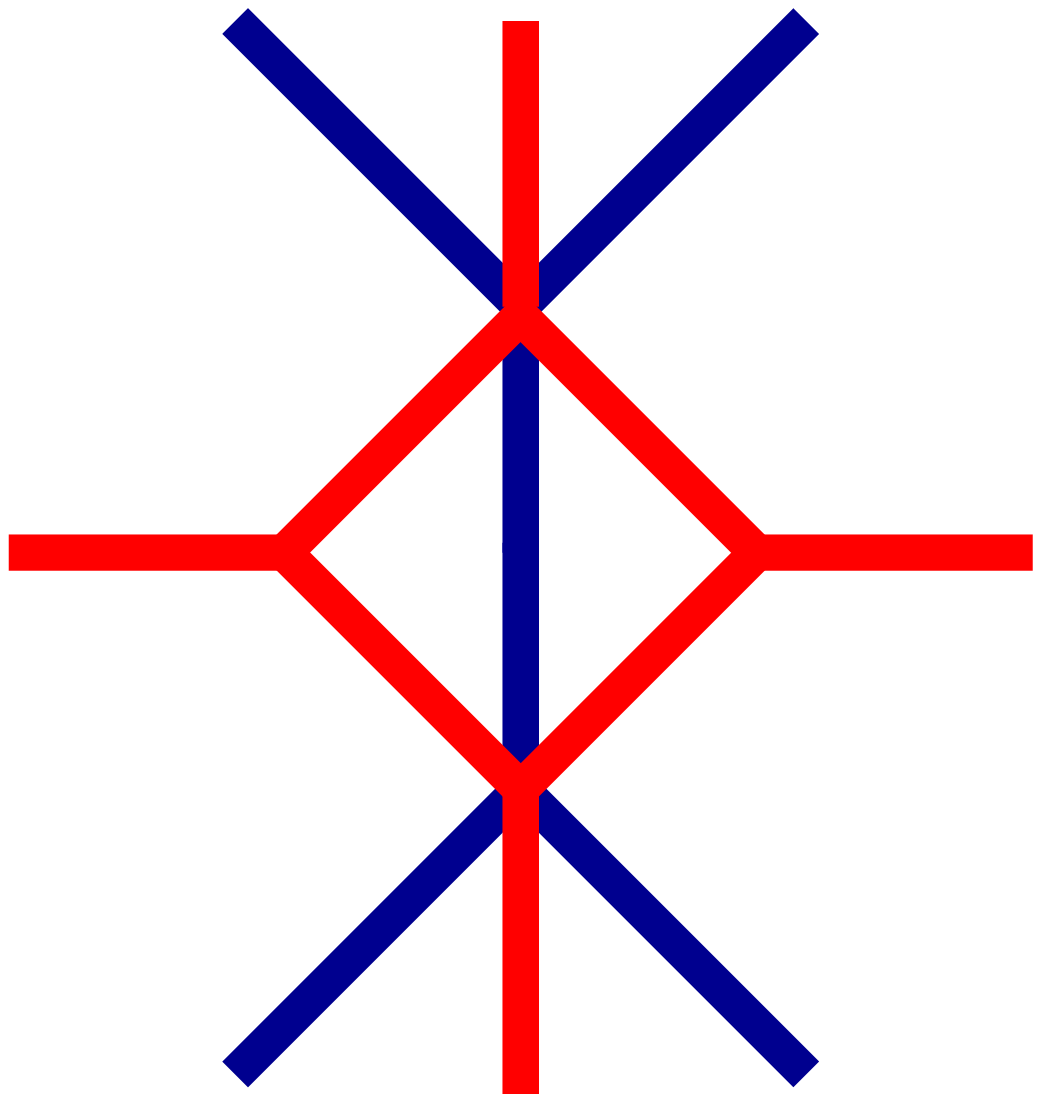}\
=\
\figins{-30}{0.85}{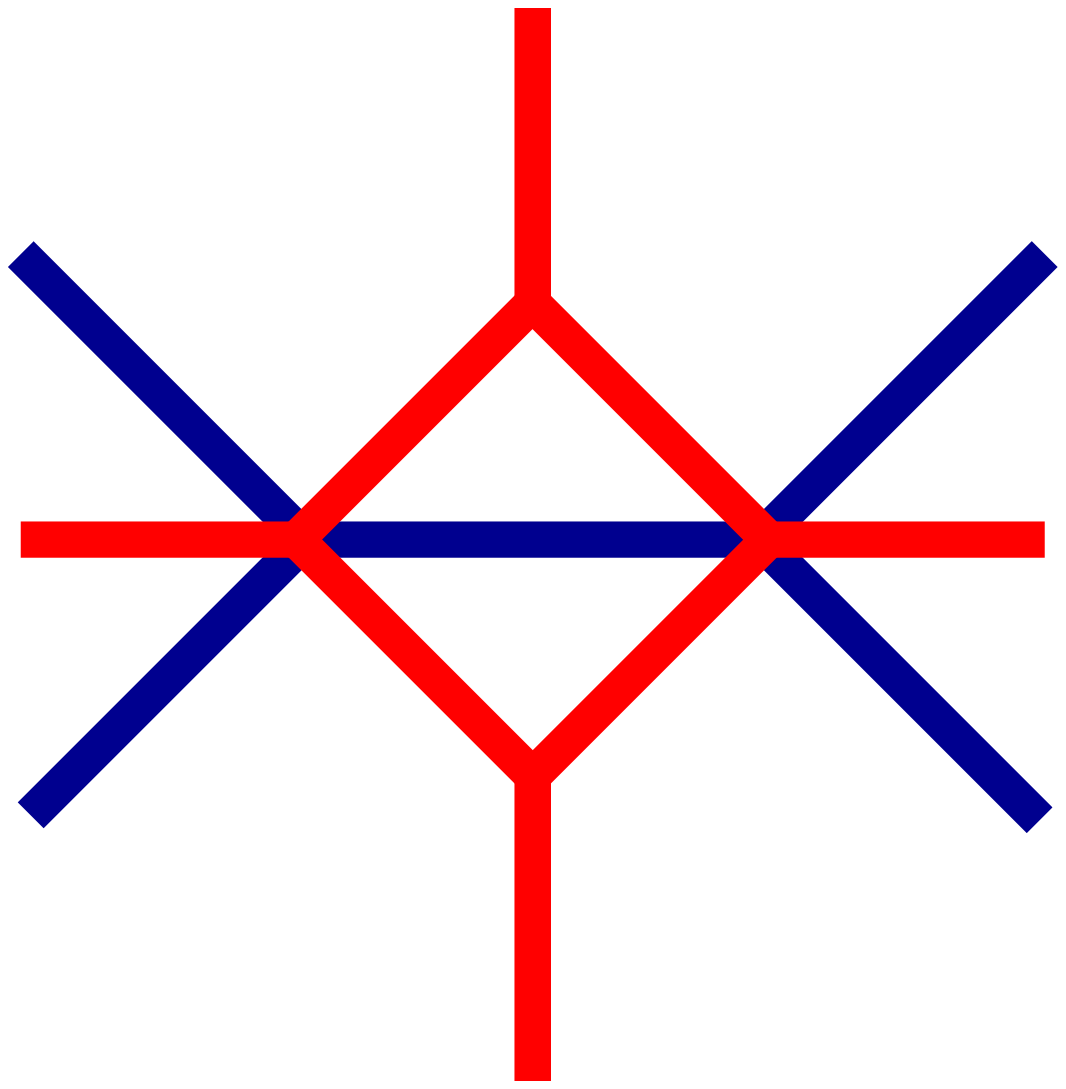}
\end{equation}

\begin{equation}\label{eq:slidenext}
\figins{-17}{0.55}{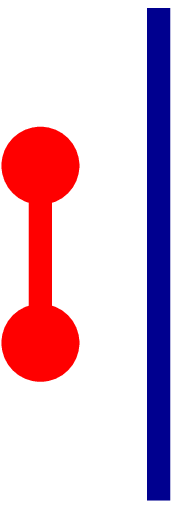}\
-\
\figins{-17}{0.55}{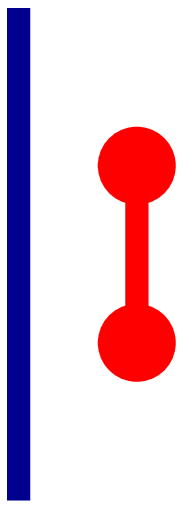}\
=\
\frac{1}{2}
\Biggl(\
\figins{-17}{0.55}{edge-startenddot}\
-\
\figins{-17}{0.55}{startenddot-edge}\
\Biggr)
\end{equation}

%%%%%%%%%%%%%%%%%%%%%%%
\medskip
\n\emph{Relations involving three colors:}
(adjacency is determined by the vertices which appear)
\begin{equation}\label{eq:slide6v}
\figins{-18}{0.6}{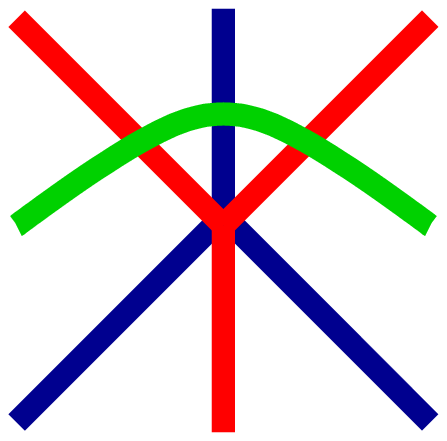}\
=\
\figins{-18}{0.6}{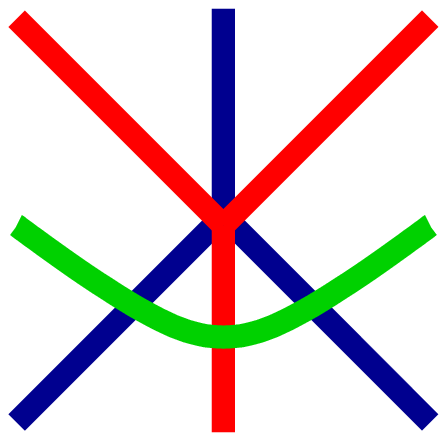}
\end{equation}

\begin{equation}\label{eq:slide4v}
\figins{-18}{0.6}{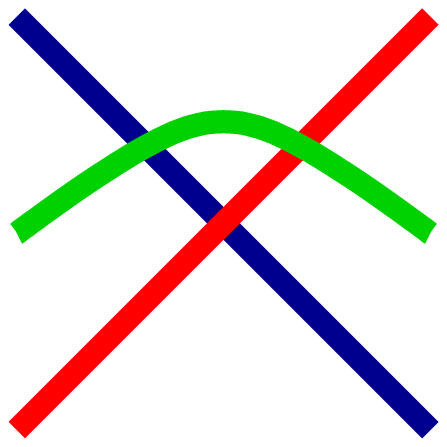}\
=\
\figins{-18}{0.6}{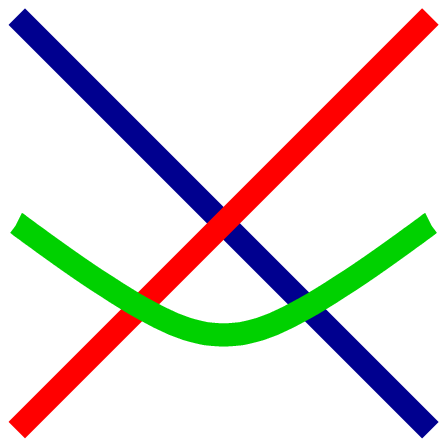}
\end{equation}

\begin{equation}\label{eq:dumbdumbsquare}
\figins{-30}{0.85}{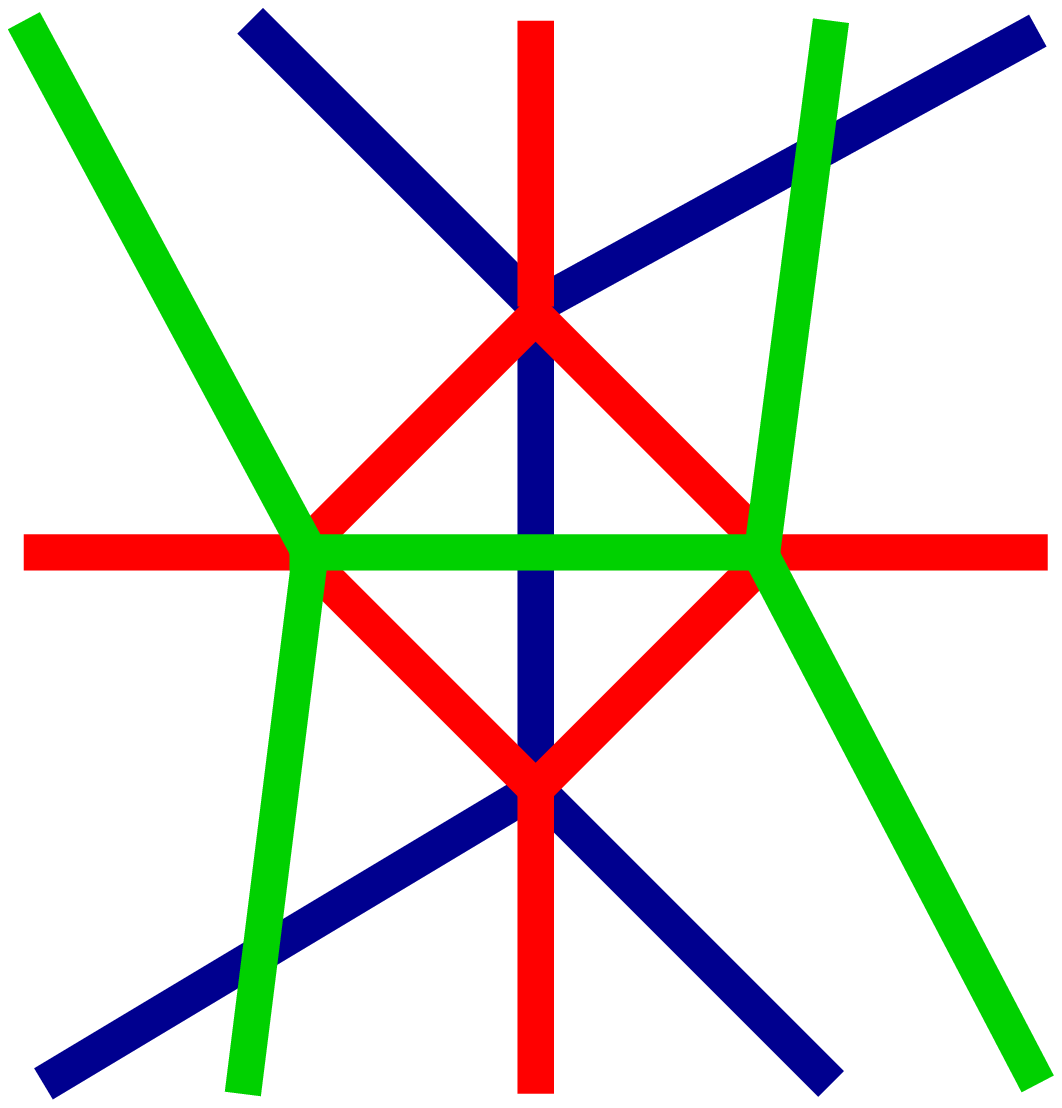}\
=\
\figins{-30}{0.85}{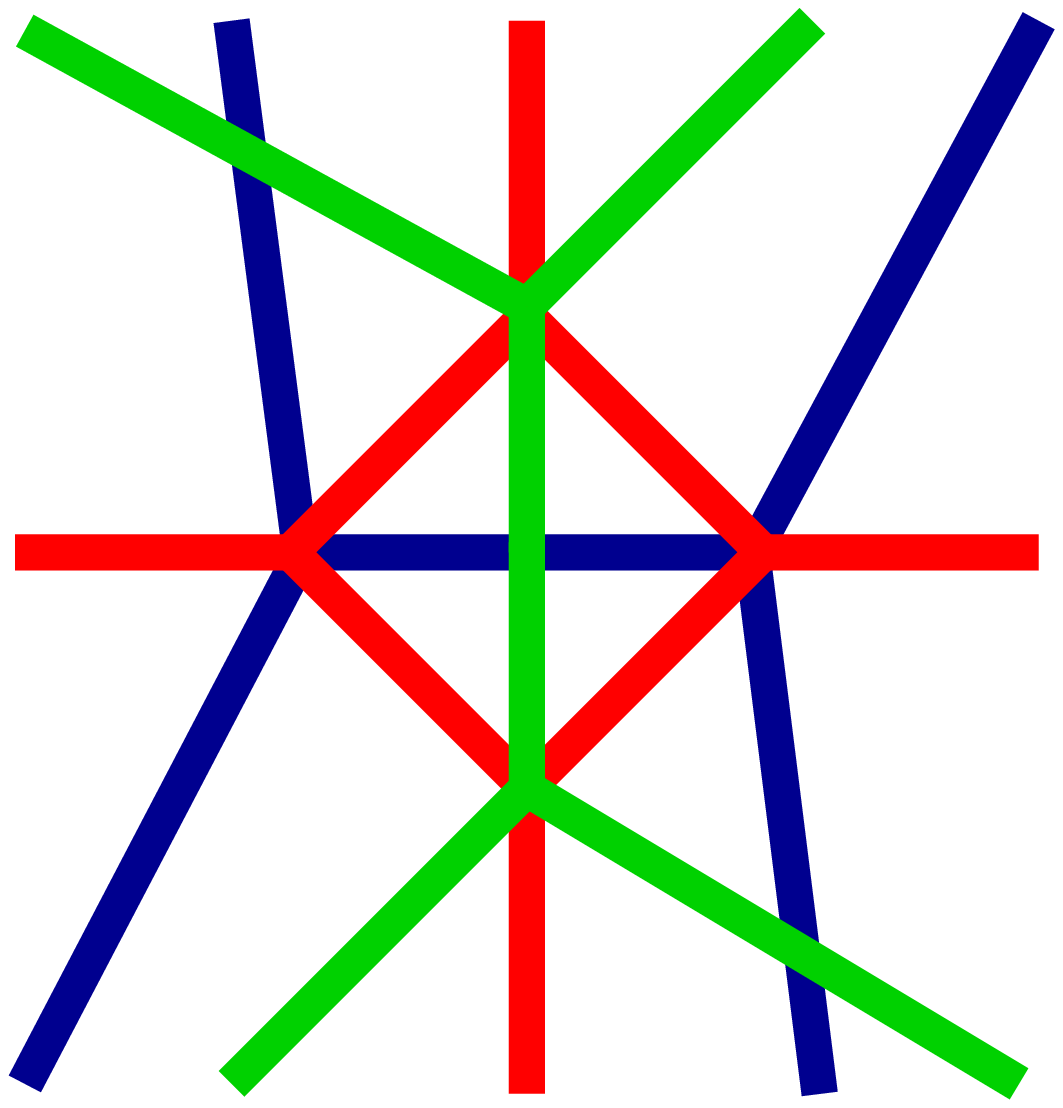}
\end{equation}

Introduce a $q$-grading on $\mathcal{SC}_1$ declaring that dots have degree $1$, 
trivalent vertices have degree $-1$ 
and $4$- and $6$-valent vertices have degree $0$.

\begin{defn}
The category $\mathcal{SC}_2$ is the category containing all direct sums and grading 
shifts of objects in 
$\mathcal{SC}_1$ and whose morphisms are the grading preserving morphisms from 
$\mathcal{SC}_1$.
\end{defn}

\begin{defn}
The category $\mathcal{SC}$ is the Karoubi envelope of the category $\mathcal{SC}_2$.
\end{defn}

Elias and Khovanov's main result in~\cite{EKh} is the following theorem.
\begin{thm}[Elias-Khovanov]
The category $\mathcal{SC}$ is equivalent to the Soergel category in~\cite{Soergel}.
\end{thm}

From Soergel's results from~\cite{Soergel} we have the following corollary.
\begin{cor}
The Grothendieck algebra of $\mathcal{SC}$ is isomorphic to the Hecke algebra.
\end{cor}

Notice that $\mathcal{SC}$ is an additive category but not abelian and we are using the 
(additive) split Grothendieck algebra.

In Sections~\ref{sec:sl2} and~\ref{sec:sl3} we will define functors from $\mathcal{SC}_1$ 
to the categories of 
$sl(2)$ and $sl(3)$ foams. 
These functors are grading preserving, so they obviously extend uniquely to $\mathcal{SC}_2$. 
By the universality of the Karoubi envelope, they also extend uniquely to functors between 
the respective Karoubi envelopes.

%%%%%%%%%%%%%%%%%%%%%%%%%%%%%%%%%%%%%%%%
%%%                                  %%%
%%%            sl(2)                 %%%
%%%                                  %%%
%%%%%%%%%%%%%%%%%%%%%%%%%%%%%%%%%%%%%%%%
\section{The $sl(2)$ case}
\label{sec:sl2}

\subsection{Clark-Morrison-Walker's category of disoriented $sl(2)$ foams}

In this subsection we review the category $\cmw$ of $sl(2)$ foams following Clark, Morrison 
and Walker's construction 
in~\cite{CMW}. This category was introduced in~\cite{CMW} to modify Khovanov's link 
homology theory making it properly functorial with respect to link cobordisms. 
Actually we will use the version with dots of Clark-Morrison-Walker's 
original construction in~\cite{CMW}. 
Recall that we obtain one from the other by replacing each dot by 
$1/2$ times a handle.

A \emph{disoriented arc} is an arc composed by oriented segments with oppositely oriented 
segments separated by a mark pointing to one of these segments.
A \emph{disoriented diagram} consists of a collection $D$ of disoriented arcs in the 
strip in $\bR^2$ bounded by the lines $y=0,1$ containing the boundary points of $D$. 
We allow diagrams containing oriented and disoriented circles.
Disoriented diagrams can be composed vertically which endows $\cmw$ with a monoidal 
structure on objects.
For example, the diagrams $1_n$ and $u_j$ for ($1 < j < n$) are disoriented diagrams: 
\begin{equation*}
\labellist
\pinlabel $\dotsm$ at 125 71
\pinlabel $\dotsm$ at 565 71
\pinlabel $\dotsm$ at 851 71
\tiny\hair 2pt
\pinlabel $1$   at   8 -15
\pinlabel $2$   at  64 -15
\pinlabel $n$   at 178 -16
\pinlabel $1$   at 507 -15
\pinlabel $j$   at 674 -15
\pinlabel $j+1$ at 736 -16
\pinlabel $n$   at 908 -15
\endlabellist
1_n\ =\
\figins{-18}{0.6}{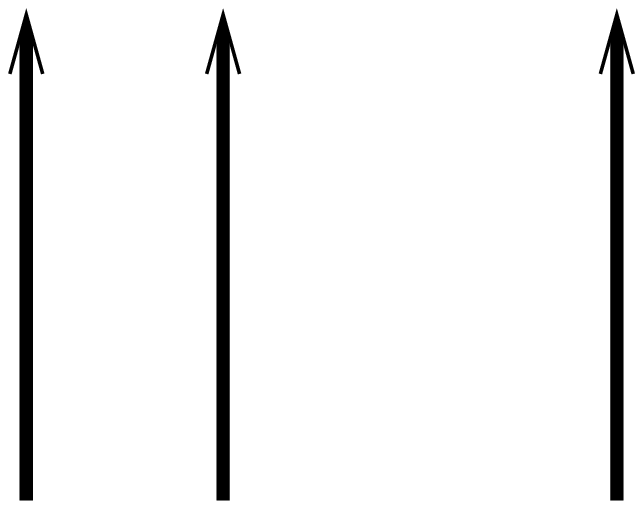}
\mspace{100mu}
u_j\ =\
\figins{-18}{0.6}{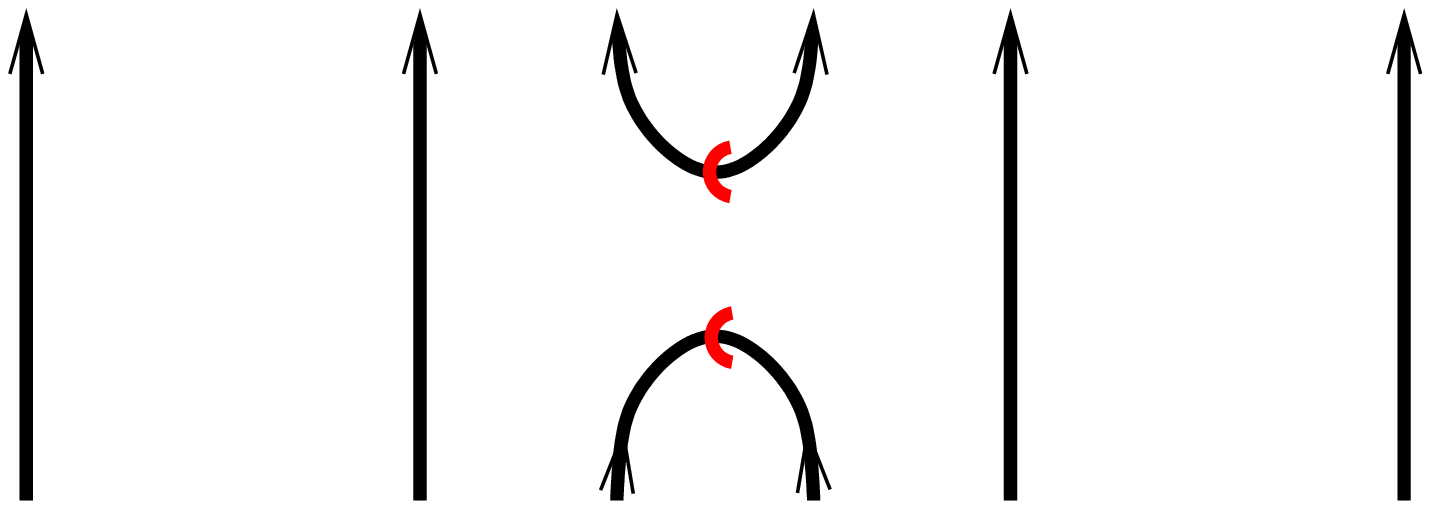}
\vspace*{2ex}
\end{equation*}

A \emph{disoriented cobordism} between disoriented diagrams is a $2d$ cobordism which 
can be deco\-rated with 
dots and with 
seams separating differently oriented regions and such that the vertical boundary of 
each cobordism is a set 
(possibly empty) 
of vertical lines. 
Disorientation seams can have one out of two possible orientations which we identify with 
a fringe. 
We read cobordisms from bottom to top. For example, 
\begin{equation*}
\labellist
\pinlabel $\dotsm$ at 180 101
\pinlabel $\dotsm$ at 450 101
\tiny\hair 2pt
\pinlabel $1$   at  -5 -15
\pinlabel $j$   at 204 -16
\pinlabel $n$   at 480 -15
\endlabellist
\figins{-32}{1.0}{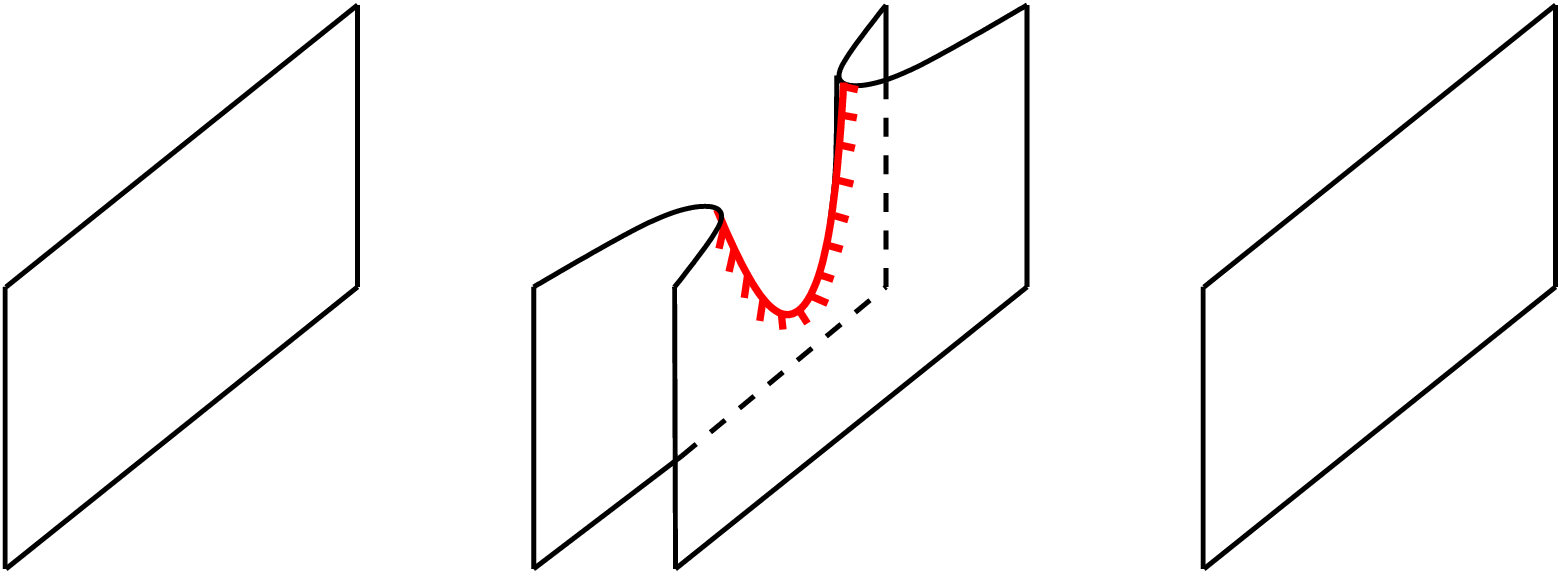}
\vspace*{2ex}
\end{equation*}
is a disoriented cobordism from $1_n$ to $u_j$.

\medskip

Cobordism composition consists of placing one cobordism on top of the other and the 
monoidal structure is given by vertical composition which corresponds to placing one 
cobordism behind the other in our pictures. 
Let $\bC[t]$ be the ring of polynomials in $t$ with coefficients in $\bC$. 
\begin{defn}
The category $\cmw$ is the category whose objects are disoriented diagrams, 
and whose morphisms are $\bC[t]$-linear combinations of isotopy classes of 
disoriented cobordisms, modulo some relations,

\indent$\bullet$ the disorientation relations
\begin{gather}
\label{eq:disorc}
\figins{-12}{0.4}{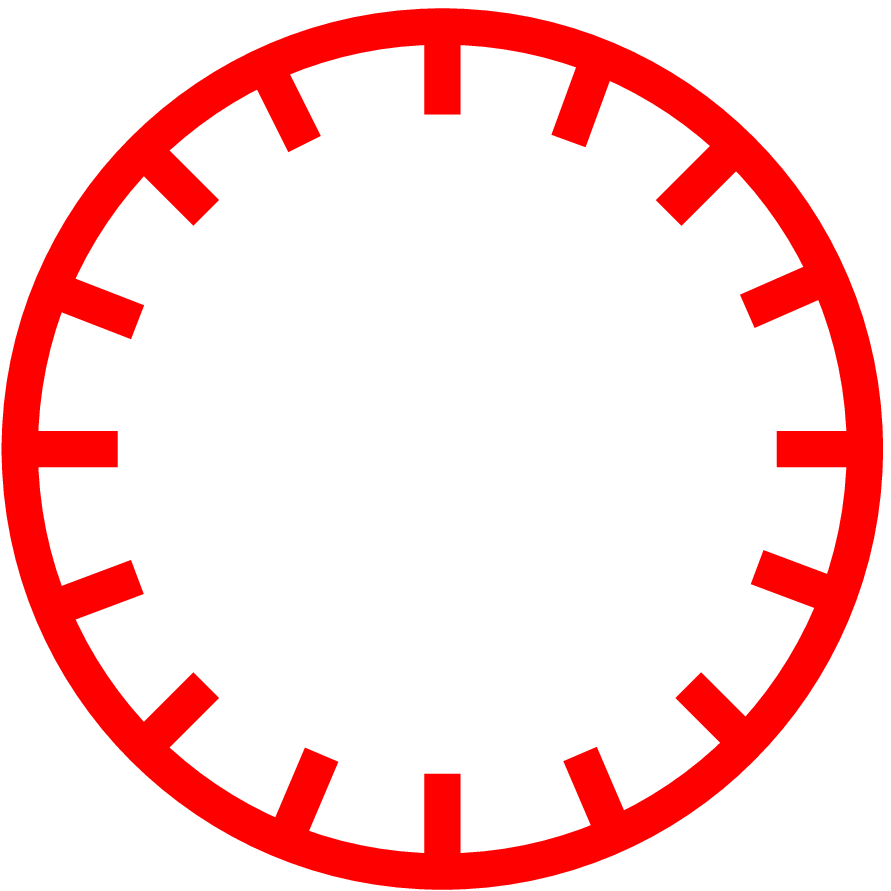}\
=\ i
\mspace{80mu}
\figins{-15}{0.48}{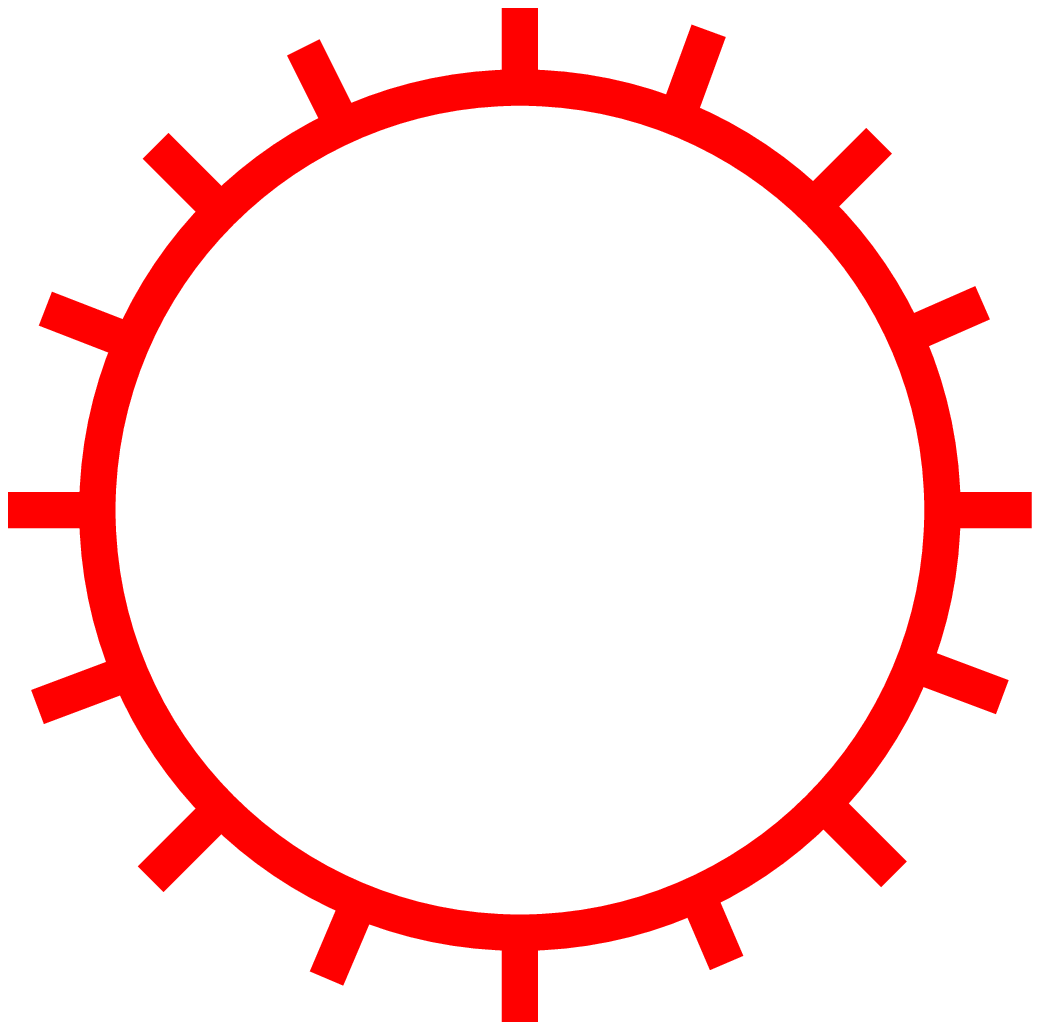}\
=\ -i
\\[1.5ex]\displaybreak[0]
\label{eq:disors}
\figins{-12}{0.4}{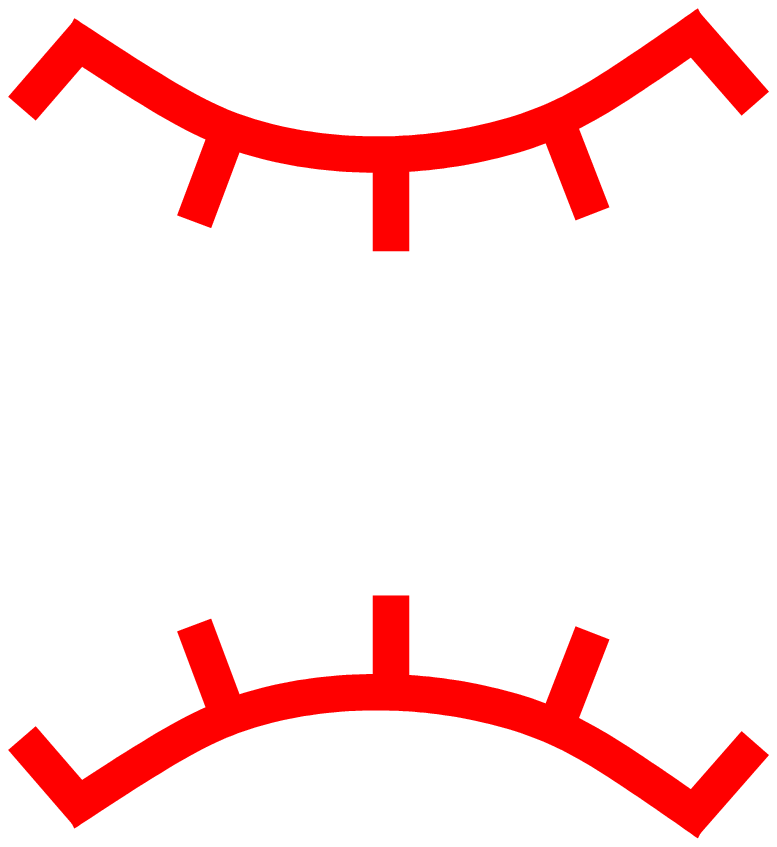}\
=\ -i\
\figins{-12}{0.4}{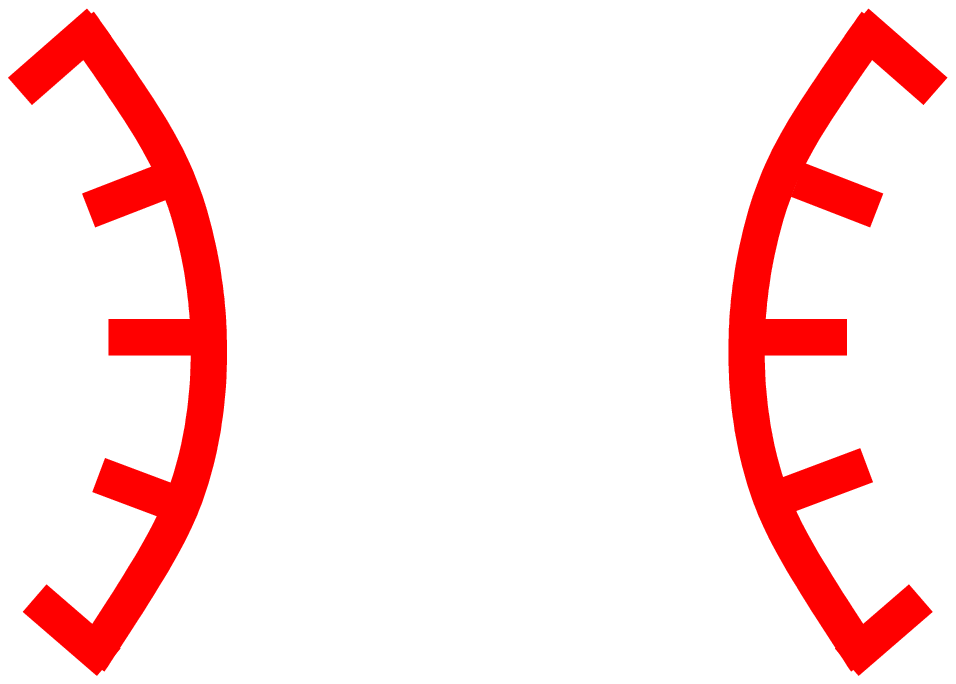}
\\[1.5ex]\displaybreak[0]
\label{eq:disorp}
\figins{-8}{0.3}{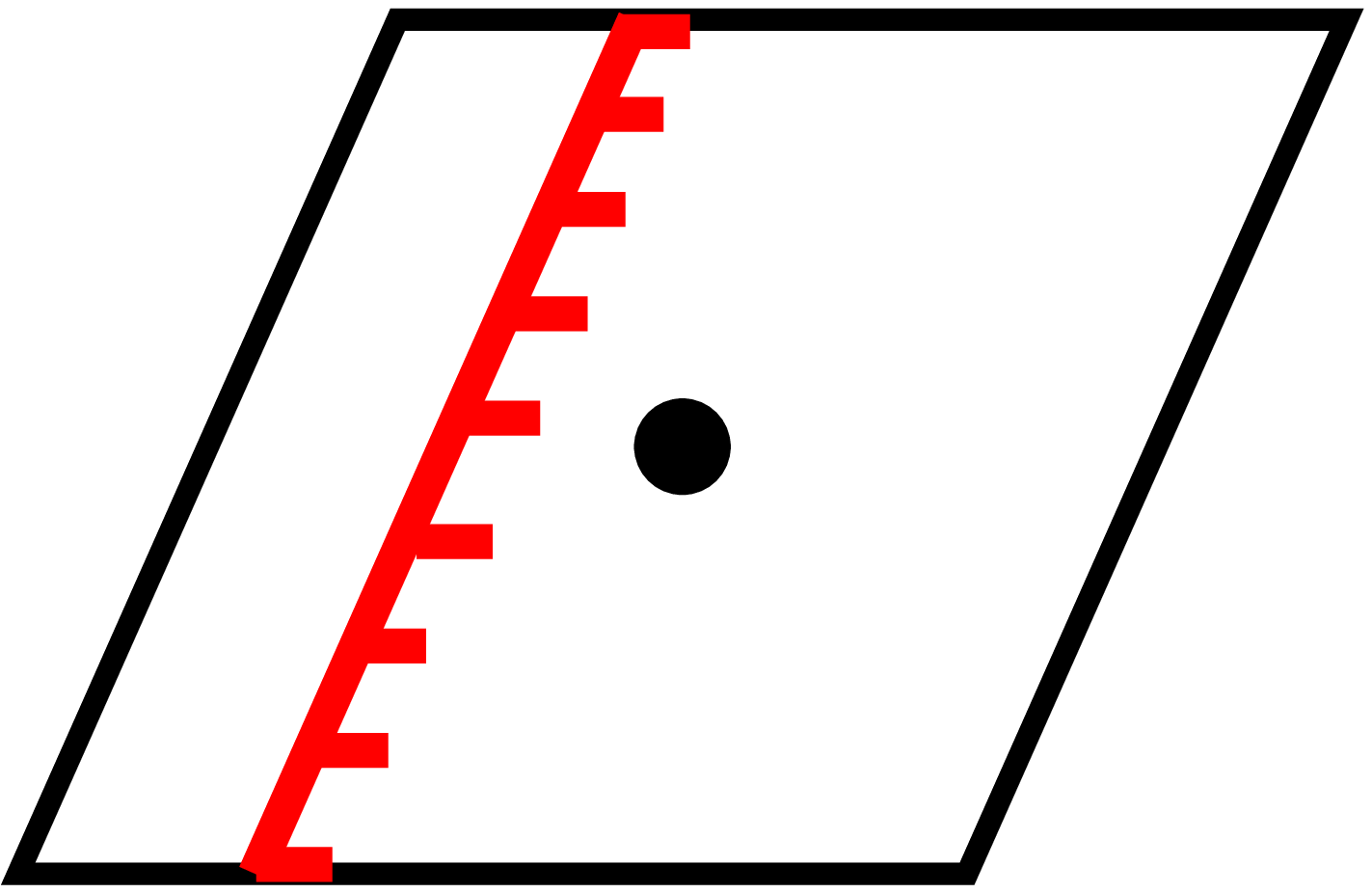}
=\ -\
\figins{-8}{0.3}{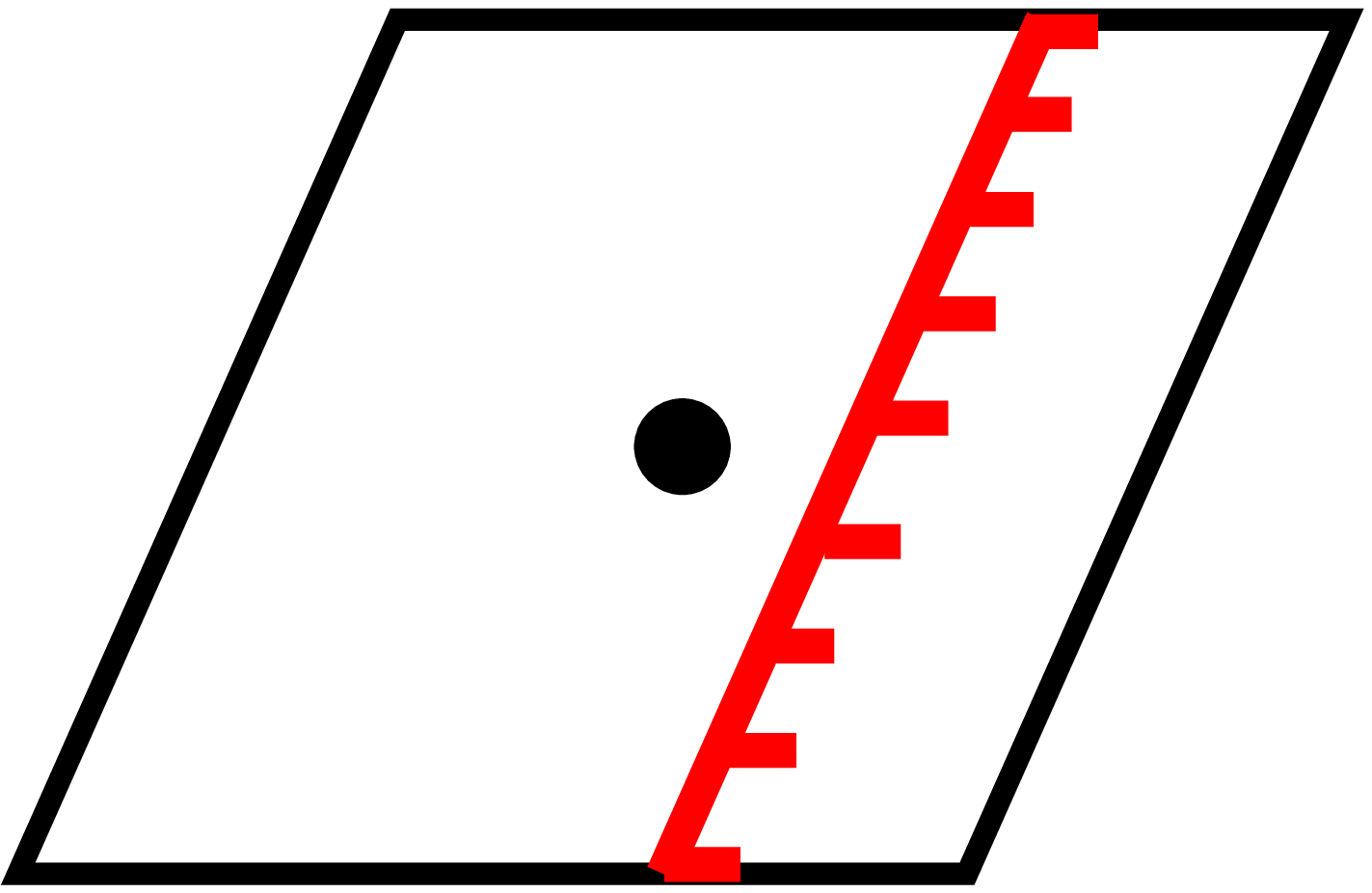}
\end{gather}
where $i$ is the imaginary unit,

\medskip

\indent$\bullet$ and the Bar-Natan (BN) relations, 
\begin{gather}
\label{eq:BN-dot}
\figins{-6}{0.25}{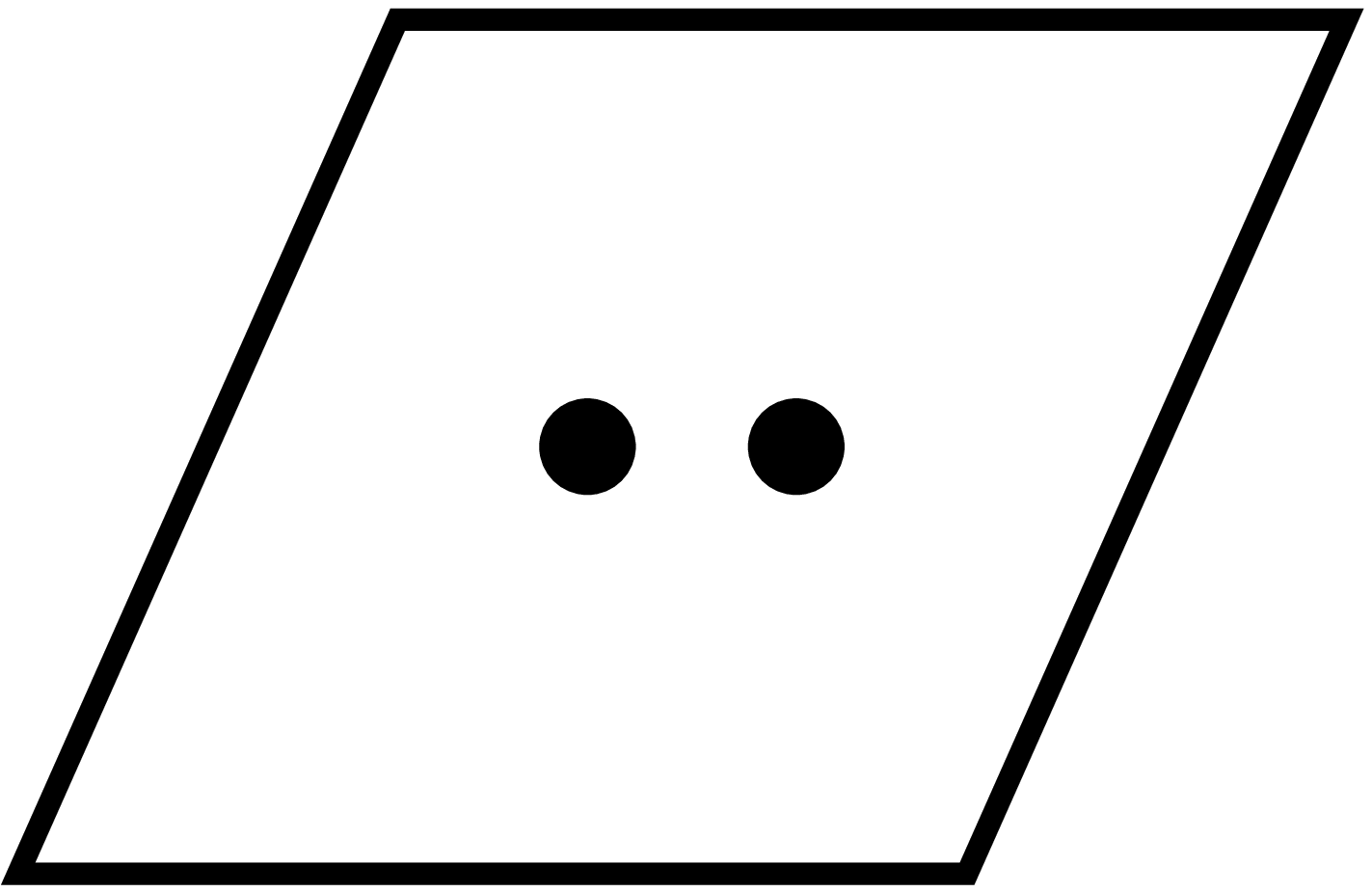}\
=\ t
\figins{-6}{0.25}{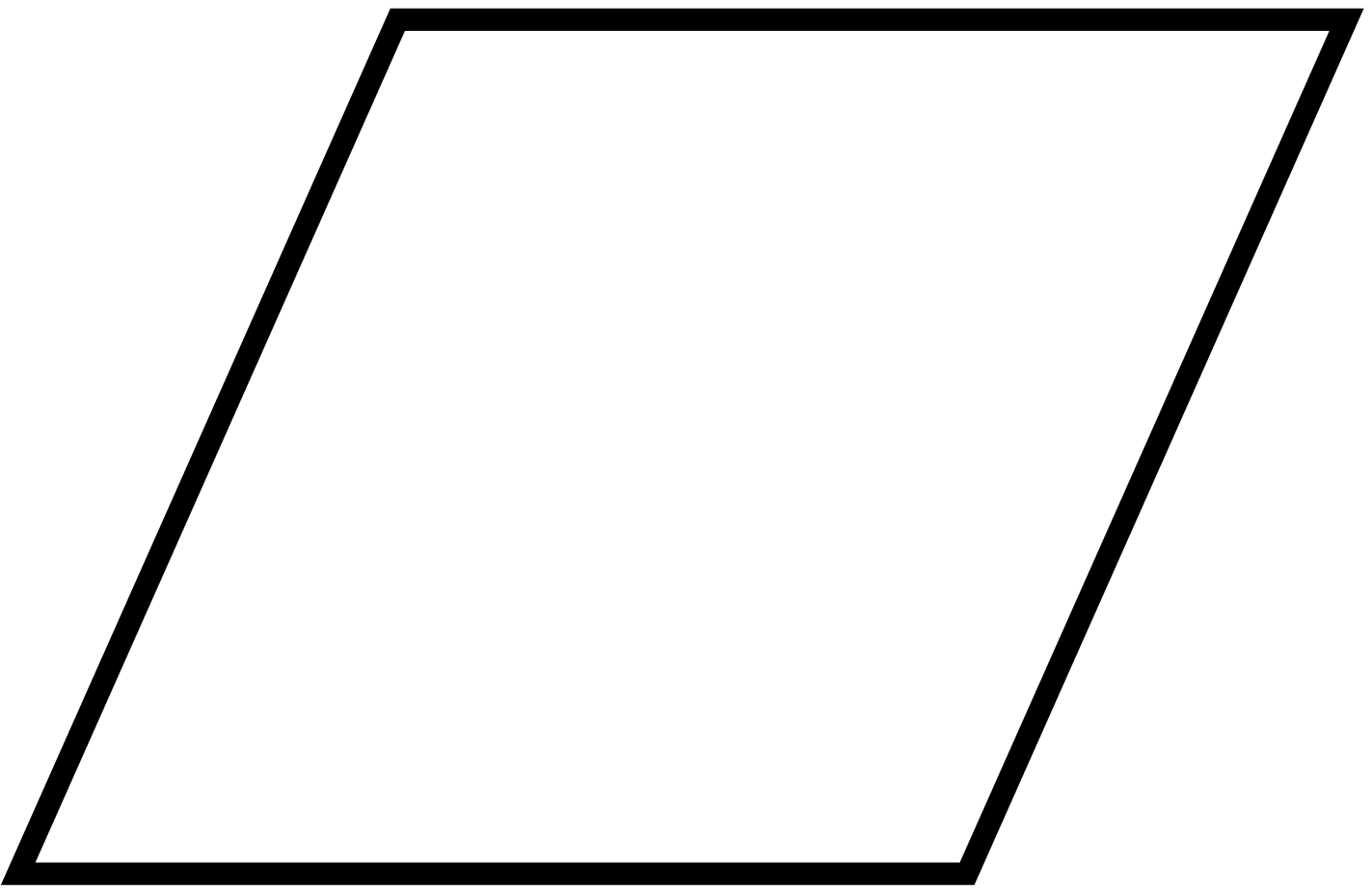} 
\\[1.5ex]\displaybreak[0]
\label{eq:BN-s}
\figins{-8}{0.3}{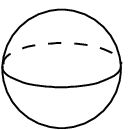}\
=\ 0
\mspace{60mu}
\figins{-8}{0.3}{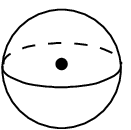}\
=\ 1
\\[1.5ex]\displaybreak[0]
\label{eq:BN-cn}
\figwhins{-22}{0.65}{0.26}{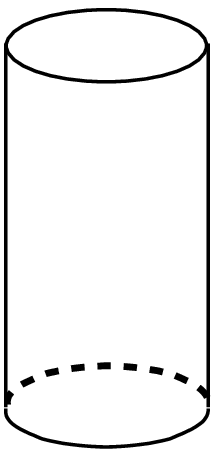}\
=\
\figwhins{-22}{0.65}{0.26}{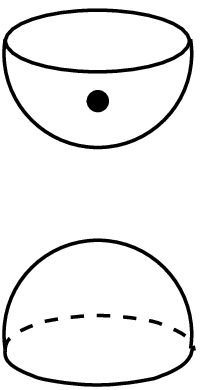}\
+\
\figwhins{-22}{0.65}{0.26}{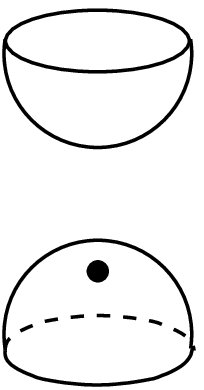}
\end{gather}
which are only valid away from the disorientations.
\end{defn}

The universal theory for the original Khovanov homology contains 
another parameter $h$, but we have to put $h=0$ in the Clark-Morrison-Walker's cobordism
theory over a field of characteristic zero. For suppose we have a cylinder with a transversal
disoriented circle. Applying~\eqref{eq:BN-cn} on one side of the disorientation 
circle followed by the disoriented relation~\eqref{eq:disorc} gives a cobordism that is 
independent of 
the side chosen to apply~\eqref{eq:BN-cn} only 
if $h=0$ over a field of characteristic zero.

Define a $q$-grading on $\bC[t]$ by $q(1)=0$ and $q(t)=4$. 
We introduce a $q$-grading on $\cmw$ as follows. Let $f$ be a cobordism with 
$\vert\bullet\vert$ dots and 
$\vert b\vert$ vertical boundary components. The $q$-grading of $f$ is given by
\begin{equation}
q(f)=-\chi (f) + 2\vert\bullet\vert + \tfrac{1}{2}\vert b\vert
\end{equation}
where $\chi$ is the Euler characteristic. 
For example the degree of a saddle is $1$ while the degree of a cap or a cup is $-1$.
The category $\cmw$ is additive and monoidal. More details about $\cmw$ can be found 
in~\cite{CMW}.

%%%%%%%%%%%%%%%%%%%%%%%%%%%%%%%%%%%%%%%%%%%%%%%%%%%%%%%%%%%%%%%
\subsection{The functor $\FKh$}

In this subsection we define a monoidal functor $\FKh$ between the categories $\mathcal{SC}$ and $\cmw$. 

\medskip
\n\emph{On objects:} 
$\FKh$ sends the empty sequence to $1_n$ and the one-term sequence $(j)$ to $u_j$
with $\FKh(jk)$ given by the vertical composite $u_{j}u_{k}$.

\medskip

\n\emph{On morphisms:}
\begin{itemize}
\item The empty diagram is sent to $n$ parallel vertical sheets:
\begin{equation*}
\labellist
\pinlabel $\dotsm$ at 185 120
\tiny\hair 2pt
\pinlabel $1$   at   0 -14
\pinlabel $2$   at  59 -14
\pinlabel $n-1$ at 219 -14
\pinlabel $n$   at 282 -14
\endlabellist
\emptyset \ \longmapsto \figins{-28}{1.0}{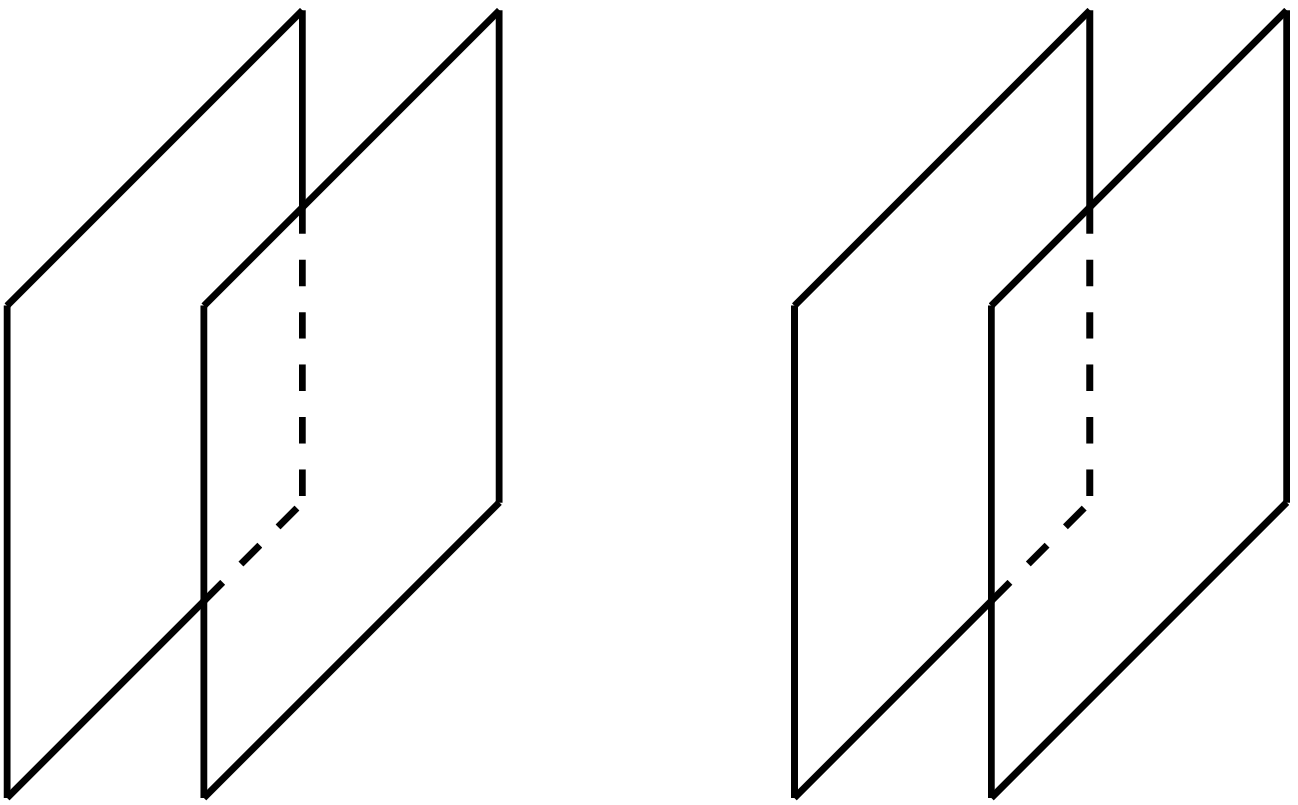}\vspace*{2ex}
\end{equation*}
\item The vertical line colored $j$ is sent to the identity cobordism of $u_j$:
\begin{equation*}
\labellist
\tiny\hair 2pt
\pinlabel $j$   at -10  60
\pinlabel $j$   at 100 -65   \pinlabel $j+1$ at 165 -66
\endlabellist
\figins{-16}{0.5}{line}\ \
\longmapsto\
\figins{-32}{1.0}{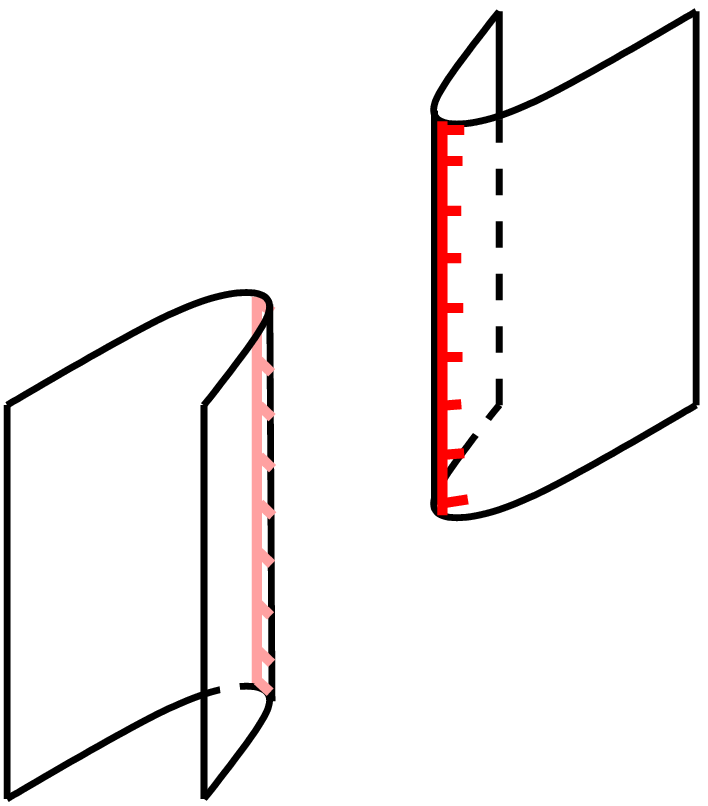}
\vspace*{2ex}
\end{equation*}
The remaining $n-2$ vertical parallel sheets on the r.h.s. are not shown for simplicity, 
a convention that we will follow from now on.

\item The \emph{StartDot} and \emph{EndDot} morphisms are sent to saddle cobordisms:
\begin{equation*}
\labellist
\tiny\hair 2pt
\pinlabel $j$   at -10  60
\pinlabel $j$   at  95 -65  \pinlabel $j+1$ at 154 -66
\endlabellist
\figins{-16}{0.5}{startdot}
\longmapsto\
\figins{-32}{1.0}{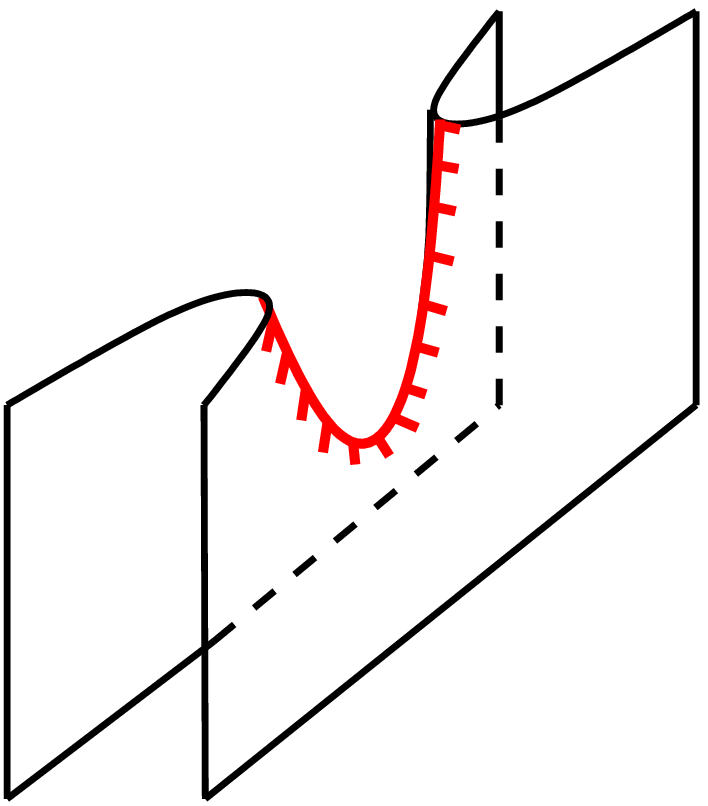}
\mspace{140mu}
\labellist
\tiny\hair 2pt
\pinlabel $j$   at -10  60
\pinlabel $j$   at  95 -65
\pinlabel $j+1$ at 154 -66
\endlabellist
\figins{-16}{0.5}{enddot}
\longmapsto\
\figins{-32}{1.0}{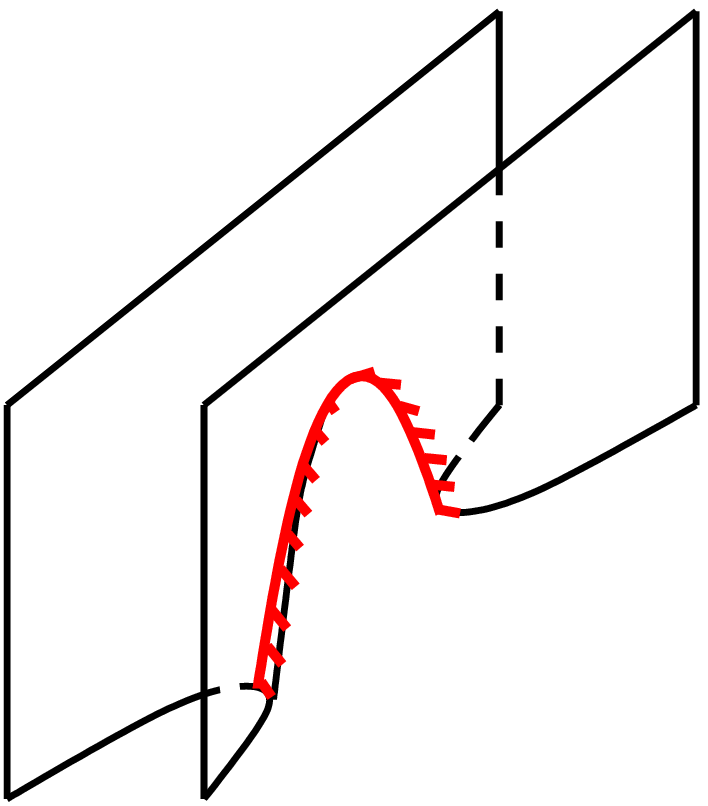}
\vspace*{2ex}
\end{equation*}

\item \emph{Merge} and \emph{Split} are sent to cup and cap cobordisms:
\begin{equation*}
\labellist
\tiny\hair 2pt
\pinlabel $j$   at  45  95
\pinlabel $j$   at 205 -65  \pinlabel $j+1$ at 274 -66
\endlabellist
\figins{-16}{0.6}{merge}
\longmapsto\
\figins{-32}{1.1}{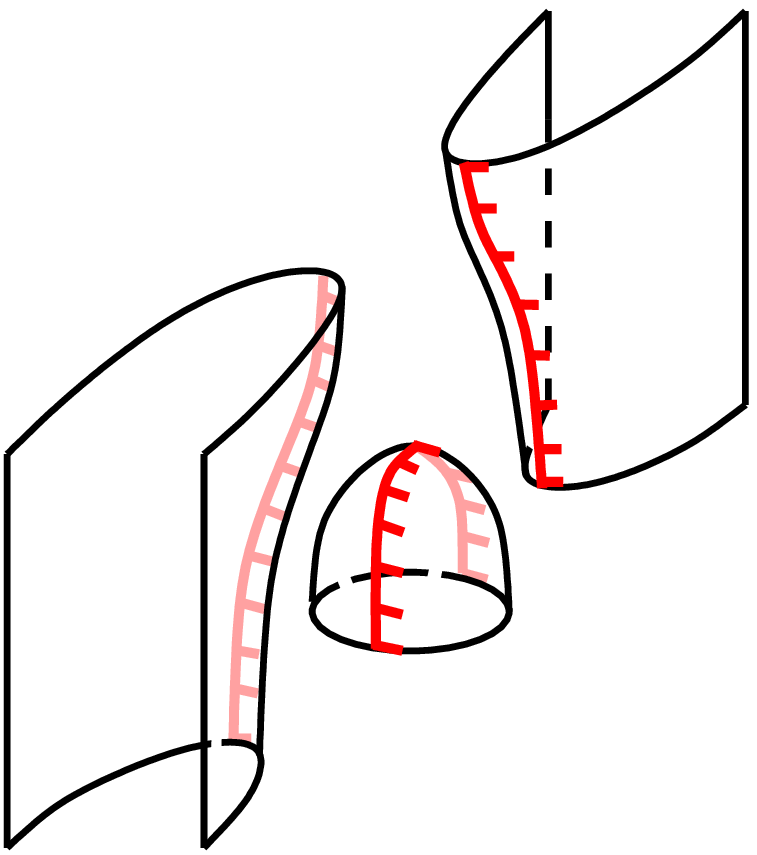}
\mspace{80mu}
\labellist
\tiny\hair 2pt
\pinlabel $j$   at  45  45
\pinlabel $j$   at 205 -65  \pinlabel $j+1$ at 279 -66
\endlabellist
\figins{-16}{0.6}{split}
\longmapsto \
\figins{-32}{1.1}{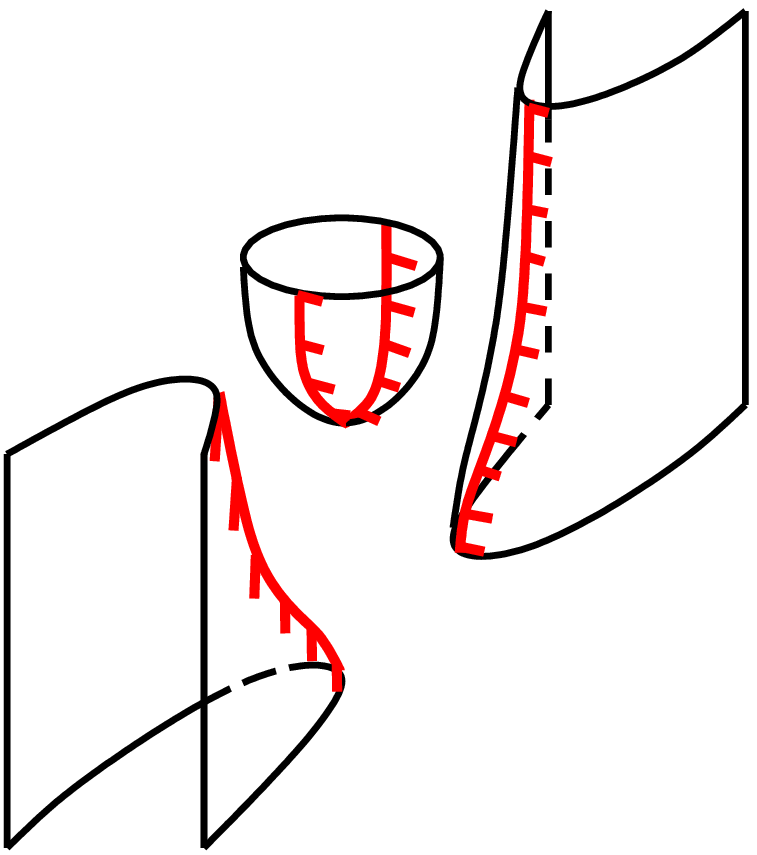}
\vspace*{2ex}
\end{equation*}

\item The \emph{4-valent vertex} with distant colors. For $j+1<k$ we have:
\begin{equation*}
\labellist
\tiny\hair 2pt
\pinlabel $k$   at  -5 -12
\pinlabel $j$   at 128 -10
\pinlabel $j$      at 205 -54
\pinlabel $j+1$    at 265 -55
\pinlabel $k$      at 365 -54
\pinlabel $k+1$    at 425 -55
\pinlabel $\dotsm$ at 325 -55
\endlabellist
\figins{-16}{0.6}{4vert}
\longmapsto\
\figins{-28}{1.0}{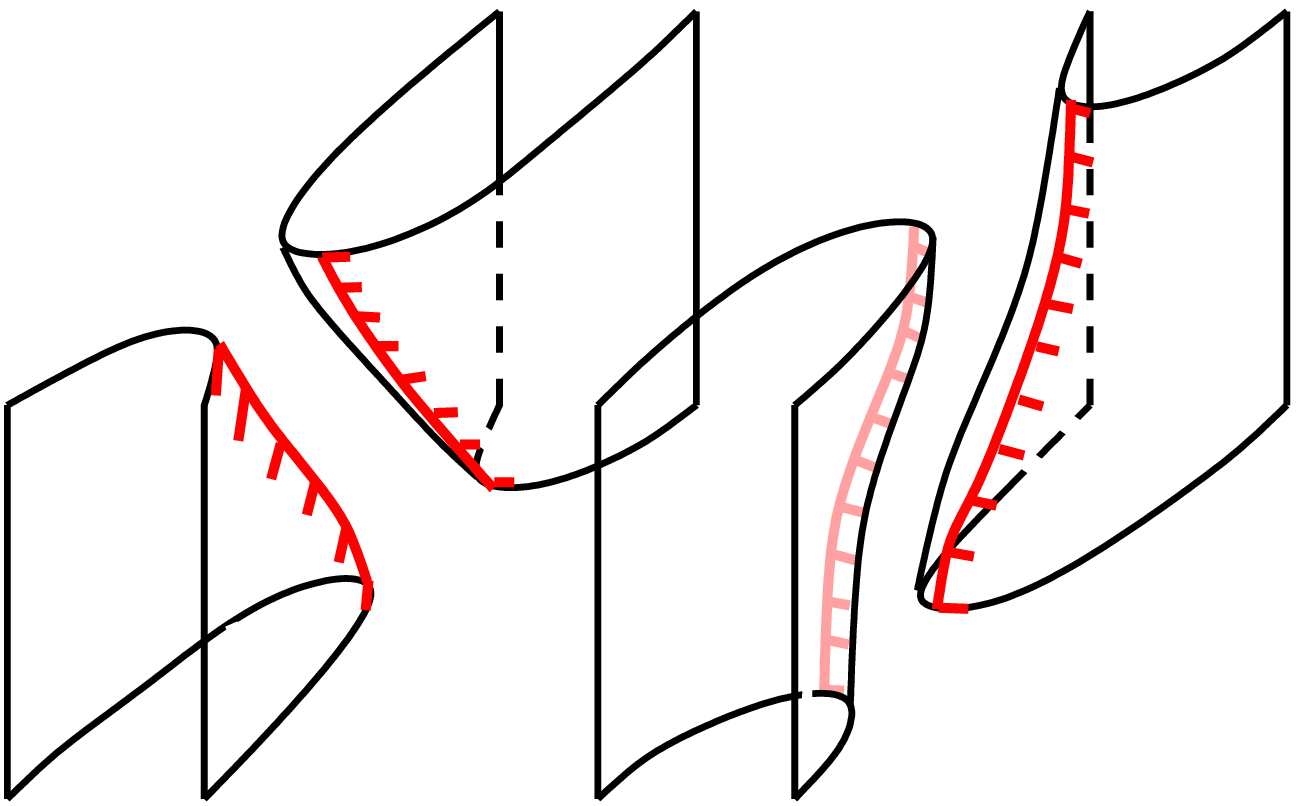}
\vspace*{2ex}
\end{equation*}
The case $j>k+1$ is given by reflection in a horizontal plane.

\item The \emph{6-valent vertices} are sent to zero:
\begin{equation*}
\figins{-18}{0.6}{6vertu}
\longmapsto\ 
0 
\end{equation*}
\end{itemize}
Notice that $\FKh$ respects the gradings of the morphisms.
Taking the quotient of $\mathcal{SC}$ by the 6-valent vertex gives a diagrammatic category 
$\mathcal{TL}$ categorifying the Temperley-Lieb algebra. According to~\cite{Elias} 
relations~\eqref{eq:dot6v} and~\eqref{eq:reid3} can be replaced by a single relation
in $\mathcal{TL}$. The functor $\FKh$ descends to a 
functor between $\mathcal{TL}$ and $\cmw$.

\begin{prop}
$\FKh$ is a monoidal functor.
\end{prop}

\begin{proof}
The assignment given by $\FKh$ clearly respects the monoidal structures of $\mathcal{SC}_1$ and $\cmw$. 
So we only need to show that $\FKh$ is a functor, i.e. it respects the relations~\eqref{eq:adj} 
to~\eqref{eq:dumbdumbsquare} of Section~\ref{sec:soergel}.
 
\medskip
%%%%%%%%%%%%%%%%%%%%%%%%%%%%%%%
\n\emph{''Isotopy relations'':} 
Relations~\eqref{eq:adj} to~\eqref{eq:v4rot} are straightforward to check and correspond to isotopies of their images 
under $\FKh$ which respect the disorientations. 
Relation~\eqref{eq:v6rot} is automatic since $\FKh$ send all terms to zero. 
For the sake of completeness we show the first equality in~\eqref{eq:adj}. We have
\begin{equation*}
\labellist
\tiny\hair 2pt
\pinlabel $j$  at -6  40        \pinlabel $j$ at 725 40
\pinlabel $j$ at 165 -52        \pinlabel $j+1$ at 214 -53
\pinlabel $j$ at 381 -52        \pinlabel $j+1$ at 433 -53
\endlabellist
\FKh\Biggl(\
\figins{-17}{0.55}{biadj-l}\
\Biggr)
=\
\figins{-32}{1.0}{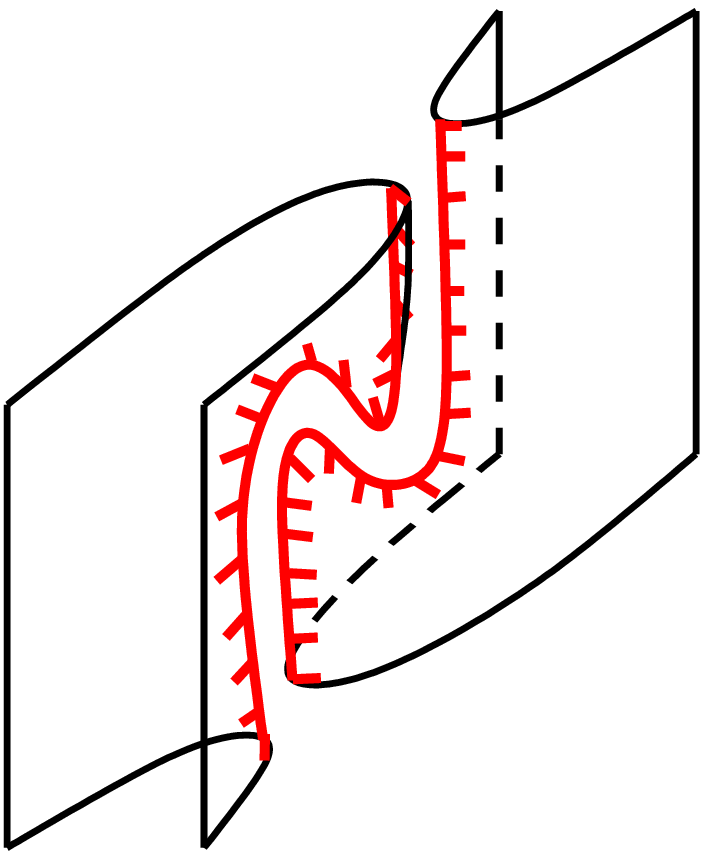}\
\cong\
\figins{-32}{1.0}{idui}\
=\
\FKh
\Biggl(\
\figins{-17}{0.55}{line}\
\Biggr)\vspace*{2ex}
\end{equation*}

%%%%%%%%%%%%%%%%%%%%%%%%%%%%%
\medskip
\n\emph{One color relations:}
For relation~\eqref{eq:dumbrot} we have 
\begin{equation*}
\FKh\Biggl(\
\figins{-14}{0.45}{dumbellh}\
\Biggr)
\cong
\FKh\Biggl(\
\figins{-14}{0.45}{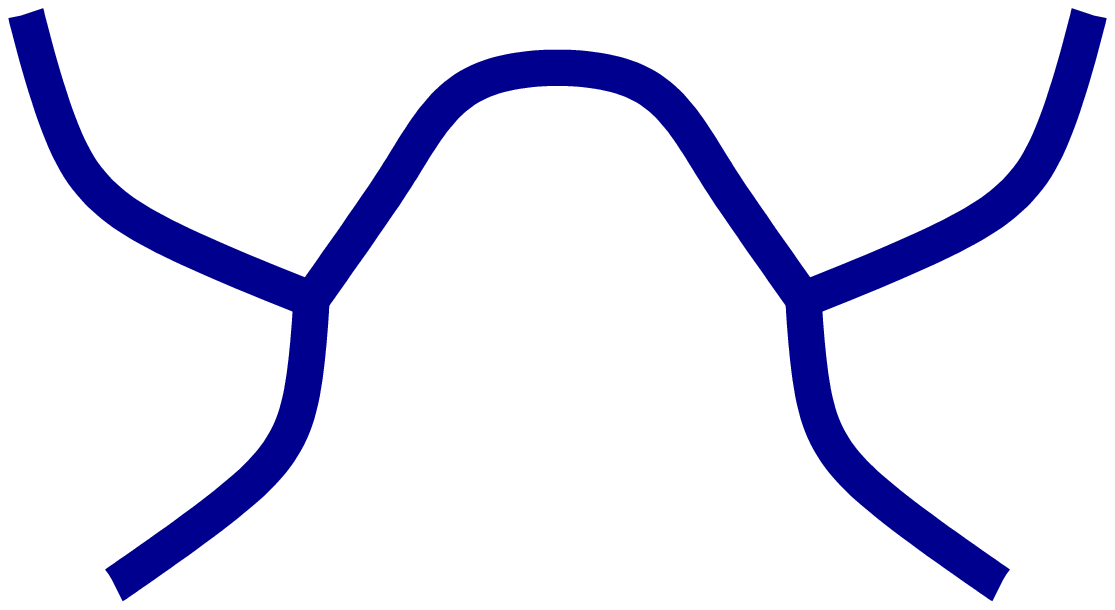}\
\Biggr)
\cong
\FKh\Biggl(\
\figins{-14}{0.45}{dumbells}\
\Biggr)\ ,
\end{equation*}
where the first equivalence follows from relations~\eqref{eq:adj} and~\eqref{eq:v3rot} and the second from isotopy of the 
cobordisms involved.

For relation~\eqref{eq:lollipop} we have
\begin{equation*}
\labellist
\tiny\hair 2pt
\pinlabel $j$ at 26 45            
\pinlabel $j$ at 240 -98 \pinlabel $j+1$ at 341 -99
\endlabellist
\FKh\Biggl(\
\figins{-16}{0.5}{lollipop-u}\
\Biggr)
=\
\figins{-32}{1.0}{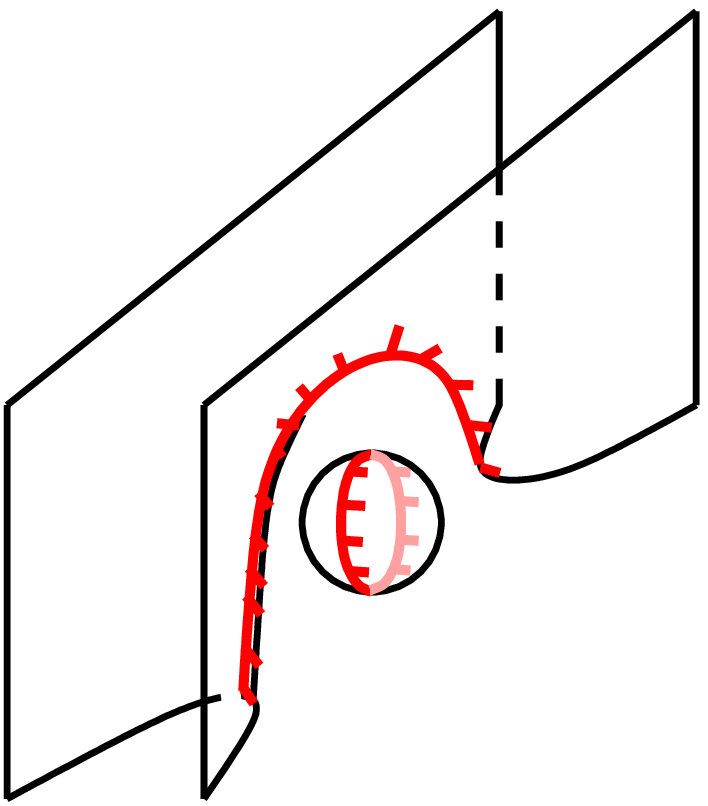}\
=\
0\mspace{30mu}\text{by relations~\eqref{eq:disorc} and~\eqref{eq:BN-s}}.
\vspace*{2ex}
\end{equation*}

Relation~\eqref{eq:deltam} requires some more work. We have
\begin{align*}
\labellist
\tiny\hair 2pt
\pinlabel $j$ at -2 22    \pinlabel $j$ at -2 130  
\pinlabel $j$ at 125 -70  \pinlabel $j+1$ at 195 -71
\pinlabel $j$ at 460 -70  \pinlabel $j+1$ at 525 -71
\pinlabel $j$ at 232 -393 \pinlabel $j+1$ at 295 -395
\pinlabel $j$ at 494 -393 \pinlabel $j+1$ at 560 -395
\endlabellist
\FKh\Biggl(\
\figins{-17}{0.55}{matches-ud}\
\Biggr)
&=\
\figins{-32}{1.0}{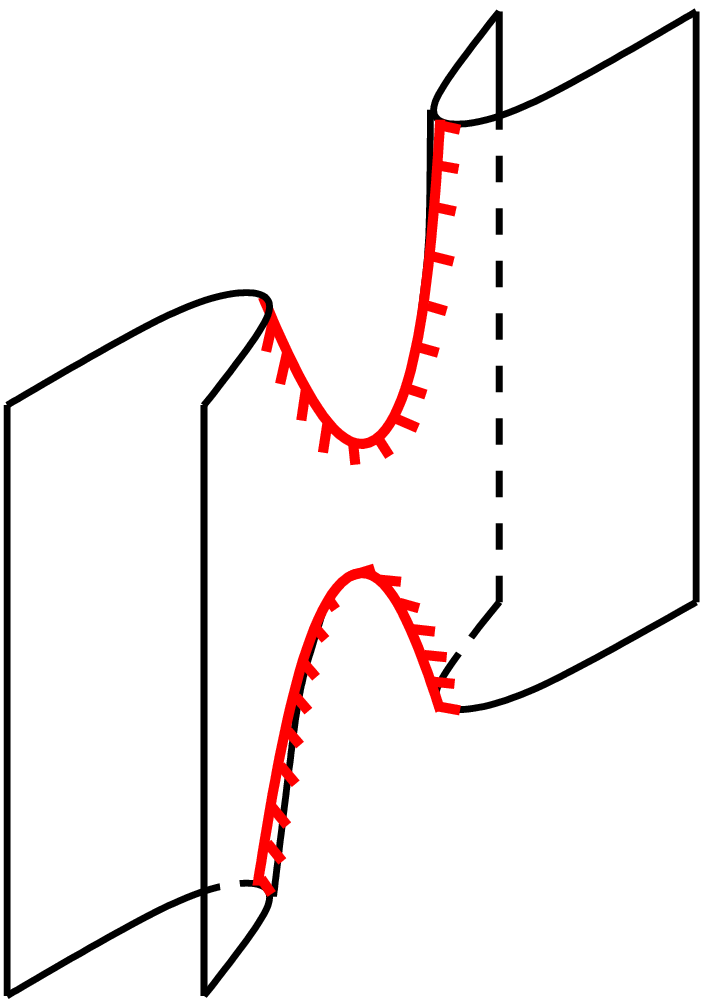}\
=\ -i\
\figins{-32}{1.0}{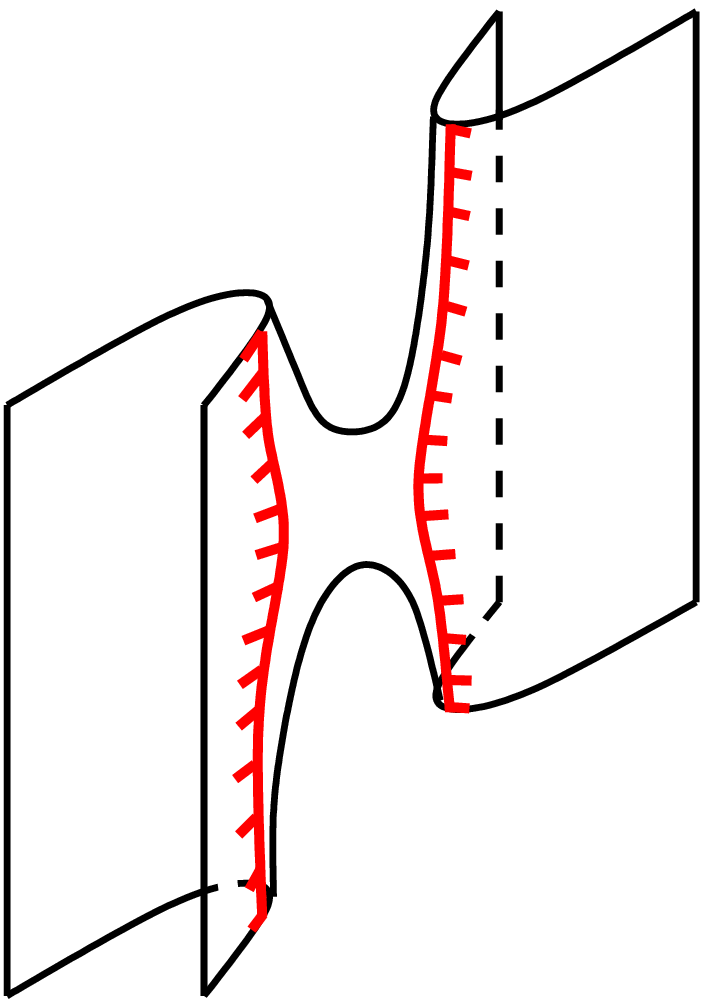}\
\\[1.5ex]
&=
-i\
\left(\
\figins{-32}{1.0}{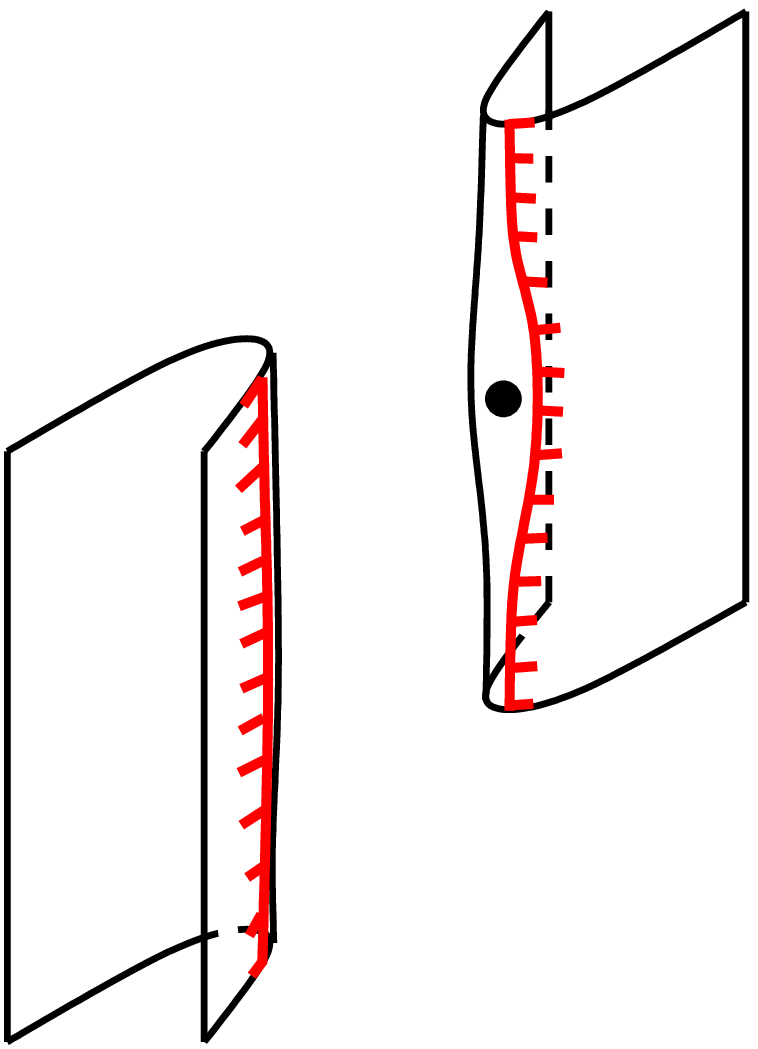}\
+
\figins{-32}{1.0}{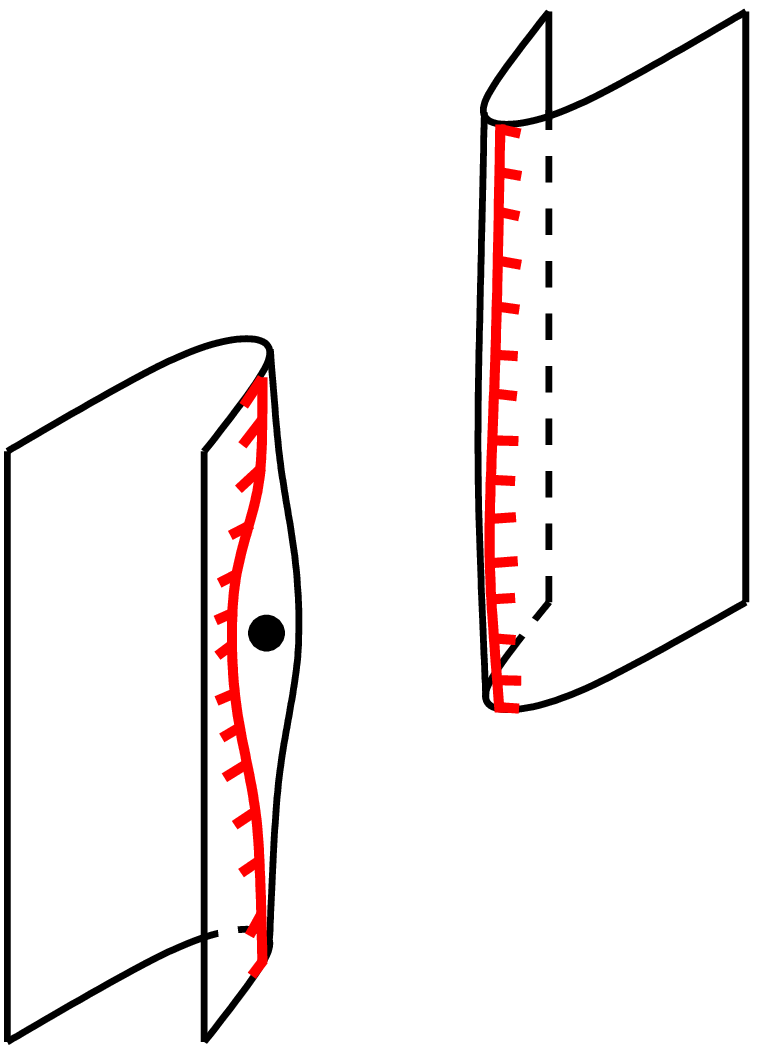}\
\right)
\vspace*{3ex}
\end{align*}

\n where the second equality follows from the disoriented relation~\eqref{eq:disors} and the third from the BN 
relation~\eqref{eq:BN-cn}. We also have
\begin{align*}
\labellist
\tiny\hair 2pt
\pinlabel $j$ at -13   16
\pinlabel $j$ at 150 -130   \pinlabel $j+1$ at  300 -131
\pinlabel $j$ at 545 -130   \pinlabel $j+1$ at  700 -131
\pinlabel $j$ at 915 -130   \pinlabel $j+1$ at 1070 -131
\endlabellist
\FKh\biggl(\
\figins{-6}{0.25}{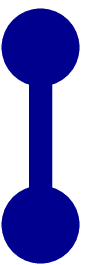}\
\biggr)\
=\
\figins{-32}{1.0}{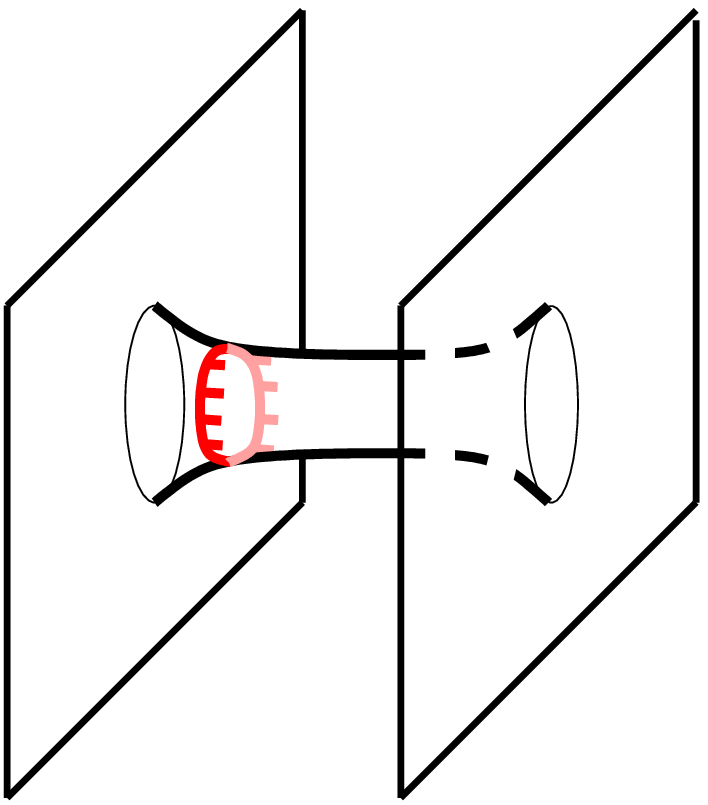}
=\
-i\
\figins{-32}{1.0}{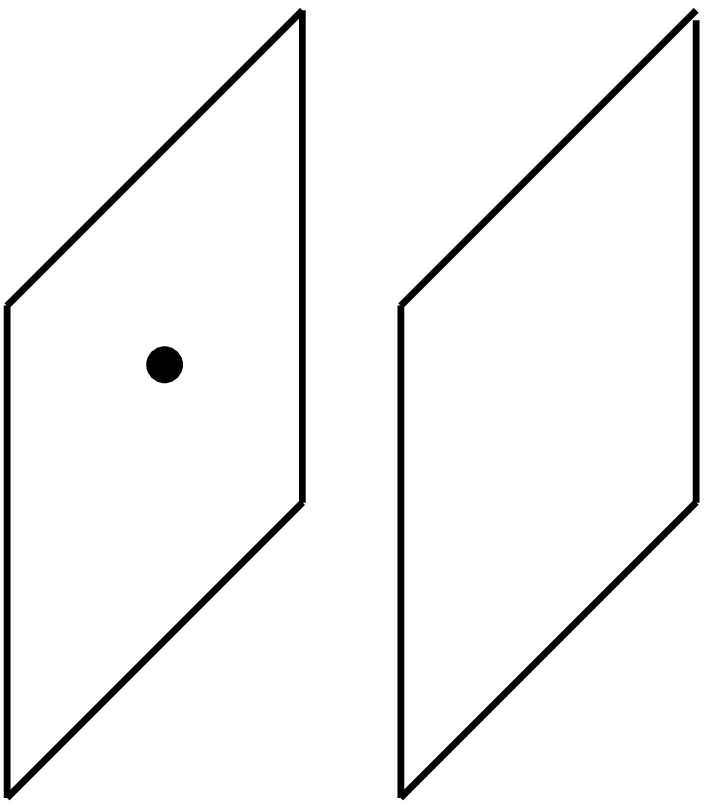}\
+\ i\
\figins{-32}{1.0}{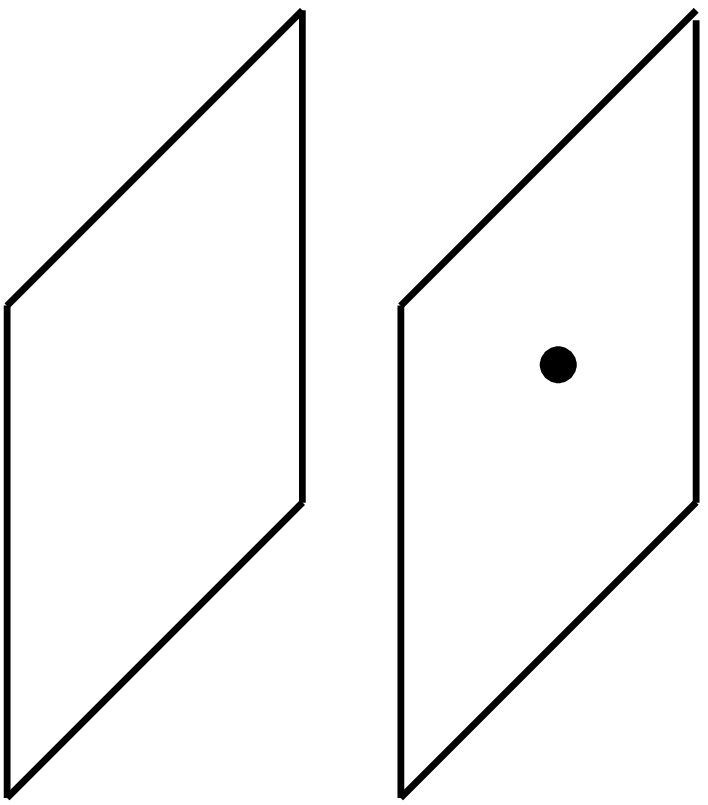}
\vspace*{2ex}
\end{align*}
and therefore
\begin{align}
\label{eq:edge-dots}
\labellist
\tiny\hair 2pt
\pinlabel $j$ at -10   60
\pinlabel $j$ at 235  -78 
\pinlabel $j+1$ at 299  -79
\endlabellist
\FKh\Biggl(\
\figins{-17}{0.55}{edge-startenddot}\
\Biggr)
&=
\ -2i\
\figins{-32}{1.0}{sheet-dot}
\\[1.5ex]\intertext{and}\displaybreak[0]
\label{eq:dots-edge}
\labellist
\tiny\hair 2pt
\pinlabel $j$ at 65   60  
\pinlabel $j$ at 235  -78  \pinlabel $j+1$ at 299  -79
\endlabellist
\FKh\Biggl(\
\figins{-17}{0.55}{startenddot-edge}\
\Biggr)
&=
\ -2i\
\figins{-32}{1.0}{dot-sheet}\ .
\end{align}
We thus have that
\begin{equation*}
\FKh
\Biggl(
\figins{-17}{0.55}{startenddot-edge}\
\Biggr)
+\
\FKh
\Biggl(\
\figins{-17}{0.55}{edge-startenddot}\
\Biggr)
=\ 2\ 
\FKh
\Biggl(\
\figins{-17}{0.55}{matches-ud}\
\Biggr).
\end{equation*}

%%%%%%%%%%%%%%%%%%%%%%%%%%%%
\medskip
\n\emph{Two distant colors:} 
Relations~\eqref{eq:reid2dist} to~\eqref{eq:slide3v} correspond to isotopies of the 
cobordisms involved and are straightforward to check.

%%%%%%%%%%%%%%%%%%%%%%%%%
\medskip
\n\emph{Adjacent colors:}
We prove the case where ''blue`` corresponds to $j$ and ''red`` corresponds to $j+1$. 
The relations with colors reversed are proved the same way.
To prove relation~\eqref{eq:dot6v} we first notice that
\begin{equation*}
\labellist
\tiny\hair 2pt
\pinlabel $j$   at 270 -121
\pinlabel $j+1$ at 352 -122
\pinlabel $j+2$ at 456 -122
\endlabellist
\FKh\Biggl(\
\figins{-10}{0.3}{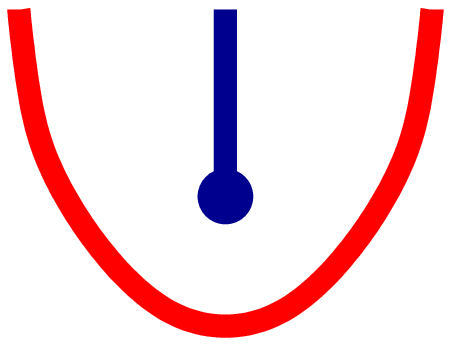}\
\Biggr)\
\cong\
\figins{-32}{1.0}{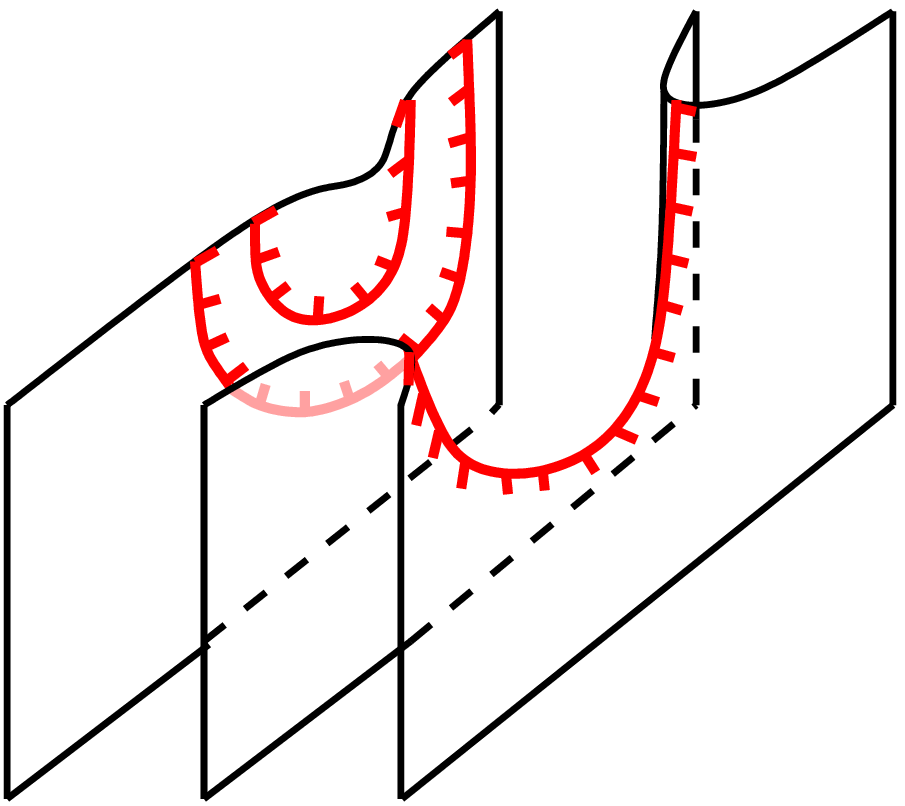}
\end{equation*}
which means that
\begin{equation*}
\labellist
\tiny\hair 2pt
\pinlabel $j$   at 270 -121
\pinlabel $j+1$ at 352 -122
\pinlabel $j+2$ at 456 -122
\endlabellist
\FKh\Biggl(\
\figins{-12}{0.4}{capcupdot}\
\Biggr)\
\cong\
\figins{-34}{1.2}{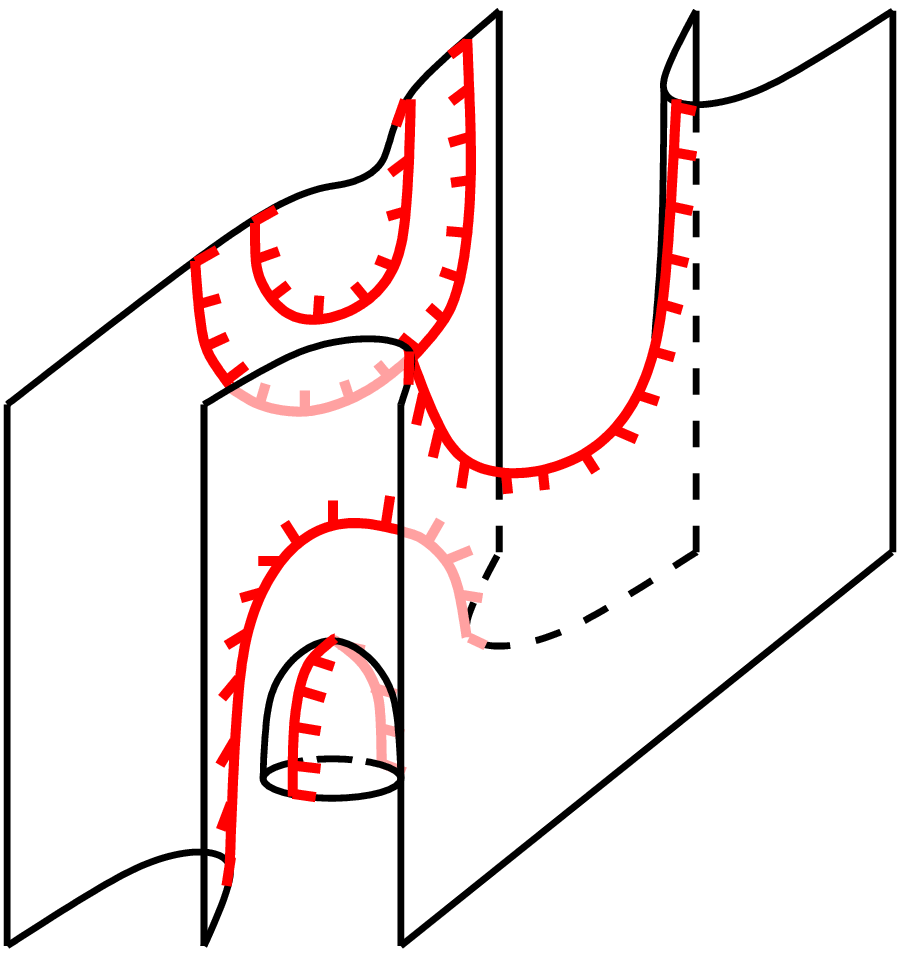}
\vspace*{2ex}
\end{equation*}

On the other side we have
\begin{equation*}
\labellist
\tiny\hair 2pt
\pinlabel $j$   at 252 -118
\pinlabel $j+1$ at 318 -119
\pinlabel $j+2$ at 400 -119
\endlabellist
\FKh\Biggl(\
\figins{-12}{0.4}{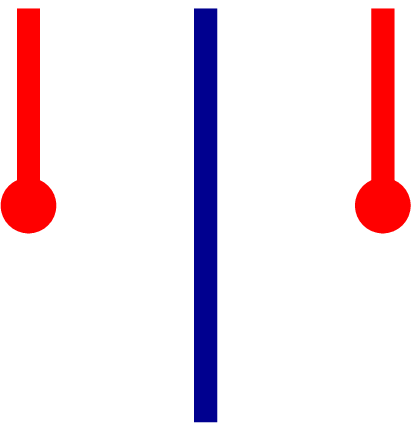}\
\Biggr)\
\cong\
\figins{-34}{1.2}{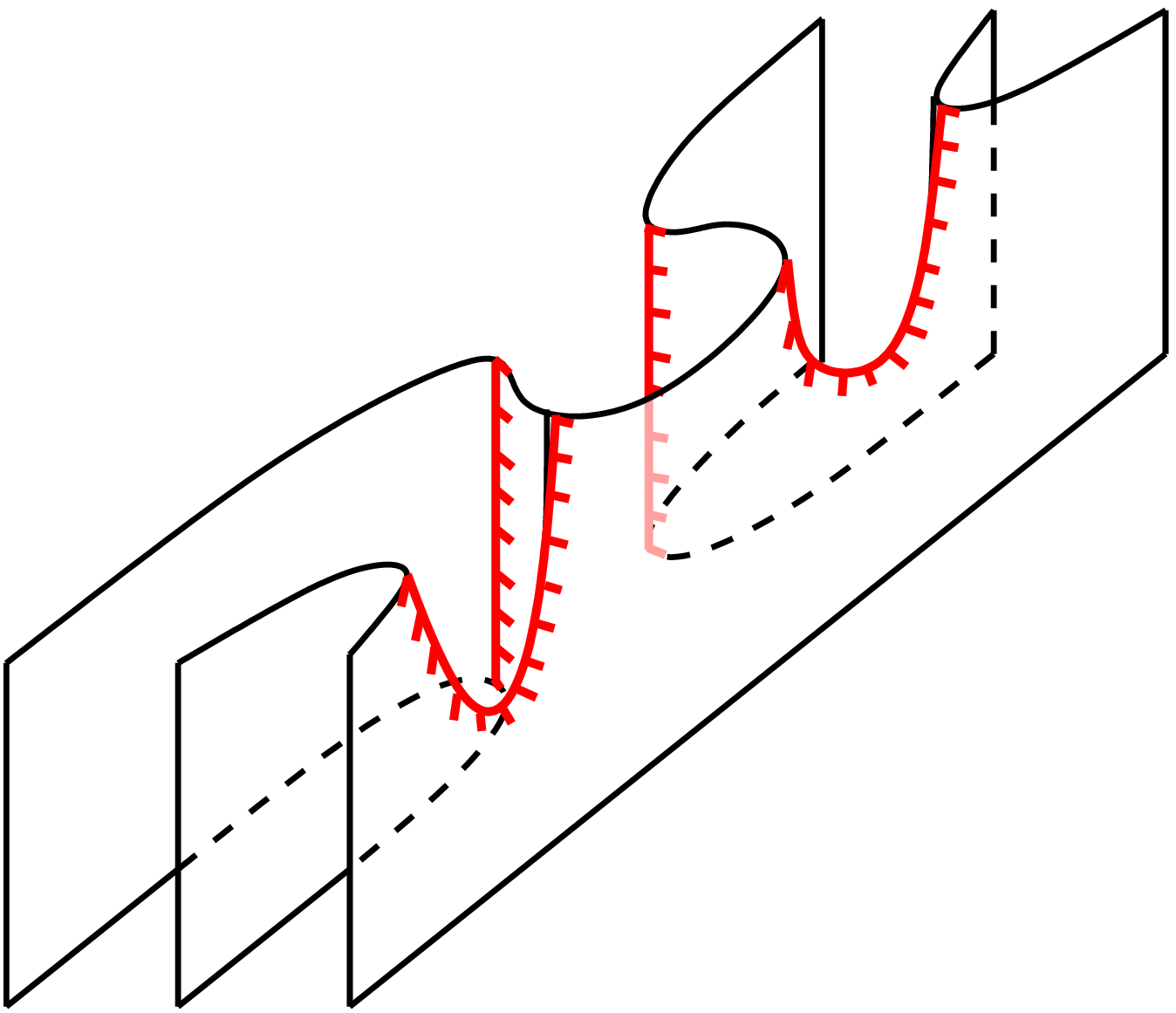}
\vspace*{2ex}
\end{equation*}
which, using isotopies and the disorientation relation~\eqref{eq:disors} twice, can be seen to be equivalent to 
\begin{equation*}
\labellist
\tiny\hair 2pt
\pinlabel $j$   at -10 -19
\pinlabel $j+1$ at  46 -20
\pinlabel $j+2$ at 126 -20
\endlabellist
-\
\figins{-32}{1.0}{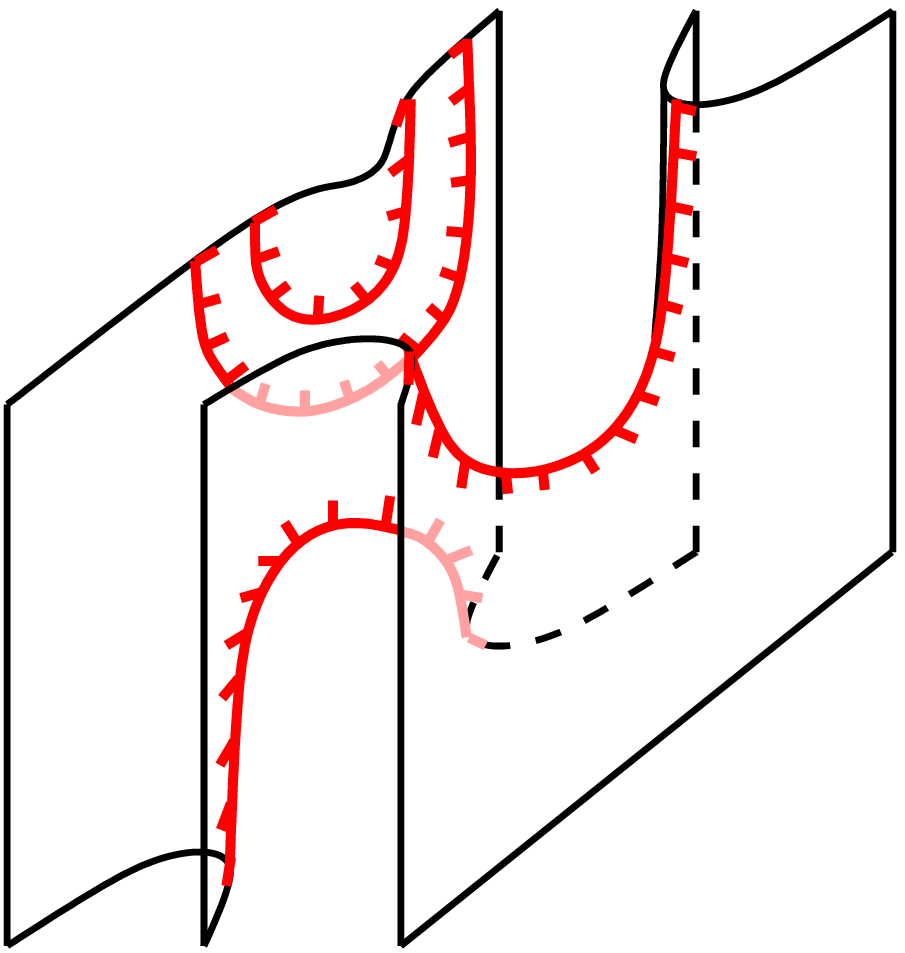}
\vspace*{2ex}
\end{equation*}
which equals
\begin{equation*}
-\FKh
\Biggl(\
\figins{-12}{0.4}{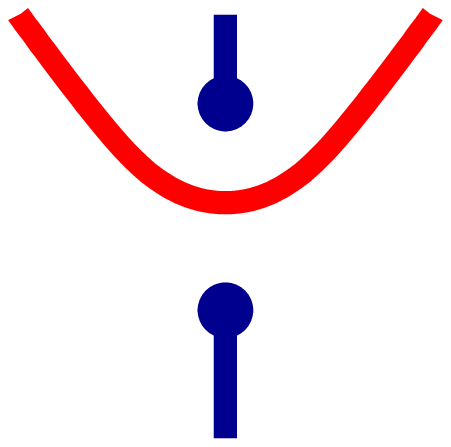}\
\Biggr).
\end{equation*}

This implies that
\begin{equation*}
0=
\FKh\Biggl(\
\figins{-17}{0.55}{6vertdotd}\
\Biggr)
=\
\FKh\Biggl(\
\figins{-17}{0.55}{mergedots}\
\Biggr)
+\
\FKh\Biggl(\
\figins{-17}{0.55}{capcupdot}\
\Biggr) .
\end{equation*}

We now prove relation~\eqref{eq:reid3}. We have isotopy equivalences
\begin{align*}
\labellist
\tiny\hair 2pt
\pinlabel $j$ at 265   -61  \pinlabel $j+1$ at 339   -62  \pinlabel $j+2$ at 439   -62
\endlabellist
\FKh
\left(\
\figins{-26}{0.8}{id-r3}\
\right)
& \cong\
\figins{-34}{1.2}{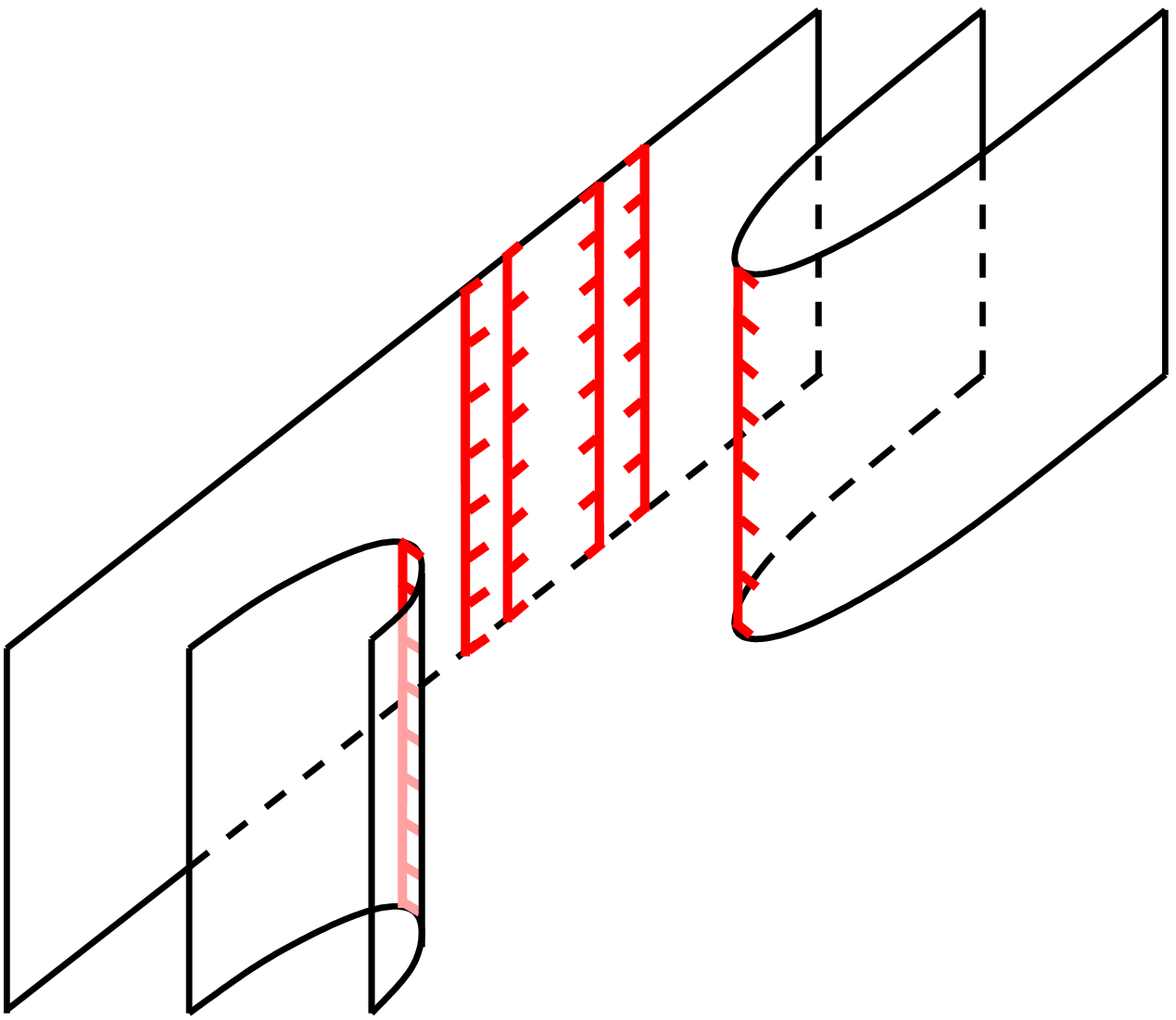}
\intertext{and}
%\\[1.5ex]
\displaybreak[0]
\labellist
\tiny\hair 2pt
\pinlabel $j$ at 265   -61  \pinlabel $j+1$ at 339   -62  \pinlabel $j+2$ at 439   -62
\endlabellist
\FKh
\left(\
\figins{-26}{0.8}{dumbell-dd}\
\right)
& \cong\
\figins{-34}{1.2}{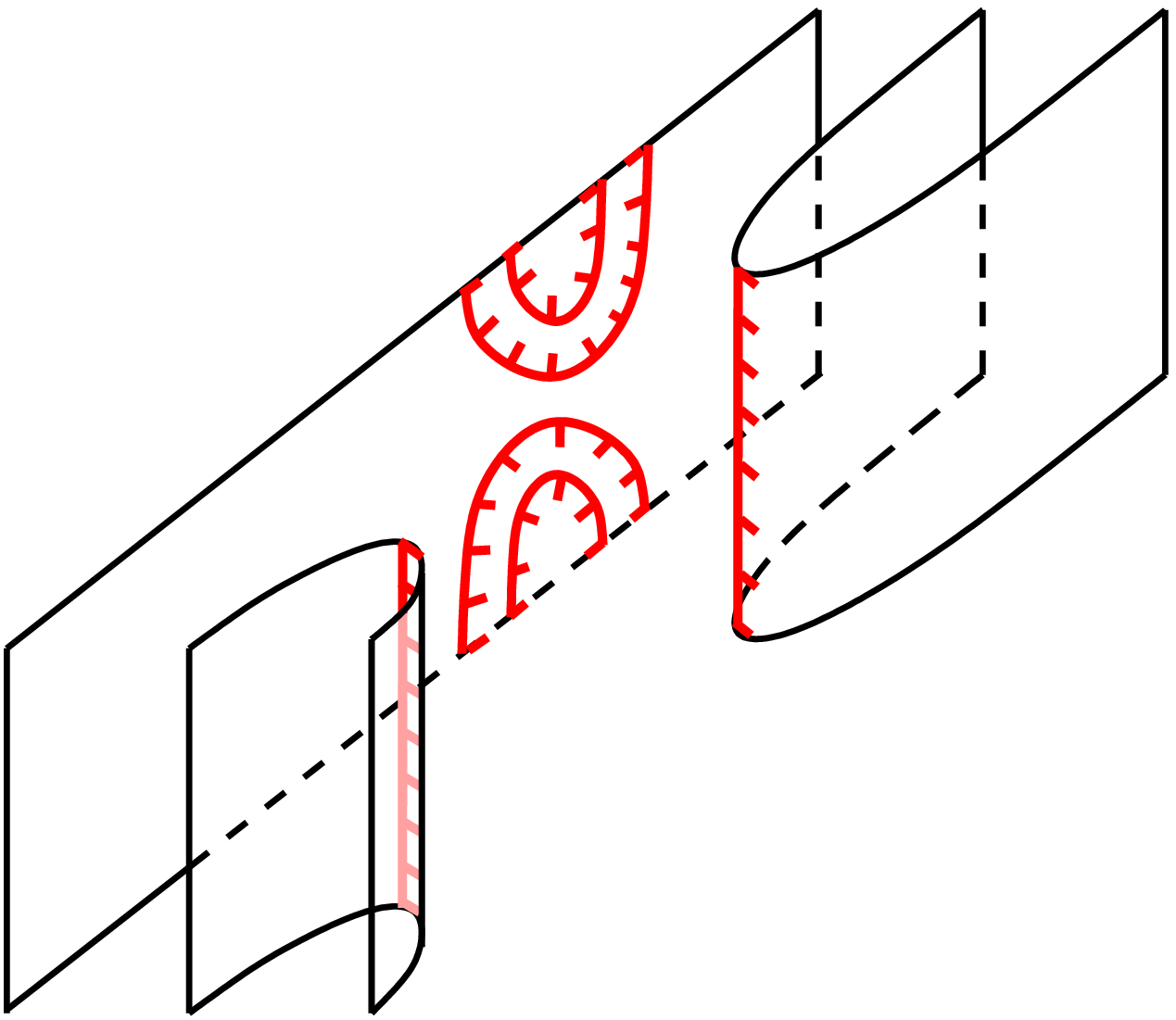}
\\[1.5ex]\displaybreak[0]
\labellist
\tiny\hair 2pt
\pinlabel $j$ at -10 -21  \pinlabel $j+1$ at 54 -22  \pinlabel $j+2$ at 129 -22
\endlabellist
&=\
-\
\figins{-34}{1.2}{sheet-disv-idui}
\vspace*{2ex}
\end{align*}
Therefore we see that 
\begin{equation*}
0=\
\FKh
\left(\
\figins{-26}{0.8}{reid3}\
\right)
=\
\FKh
\left(\
\figins{-26}{0.8}{id-r3}\
\right)
+\
\FKh
\left(\
\figins{-26}{0.8}{dumbell-dd}\
\right).
\end{equation*}

The functor $\FKh$ sends both sides of relation~\eqref{eq:dumbsq} to zero and so there is nothing to prove here. To prove relation~\eqref{eq:slidenext} we start with the equivalence
\begin{align*}
\labellist
\tiny\hair 2pt
\pinlabel $j$ at 165  -77   \pinlabel $j+1$ at 219  -78   \pinlabel $j+2$ at 295  -78
\pinlabel $j$ at 225 -397   \pinlabel $j+1$ at 279 -398   \pinlabel $j+2$ at 355 -398
\pinlabel $j$ at 655 -397   \pinlabel $j+1$ at 709 -398   \pinlabel $j+2$ at 785 -398
\endlabellist
\FKh\Biggl(\
\figins{-17}{0.55}{sedot-edge-d}\
\Biggr)\
&=\
\figins{-32}{1.0}{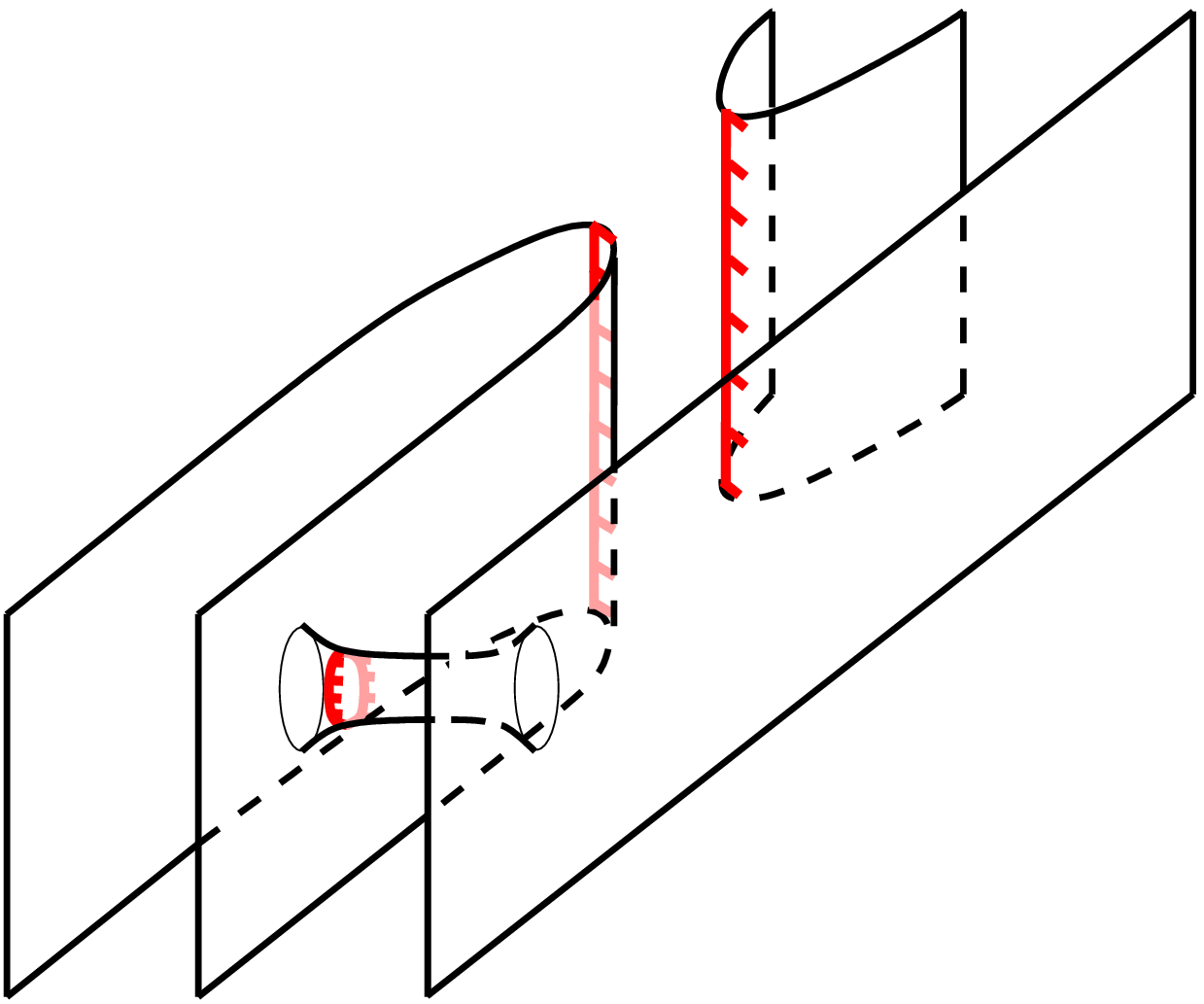}
\\[1.5ex]\displaybreak[0]
&\cong\
-i\
\figins{-32}{1.0}{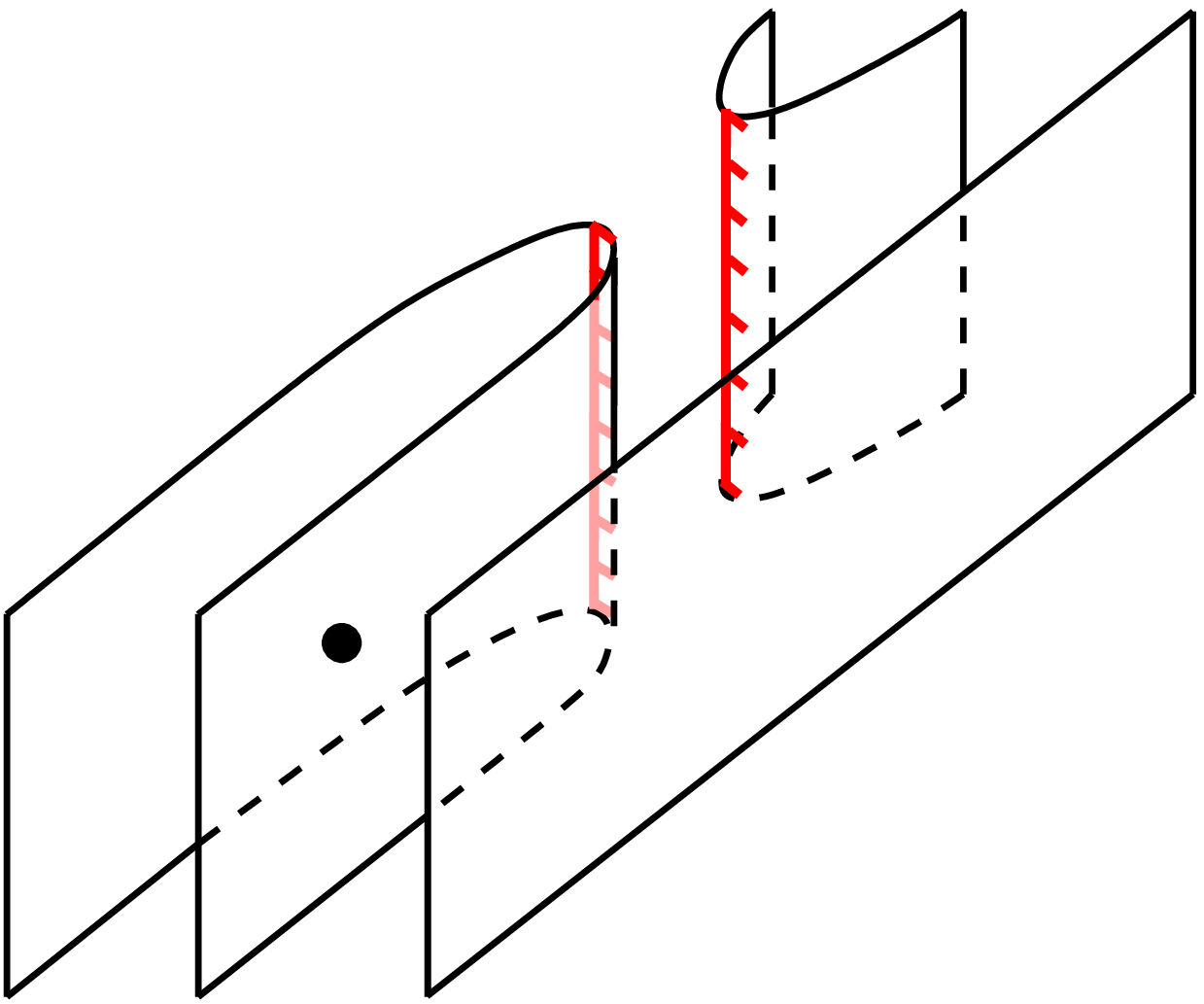}\
+\
i\
\figins{-32}{1.0}{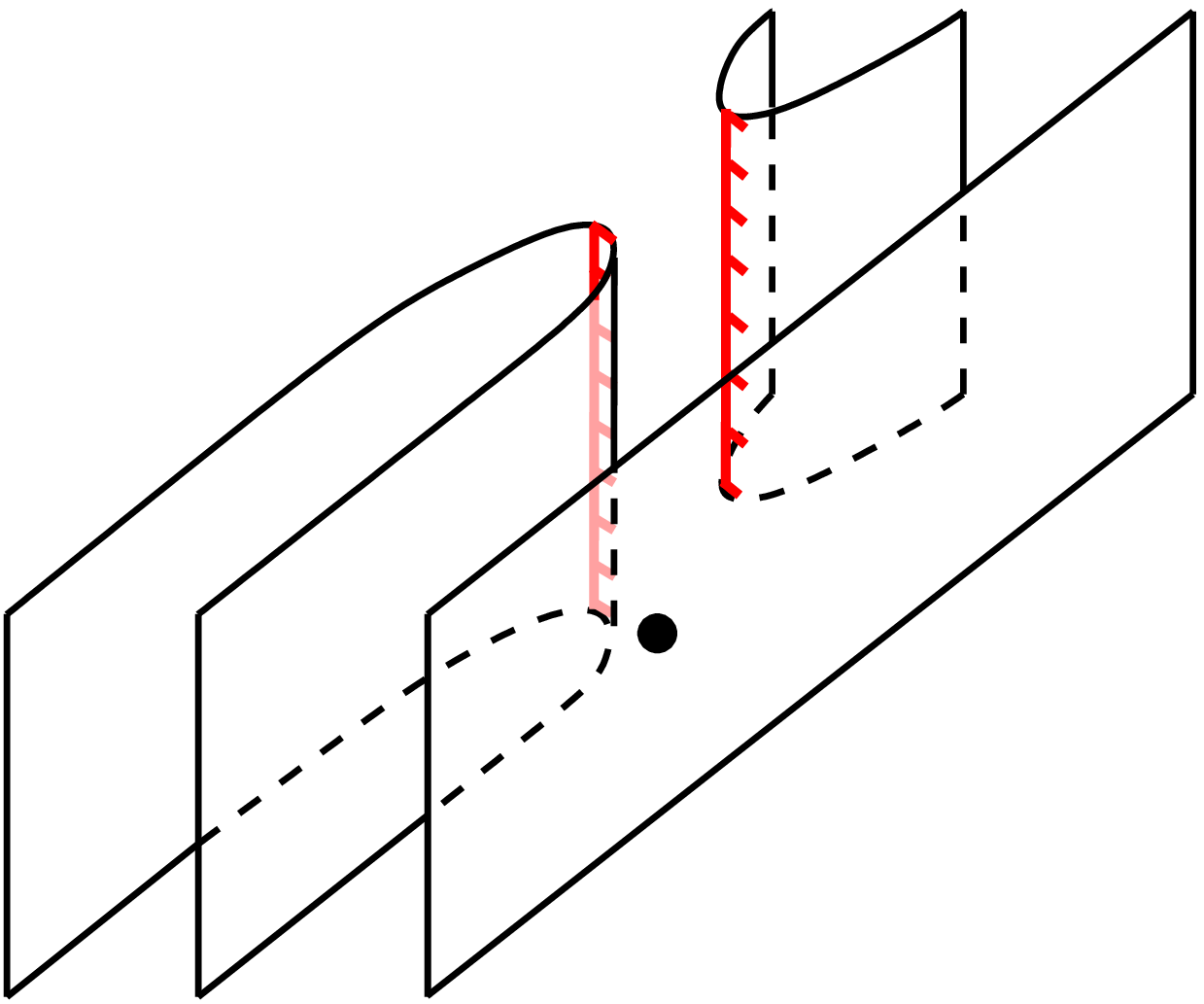}
\vspace*{2ex}
\end{align*}

\bigskip

\n which is a consequence of the neck cutting relation~\eqref{eq:BN-cn} and the disorientation 
relations~\eqref{eq:disorc} and~\eqref{eq:disorp}. We also have
\begin{align*}
\labellist
\tiny\hair 2pt
\pinlabel $j$ at 225 -77   \pinlabel $j+1$ at 279 -78   \pinlabel $j+2$ at 355 -78
\pinlabel $j$ at 655 -77   \pinlabel $j+1$ at 714 -78   \pinlabel $j+2$ at 790 -78
\endlabellist
\FKh\Biggl(\
\figins{-17}{0.55}{edge-sedot-d}\
\Biggr)\
\cong\
-i\
\figins{-32}{1.0}{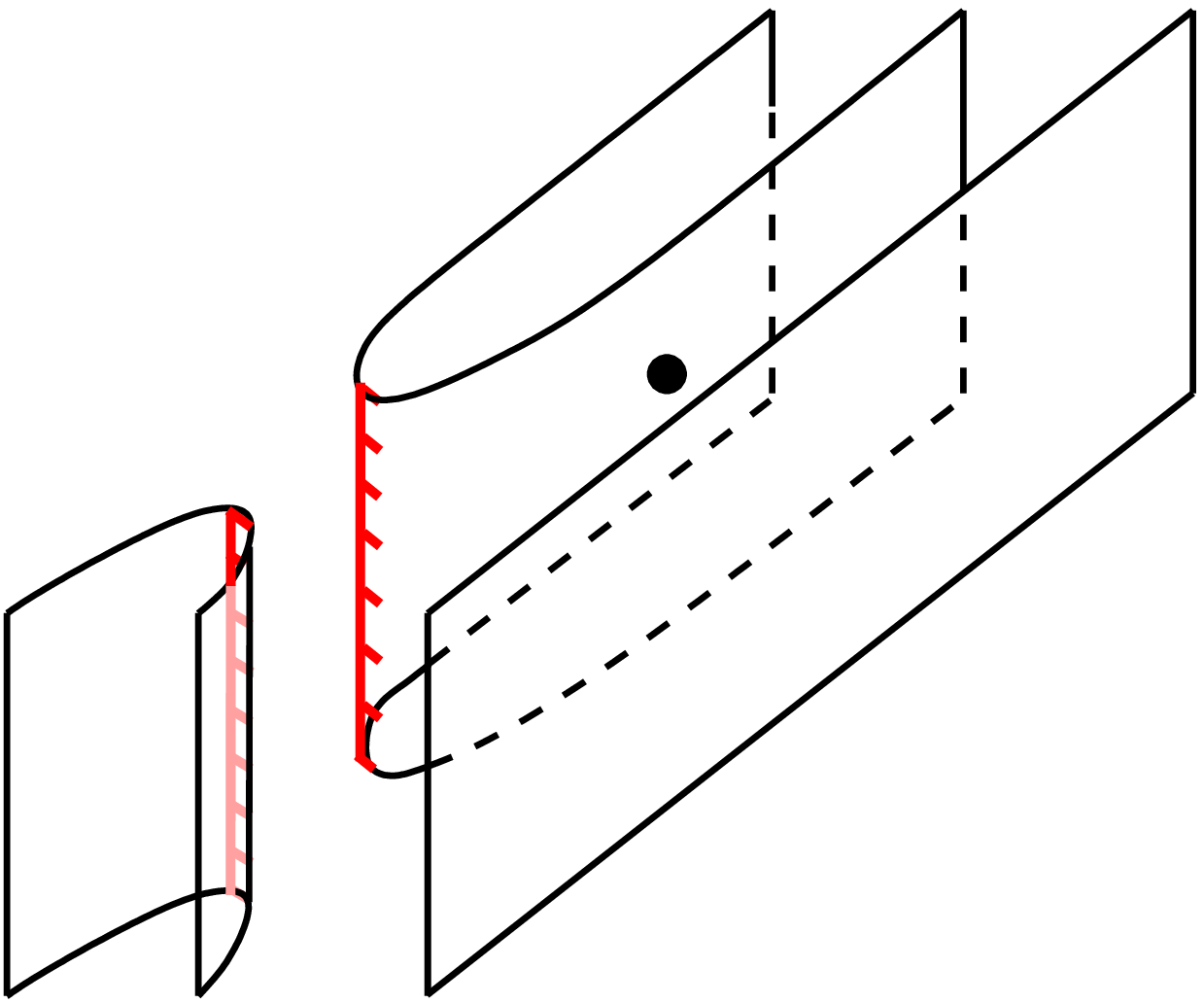}\
+\
i\
\figins{-32}{1.0}{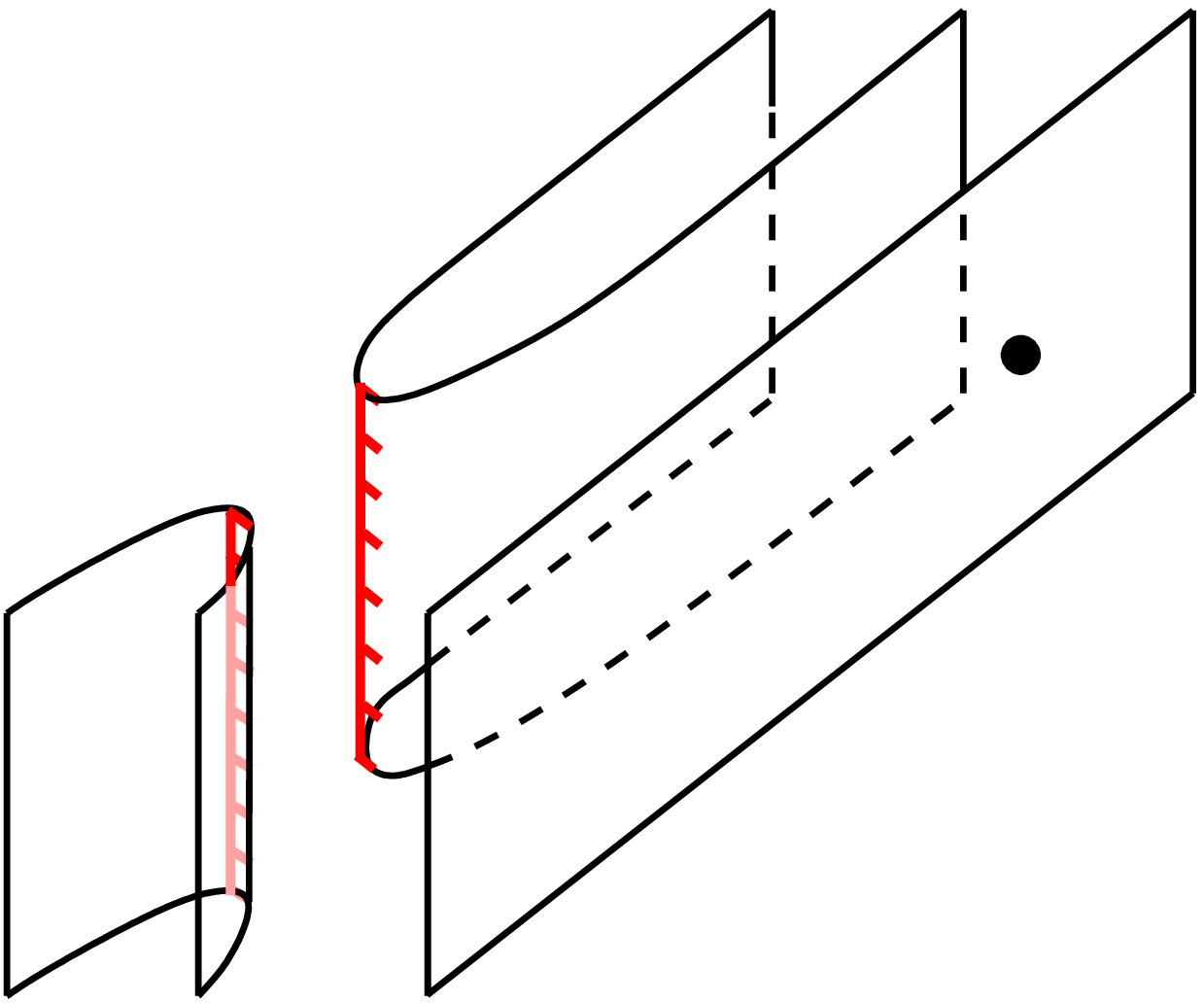}\ .
\vspace*{2ex}
\end{align*}

\bigskip

\n Comparing with Equations~\eqref{eq:edge-dots} and~\eqref{eq:dots-edge} and using the disoriented 
relation~\eqref{eq:disorp} we get
\begin{equation*}
\FKh\Biggl(\
\figins{-17}{0.55}{sedot-edge-d}\
\Biggr)\
-\
\FKh
\Biggl(\
\figins{-17}{0.55}{edge-sedot-d}\
\Biggr)\
=\
\frac{1}{2}
\FKh
\Biggl(\
\figins{-18}{0.55}{edge-startenddot}\
\Biggr)\
-\
\frac{1}{2}
\FKh
\Biggl(\
\figins{-18}{0.55}{startenddot-edge}\
\Biggr).
\end{equation*}

%%%%%%%%%%%%%%%%%%%%%%%%%%%%%%%%%%%%%%%%%%
\medskip
\n\emph{Relations involving three colors:} 
Functor $\FKh$ sends to zero both sides of relations~\eqref{eq:slide6v} 
and~\eqref{eq:dumbdumbsquare}. 
Relation~\eqref{eq:slide4v} follows from isotopies of the cobordisms involved.
\end{proof}

%%%%%%%%%%%%%%%%%%%%%%%%%%%%%%%%%%%%%%%%
%%%                                  %%%
%%%            sl(3)                 %%%
%%%                                  %%%
%%%%%%%%%%%%%%%%%%%%%%%%%%%%%%%%%%%%%%%%
\section{The $sl(3)$ case}
\label{sec:sl3}

\subsection{The category $\foamt$ of $sl(3)$ foams}\label{ssec:sl3}

In this subsection we review the category $\foamt$ of $sl(3)$-foams introduced by the 
author and Mackaay in~\cite{MV1}. 
This category was introduced to universally deform Khovanov's construction in~\cite{Kh-sl3} 
leading to the $sl(3)$-link 
homolo\-gy theory.

We follow the conventions and notation of~\cite{MV1}. Recall that a \emph{web} is a 
trivalent planar graph, where near 
each vertex all edges are oriented away from it or all edges are oriented towards it. 
We also allow webs without 
vertices, which are oriented loops. A \emph{pre-foam} is a cobordism with singular arcs 
between two webs. 
A singular arc in a pre-foam $f$ is the set of points of $f$ which have a neighborhood 
homeomorphic to the letter Y 
times an interval. Singular arcs are disjoint. 
Interpreted as morphisms, we read pre-foams from bottom to top by convention, foam 
composition consists of placing one 
pre-foam on top of the other. The orientation of the singular arcs is by convention as in 
the \emph{zip} and the 
\emph{unzip},
\begin{equation*}
\figins{-20}{0.5}{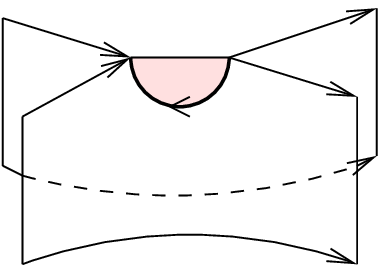}
\mspace{35mu}\text{and}\mspace{35mu}
\figins{-20}{0.5}{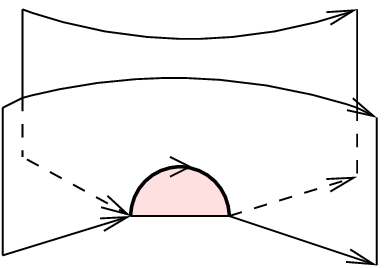}\ ,
\end{equation*}
respectively. Pre-foams can have dots which can move freely on the facet to which they 
belong but are not allowed to 
cross singular arcs. A \emph{foam} is an isotopy class of pre-foams. Let $\bC[a,b,c]$ be the 
ring of polynomials 
in $a,b,c$ with coefficients in $\bC$.

We impose the set of relations $\ell=(3D,CN,S,\Theta)$ on foams, as well as the \emph{closure relation}, 
which are explained below. 
\begin{gather*}
%3Dot reduction
\figins{-7}{0.25}{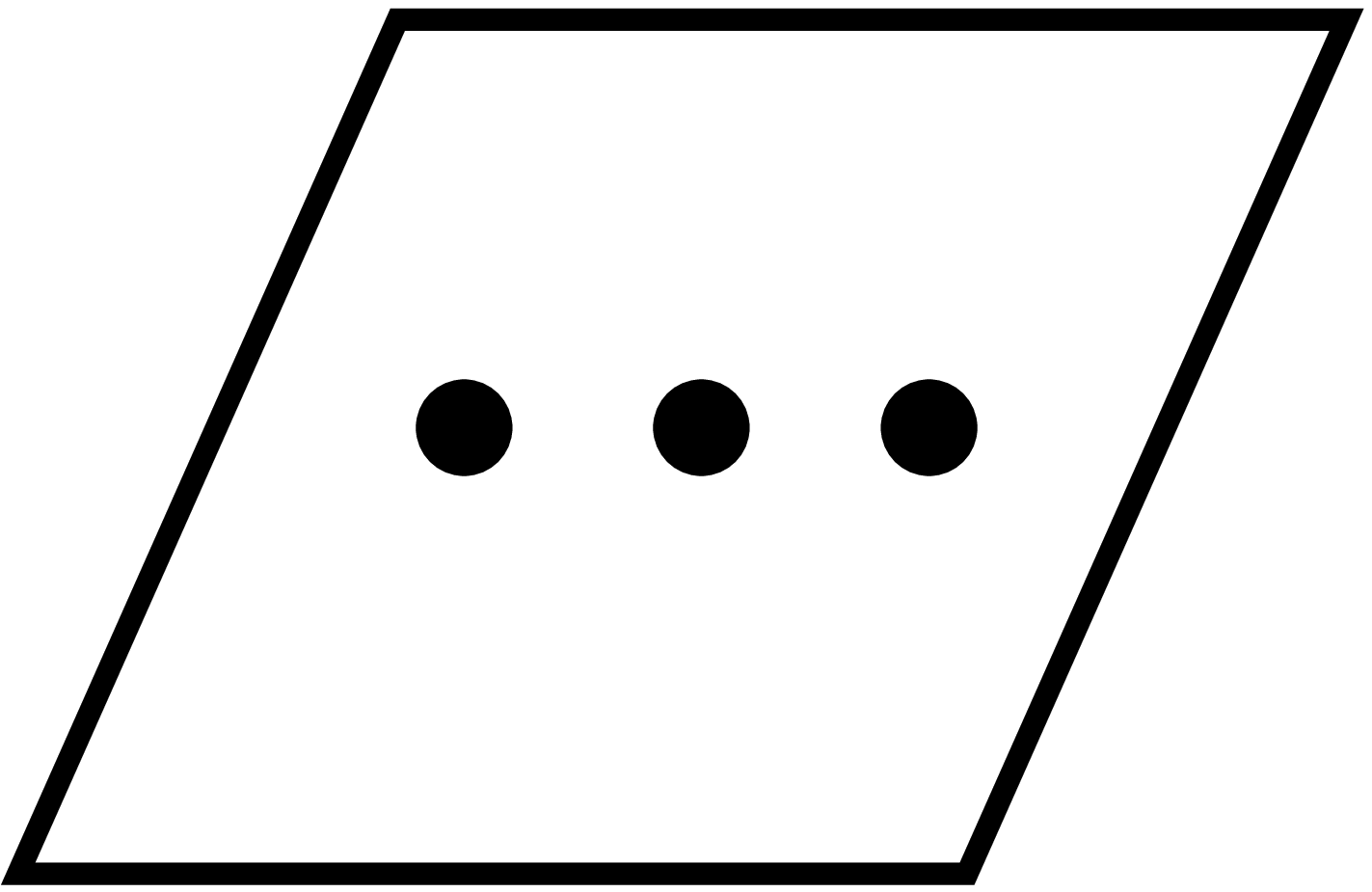}
= a
\figins{-7}{0.25}{plandd}
+ b
\figins{-7}{0.25}{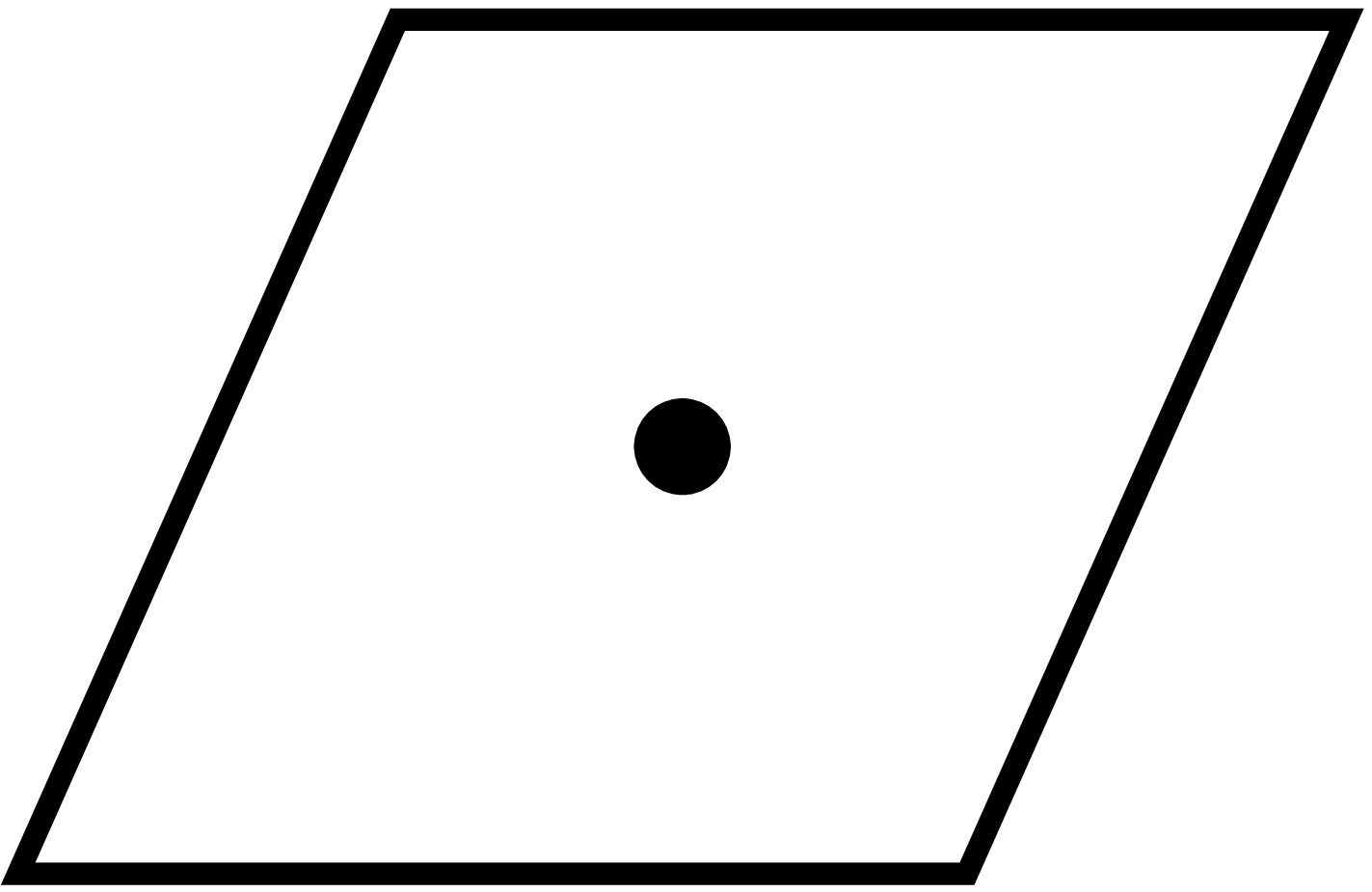}
+ c
\figins{-7}{0.25}{plan}
\tag{3D}
\\[1.5ex]\displaybreak[0]
%cutting neck
 \figwhins{-22}{0.65}{0.26}{cyl}=
-\figwhins{-22}{0.65}{0.26}{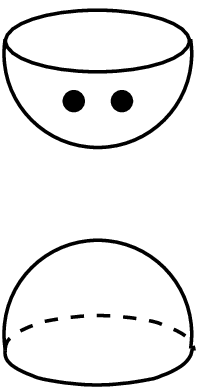}
-\figwhins{-22}{0.65}{0.26}{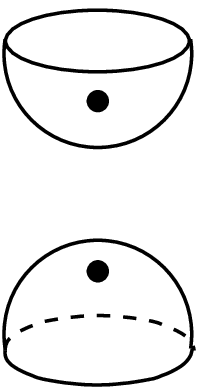}
-\figwhins{-22}{0.65}{0.26}{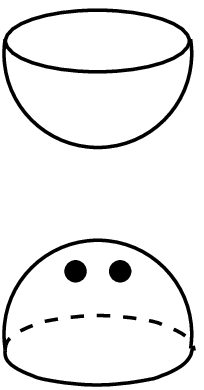}
+a 
\left(\
\figwhins{-22}{0.65}{0.26}{cnecka1}+
\figwhins{-22}{0.65}{0.26}{cnecka2}\
\right)
+b
\figwhins{-22}{0.65}{0.26}{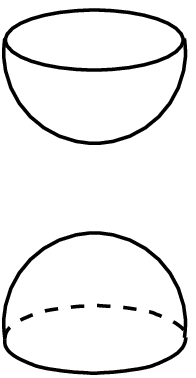}
\tag{CN}\label{eq:cn}
\\[1.5ex]\displaybreak[0]
%S-relation
\figins{-8}{0.3}{sundot}=
\figins{-8}{0.3}{sdot}=0,\quad
\figins{-8}{0.3}{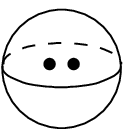}=-1
\tag{S}
\\[1.5ex]\displaybreak[0]
\labellist
\small\hair 2pt
\pinlabel $\alpha$ at 3 33
\pinlabel $\beta$ at -3 17
\pinlabel $\delta$ at 3 5
\endlabellist
\centering
\figins{-10}{0.4}{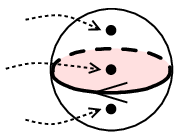} 
=\begin{cases}
\ \ 1 & (\alpha,\beta,\delta)=(1,2,0)\text{ or a cyclic permutation} \\ 
-1 & (\alpha,\beta,\delta)=(2,1,0)\text{ or a cyclic permutation} \\ 
\ \ 0 & \text{else}
\end{cases}
\tag{$\Theta$}\label{eq:theta}
\end{gather*}

The \emph{closure relation} says that any 
$\bC[a,b,c]$-linear combination of foams, all of which have the same boundary,  
is equal to zero if and only if any common way of closing these foams yields a 
$\bC[a,b,c]$-linear combination of closed foams whose evaluation is zero.

Using the relations $\ell$ one can prove the identities
below (for detailed proofs see \cite{MV1}).
\begin{align}
%%bamboo
\figwhins{-22}{0.65}{0.30}{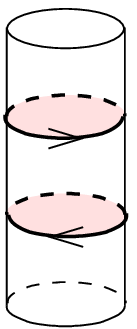}
& = -\ 
\figwhins{-22}{0.65}{0.30}{cneckb}
\tag{Bamboo}\label{eq:bamboo}
\\[1.2ex]\displaybreak[0]
%removing disk
\figwhins{-22}{0.65}{0.30}{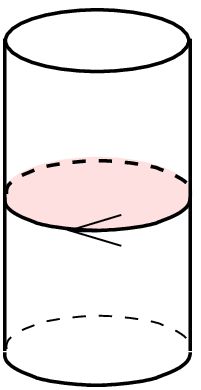}
& =\
\figwhins{-22}{0.65}{0.30}{cnecka1}-
\figwhins{-22}{0.65}{0.30}{cnecka2}
 \tag{RD}\label{eq:rd}
\\[1.2ex]\displaybreak[0]
%bubble
\figins{-12}{0.4}{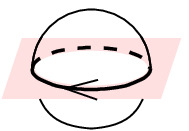} 
& =\ 0
\tag{Bubble}\label{eq:bubble}
\\[1.2ex]\displaybreak[0]
%digon removal
\figins{-20}{0.6}{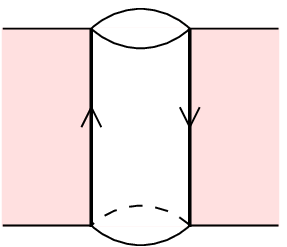}
& = 
\figins{-26}{0.75}{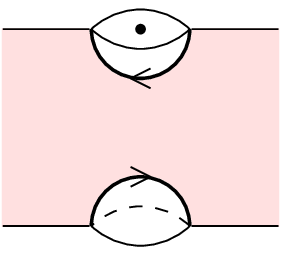}-
\figins{-26}{0.75}{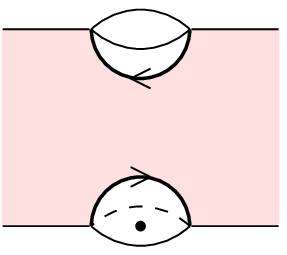}
\tag{DR}\label{eq:dr}
\\[1.2ex]\displaybreak[0]
%Square removal
\figins{-28}{0.8}{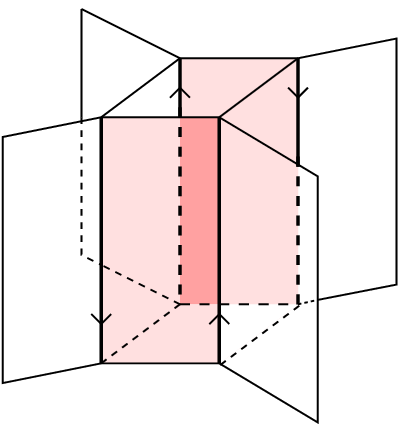}
&=
-\ \figins{-28}{0.8}{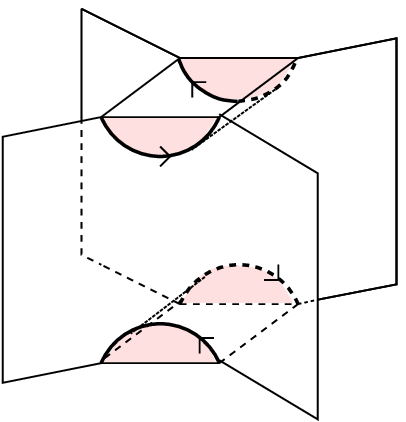}
-\figins{-28}{0.8}{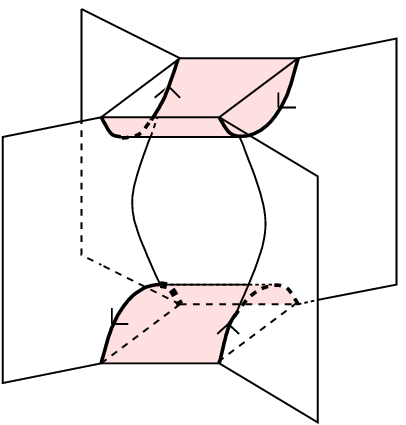}
\tag{SqR}\label{eq:sqr}
\end{align}

\begin{equation}\tag{Dot Migration}\label{eq:dotm}
\begin{split}
\figins{-22}{0.6}{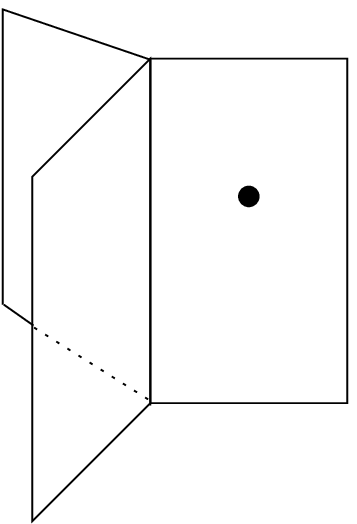}
\,+\,
\figins{-22}{0.6}{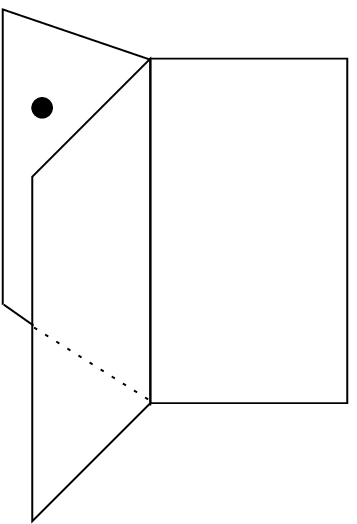}
\,+\,
\figins{-22}{0.6}{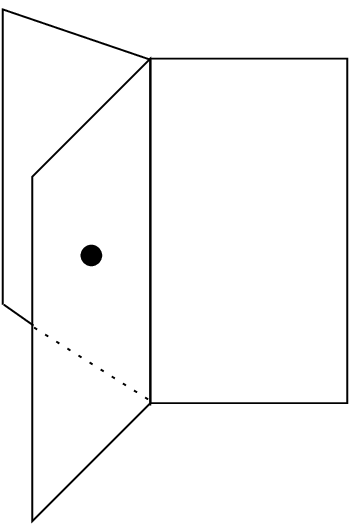}
\, &=\mspace{17mu} a\
\figins{-22}{0.6}{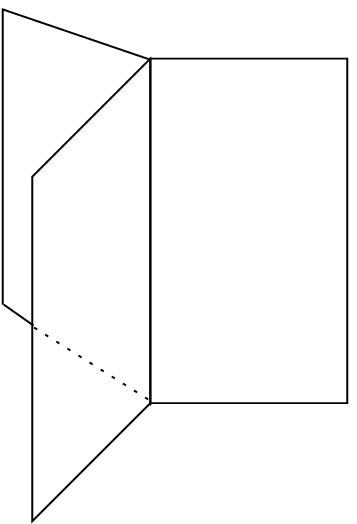} 
\\[1.2ex]
\figins{-22}{0.6}{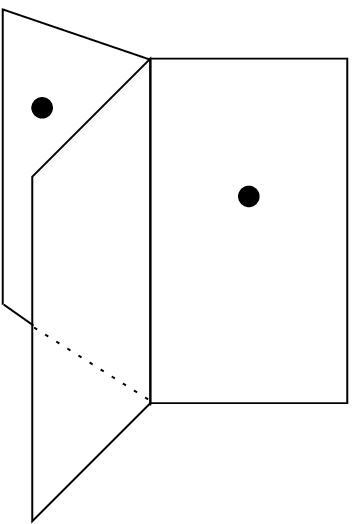}
\,+\,
\figins{-22}{0.6}{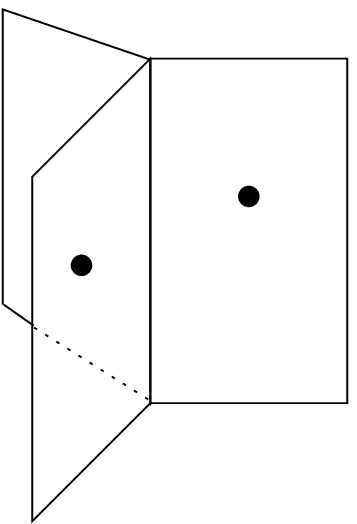}
\,+\,
\figins{-22}{0.6}{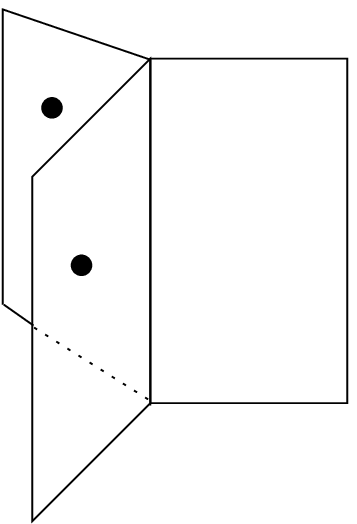}
\, &=\ -b\
\figins{-22}{0.6}{pdots000} 
\\[1.2ex]
\figins{-22}{0.6}{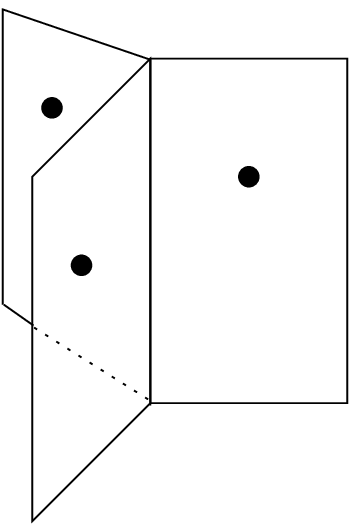}
\, &=\mspace{18.5mu} c\
\figins{-22}{0.6}{pdots000}
\end{split}
\end{equation}

In this paper we will work with open webs and open foams.
\begin{defn}
$\foamt$ is the category whose objects are webs $\Gamma$ inside a horizontal strip in 
$\bR^2$ bounded by the lines 
$y=0,1$ containing the boundary points of $\Gamma$ and whose morphisms are 
$\bC[a,b,c]$-linear combinations of foams inside that strip times the unit interval such 
that the vertical boundary of 
each foam is a set (possibly empty) of vertical lines.
\end{defn}

For example, the diagrams $1_n$ and $v_j$ are objects of $\foamt$:
\begin{equation*}
\labellist
\pinlabel $\dotsm$ at 125 71
\pinlabel $\dotsm$ at 564 71
\pinlabel $\dotsm$ at 853 71
\tiny\hair 2pt
\pinlabel $1$   at   6 -15
\pinlabel $2$   at  62 -15
\pinlabel $n$   at 176 -16
\pinlabel $1$   at 508 -15
\pinlabel $j$   at 675 -15
\pinlabel $j+1$ at 737 -16
\pinlabel $n$   at 906 -15
\endlabellist
1_n\ =\
\figins{-20}{0.6}{1n}
\mspace{100mu}
v_j\ =\ 
\figins{-20}{0.6}{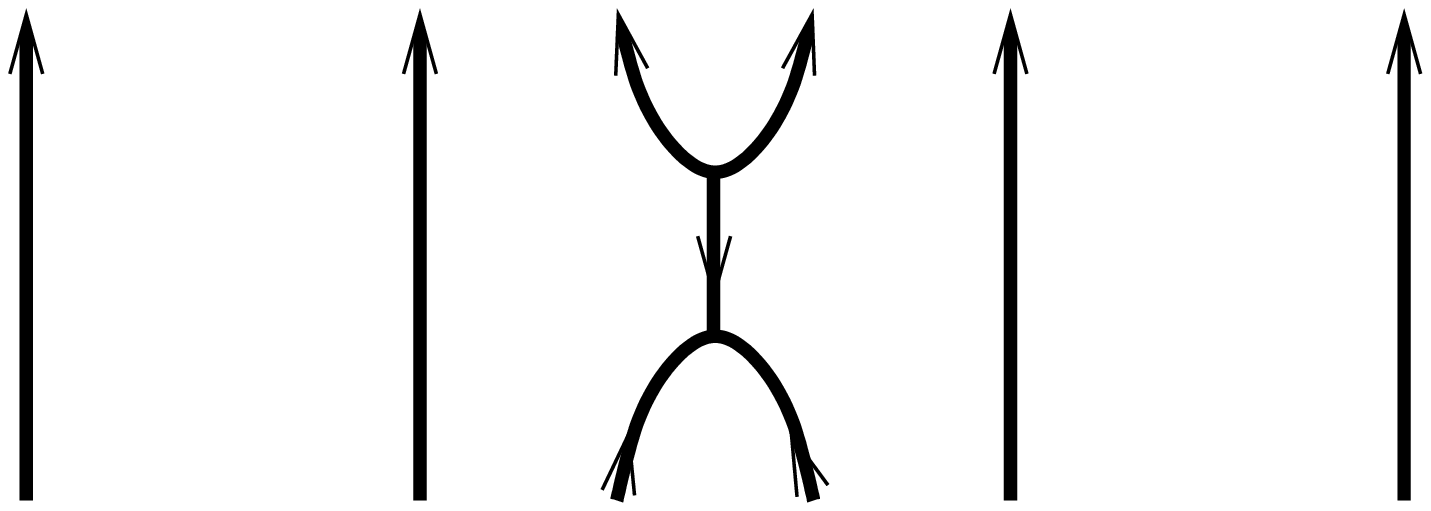}\vspace*{2ex}
\end{equation*}

The category $\foamt$ is additive and monoidal, with the monoidal structure given as in
$\cmw$.
The category $\foamt$ is also additive and graded. The $q$-grading in 
$\bC[a,b,c]$ is defined as 
$$q(1)=0,\quad q(a)=2,\quad q(b)=4,\quad q(c)=6$$
and the degree of a foam $f$ with $\vert\bullet\vert$ dots and $\vert b\vert$ vertical boundary components 
is given by
$$q(f)=-2\chi(f)+\chi(\partial f)+2\vert\bullet\vert+\vert b\vert ,$$
where $\chi$ denotes the Euler characteristic and $\partial f$ is the boundary of $f$.

%%%%%%%%%%%%%%%%%%%%%%%%%%%%%%%
\subsection{The functor $\Fmv$}
In this subsection we define a monoidal functor $\Fmv$ between the categories $\mathcal{SC}$ and $\foamt$.

\medskip

\n\emph{On objects:} 
$\Fmv$ sends the empty sequence to $1_n$ and the one-term sequence $(j)$ to $v_j$
with $\Fmv(jk)$ given by the vertical composite $v_{j}v_{k}$.

\medskip
\n\emph{On morphisms:}

\begin{itemize}
\item As before the empty diagram is sent to $n$ parallel vertical sheets:
\begin{equation*}
\labellist
\pinlabel $\dotsm$ at 185 120
\tiny\hair 2pt
\pinlabel $1$   at   0 -14
\pinlabel $2$   at  59 -14
\pinlabel $n-1$ at 219 -14
\pinlabel $n$   at 282 -14
\endlabellist
\emptyset \ \longmapsto \figins{-28}{1.0}{nsheets}\vspace*{2ex}
\end{equation*}
\item The vertical line colored $j$ is sent to the identity foam of $v_j$:
\begin{equation*}
\labellist
\tiny\hair 2pt
\pinlabel $j$   at -10  60
\pinlabel $j$   at 100 -65
\pinlabel $j+1$ at 165 -66
\endlabellist
\figins{-16}{0.5}{line}\ \
\longmapsto\
\figins{-32}{1.0}{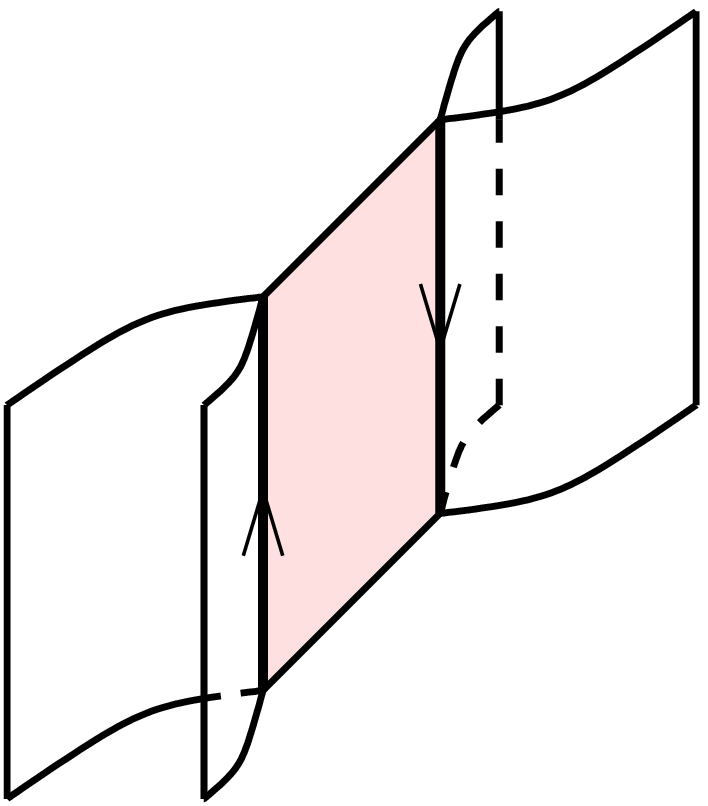}
\vspace*{2ex}
\end{equation*}

\item The \emph{StartDot} and \emph{EndDot} morphisms are sent to the zip and the unzip respectively:
\begin{equation*}
\labellist
\tiny\hair 2pt
\pinlabel $j$   at -10  60
\pinlabel $j$   at 100 -65
\pinlabel $j+1$ at 165 -66
\endlabellist
\figins{-16}{0.5}{startdot}
\longmapsto\
\figins{-32}{1.0}{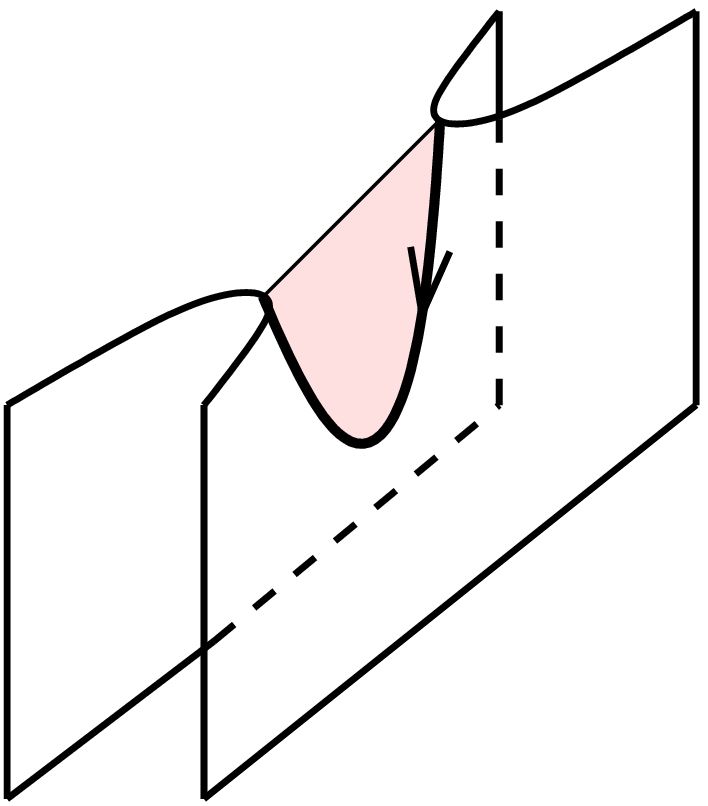}
\mspace{140mu}
\labellist
\tiny\hair 2pt
\pinlabel $j$   at -10  60
\pinlabel $j$   at 100 -65
\pinlabel $j+1$ at 165 -66
\endlabellist
\figins{-16}{0.5}{enddot}
\longmapsto\
\figins{-32}{1.0}{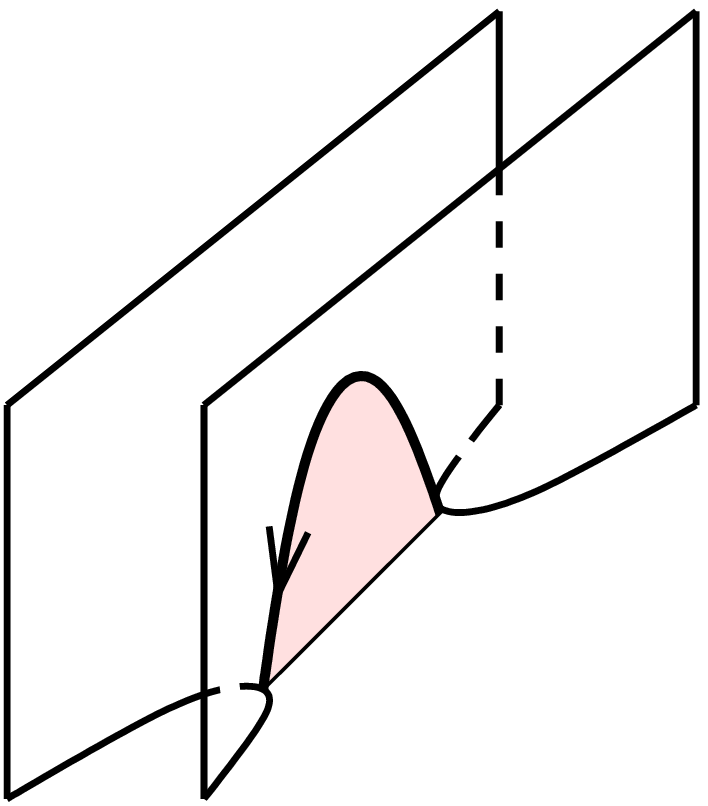}
\vspace*{2ex}
\end{equation*}

\item \emph{Merge} and \emph{Split} are sent to the digon creation and annihilation respectively:
\begin{equation*}
\labellist
\tiny\hair 2pt
\pinlabel $j$   at  45  95
\pinlabel $j$   at 205 -65
\pinlabel $j+1$ at 275 -66
\endlabellist
\figins{-16}{0.6}{merge}
\longmapsto\
\figins{-32}{1.1}{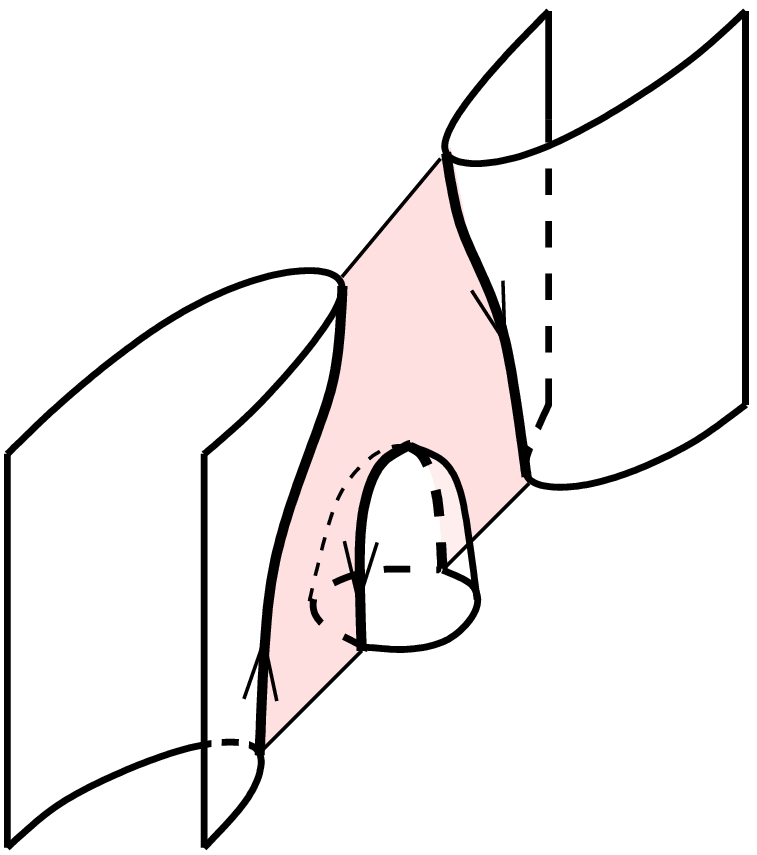}
%\vspace*{2ex}
%\end{equation*}
\mspace{80mu}
%\begin{equation*}
\labellist
\tiny\hair 2pt
\pinlabel $j$   at  45  45
\pinlabel $j$   at 205 -65  
\pinlabel $j+1$ at 275 -66
\endlabellist
\figins{-16}{0.6}{split}
\longmapsto \
\figins{-32}{1.1}{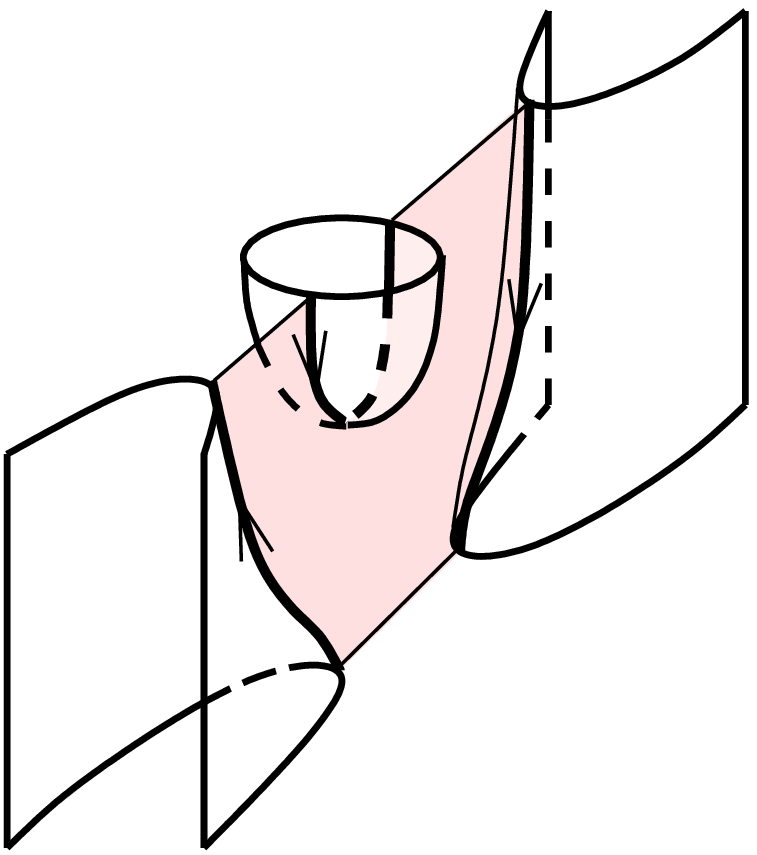}
\vspace*{2ex}
\end{equation*}

\item The \emph{4-valent vertex} with distant colors. For $j+1<k$ we have:
\begin{equation*}
\labellist
\tiny\hair 2pt
\pinlabel $k$   at  -5 -12
\pinlabel $j$   at 128 -10
\pinlabel $j$      at 205 -61
\pinlabel $j+1$    at 270 -62
\pinlabel $k$      at 385 -61
\pinlabel $k+1$    at 445 -62
\pinlabel $\dotsm$ at 335 -62
\endlabellist
\figins{-16}{0.6}{4vert}
\longmapsto\
\figins{-32}{1.1}{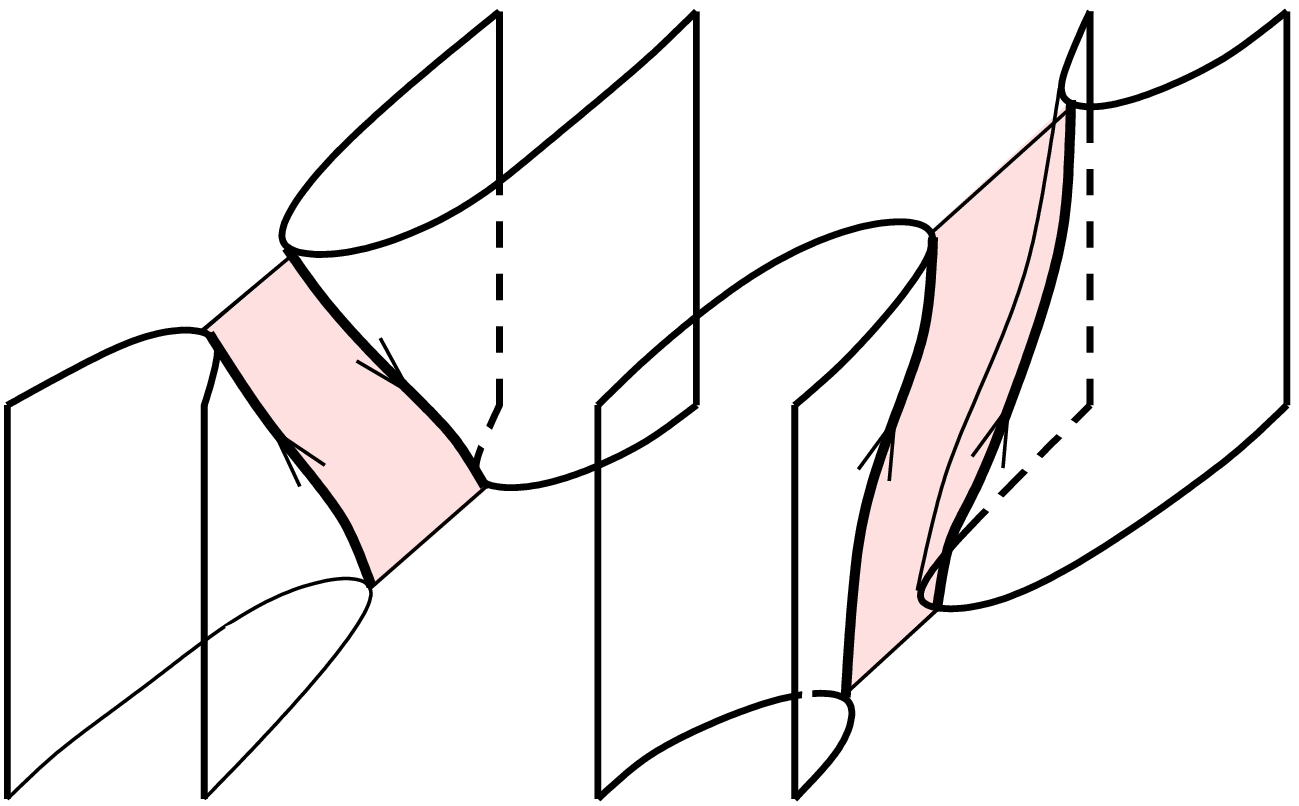}
\vspace*{2ex}
\end{equation*}
The case $j>k+1$ is given by reflection around a horizontal plane.

\item For the \emph{6-valent vertex} we have:
\begin{equation*}
\labellist
\tiny\hair 2pt
\pinlabel $j+1$ at -5 -16  \pinlabel $j$ at 67 -15
\pinlabel $j$ at 243 -115  \pinlabel $j+1$ at 320 -116  \pinlabel $j+2$ at 395 -116
\endlabellist
\figins{-18}{0.6}{6vertu}
\longmapsto\ 
-\
\figins{-52}{1.8}{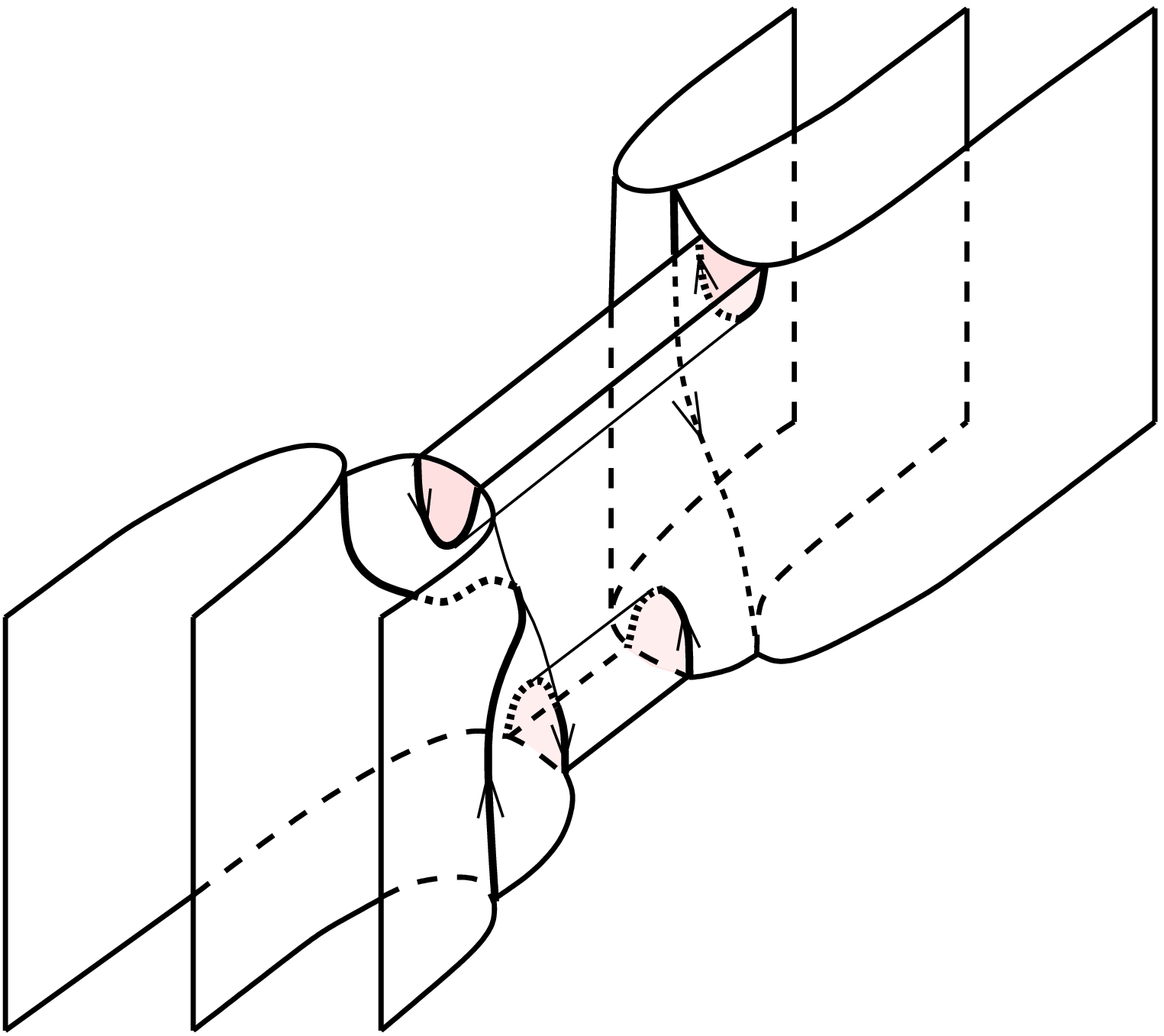}
\vspace*{2ex}
\end{equation*}
\end{itemize}

The case with the colors switched is given by reflection in a vertical plane.

Notice that $\Fmv$ respects the gradings of the morphisms.

\begin{prop}
$\Fmv$ is a monoidal functor.
\end{prop}
%\proof
\begin{proof}
The assignement given by $\Fmv$ clearly respects the monoidal structures of 
$\mathcal{SC}_1$ and $\foamt$. 
To prove it is a monoidal functor we need only to show that it is actually a functor, 
that is, it respects the 
relations~\eqref{eq:adj} to~\eqref{eq:dumbdumbsquare} of Section~\ref{sec:soergel}.

%%%%%%%%%%%%%%%%%%%%%%%%%%%%%%%
\medskip
\n\emph{''Isotopy relations'':} 
Relations~\eqref{eq:adj} to~\eqref{eq:v6rot} correspond to isotopies of their images under $\Fmv$ and we leave 
its check to the reader.

%%%%%%%%%%%%%%%%%%%%%%%%%%%%%
\medskip
\n\emph{One color relations:}
Relation~\eqref{eq:dumbrot} is straightforward and left to the reader. For relation~\eqref{eq:lollipop} we have
\begin{equation*}
\labellist
\tiny\hair 2pt
\pinlabel $j$ at 26 45            
\pinlabel $j$ at 240 -98 \pinlabel $j+1$ at 345 -99
\endlabellist
\Fmv\Biggl(\
\figins{-16}{0.5}{lollipop-u}\
\Biggr)
=\
\figins{-32}{1.0}{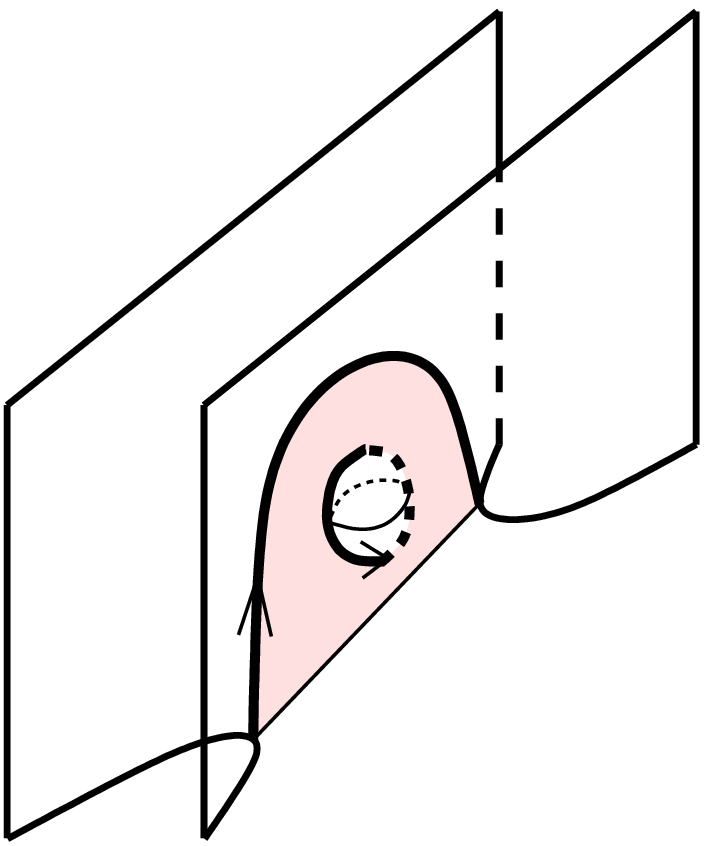}\
=\
0\ ,
\vspace*{2ex}
\end{equation*}
the last equality following from the~\eqref{eq:bubble} relation.

For relation~\eqref{eq:deltam} we have
\begin{equation*}
\labellist
\tiny\hair 2pt
\pinlabel $j$ at -2 22    \pinlabel $j$ at -2 130  
\pinlabel $j$ at 130 -70 \pinlabel $j+1$ at 205 -71
\pinlabel $j$ at 405 -70 \pinlabel $j+1$ at 485 -71
\pinlabel $j$ at 725 -70 \pinlabel $j+1$ at 805 -71
\endlabellist
\Fmv\Biggl(\
\figins{-17}{0.55}{matches-ud}\
\Biggr)
=\
\figins{-32}{1.0}{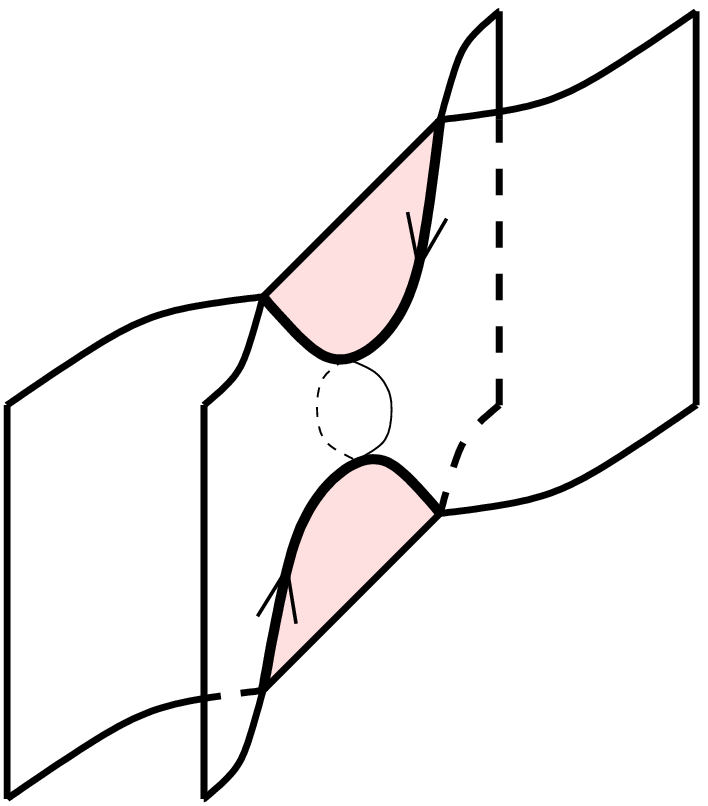}\
=\
\figins{-32}{1.0}{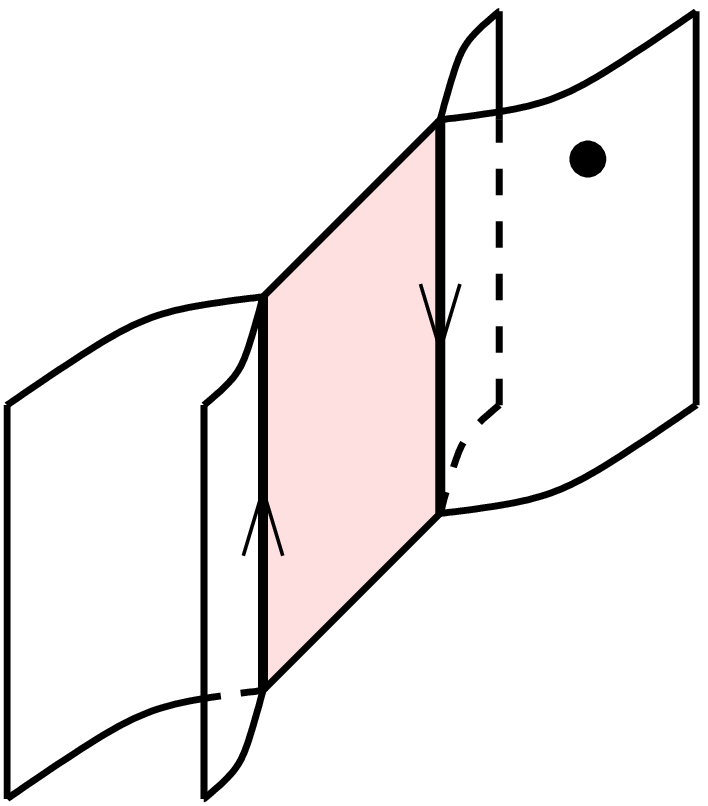}\
-\
\figins{-32}{1.0}{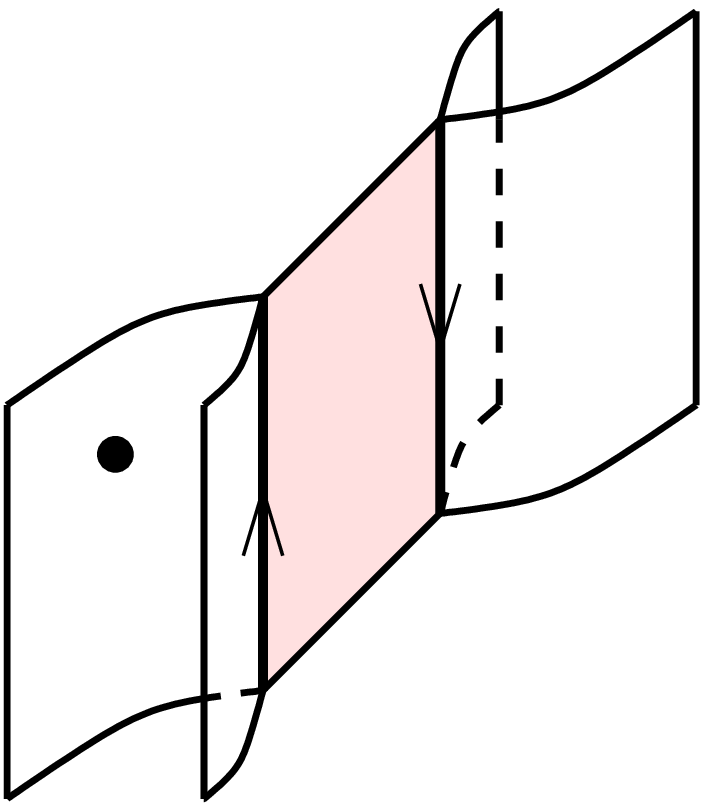}
\vspace*{2ex}
\end{equation*}

\n where the second equality follows from the~\eqref{eq:dr} relation. We also have
\begin{equation*}
\labellist
\tiny\hair 2pt
\pinlabel $j$ at -13   16
\pinlabel $j$ at 140 -130   \pinlabel $j+1$ at 220 -131
\pinlabel $j$ at 440 -130   \pinlabel $j+1$ at 590 -131
\pinlabel $j$ at 785 -130   \pinlabel $j+1$ at 940 -131
\endlabellist
\Fmv\biggl(\
\figins{-6}{0.25}{startenddot}\
\biggr)\
=\
\figins{-32}{1.0}{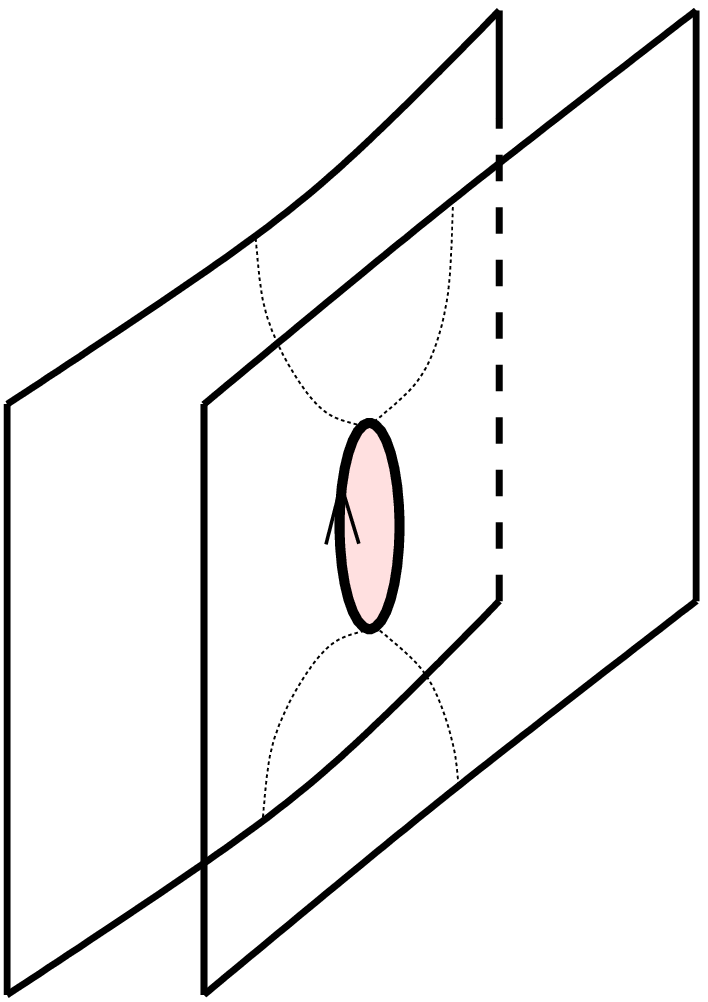}\
=\
\figins{-32}{1.0}{sheet-dot-l}\
-\
\figins{-32}{1.0}{dot-sheet-l}\ ,
\vspace*{2ex}
\end{equation*}
which is given by~\eqref{eq:rd}.
Using~\eqref{eq:dotm} one obtain
\begin{align}
\label{eq:edge-dots3}
\labellist
\tiny\hair 2pt
\pinlabel $j$ at -10   60  \pinlabel $j$ at    60 -265
\pinlabel $j$ at 190  -78  \pinlabel $j+1$ at 265  -79
\pinlabel $j$ at 510  -78  \pinlabel $j+1$ at 585  -79
\pinlabel $j$ at 860  -78  \pinlabel $j+1$ at 935  -79
\pinlabel $j$ at 220 -400  \pinlabel $j+1$ at 300 -401
\pinlabel $j$ at 540 -400  \pinlabel $j+1$ at 620 -401
\pinlabel $j$ at 890 -400  \pinlabel $j+1$ at 970 -401
\endlabellist
\Fmv\Biggl(\
\figins{-17}{0.55}{edge-startenddot}\
\Biggr)
&=
\ 2\
\figins{-32}{1.0}{idvi-dot}\
+\
\figins{-32}{1.0}{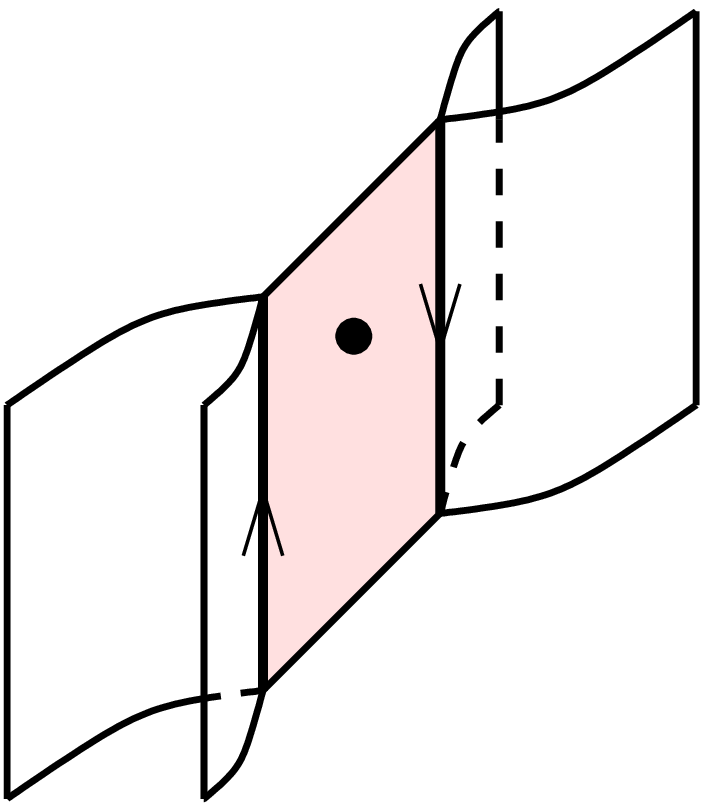}\
+\
a\
\figins{-32}{1.0}{idvi}
\\[1.5ex]\displaybreak[0]
\label{eq:dots-edge3}
\Fmv\Biggl(\
\figins{-17}{0.55}{startenddot-edge}\
\Biggr)
&=
\ -2\
\figins{-32}{1.0}{dot-idvi}\
-\
\figins{-32}{1.0}{idvid}\
-\
a\
\figins{-32}{1.0}{idvi}
\end{align}
and therefore, we have that
\begin{equation*}
\Fmv
\Biggl(
\figins{-17}{0.55}{startenddot-edge}\
\Biggr)
+\
\Fmv
\Biggl(\
\figins{-17}{0.55}{edge-startenddot}\
\Biggr)
=\ 2\ 
\Fmv
\Biggl(\
\figins{-17}{0.55}{matches-ud}\
\Biggr).
\end{equation*}

%%%%%%%%%%%%%%%%%%%%%%%%%%%%
\medskip
\n\emph{Two distant colors:} 
Relations~\eqref{eq:reid2dist} to~\eqref{eq:slide3v} correspond to isotopies of the foams involved and are 
straightforward to check.

%%%%%%%%%%%%%%%%%%%%%%%%%
\medskip
\n\emph{Adjacent colors:}
We prove the case where ''blue`` corresponds to $j$ and ''red`` corresponds to $j+1$. 
The relations with colors reversed are proved the same way.
To prove relation~\eqref{eq:dot6v} we first notice that
\begin{align*}
\labellist
\tiny\hair 2pt
\pinlabel $j$   at 275 -121
\pinlabel $j+1$ at 357 -122
\pinlabel $j+2$ at 461 -122
\endlabellist
\Fmv\Biggl(\
\figins{-12}{0.4}{capcupdot}\
\Biggr)\
& =\  \
\figins{-34}{1.2}{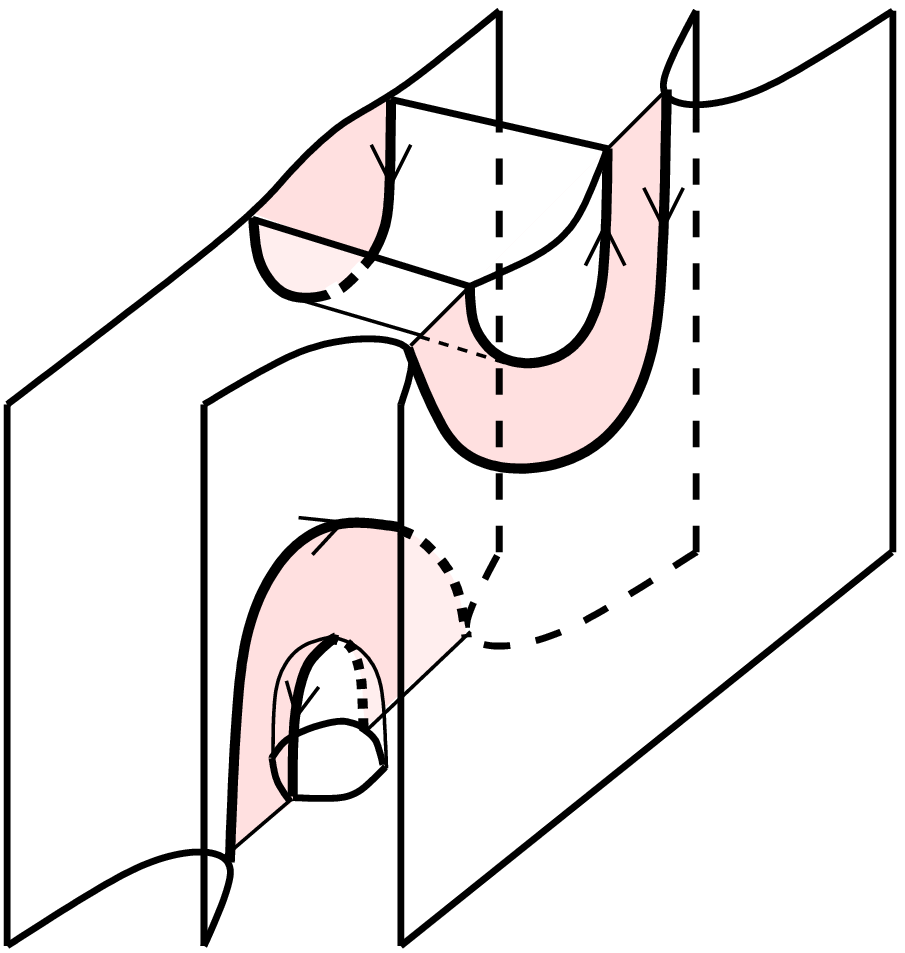}
\\ \displaybreak[0]
\intertext{and}
\labellist
\tiny\hair 2pt
\pinlabel $j$   at 275 -118
\pinlabel $j+1$ at 357 -119
\pinlabel $j+2$ at 461 -119
\endlabellist
\Fmv\Biggl(\
\figins{-12}{0.4}{mergedots}\
\Biggr)\
& =\  \
\figins{-34}{1.2}{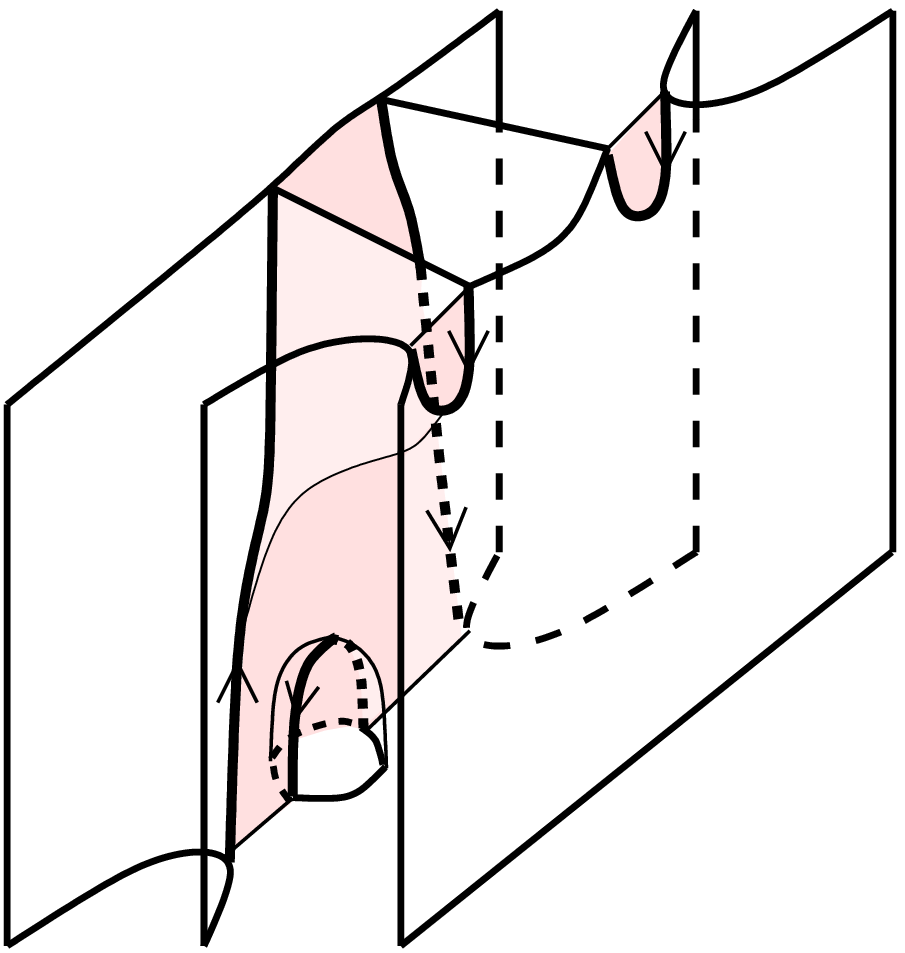}.
\vspace*{2ex}
\end{align*}

We also have an isotopy equivalence
\begin{equation*}
\labellist
\tiny\hair 2pt
\pinlabel $j$   at 315 -121
\pinlabel $j+1$ at 397 -122
\pinlabel $j+2$ at 501 -122
\endlabellist
\Fmv\Biggl(\
\figins{-12}{0.4}{6vertdotd}\
\Biggr)\
\cong\
-\
\figins{-34}{1.2}{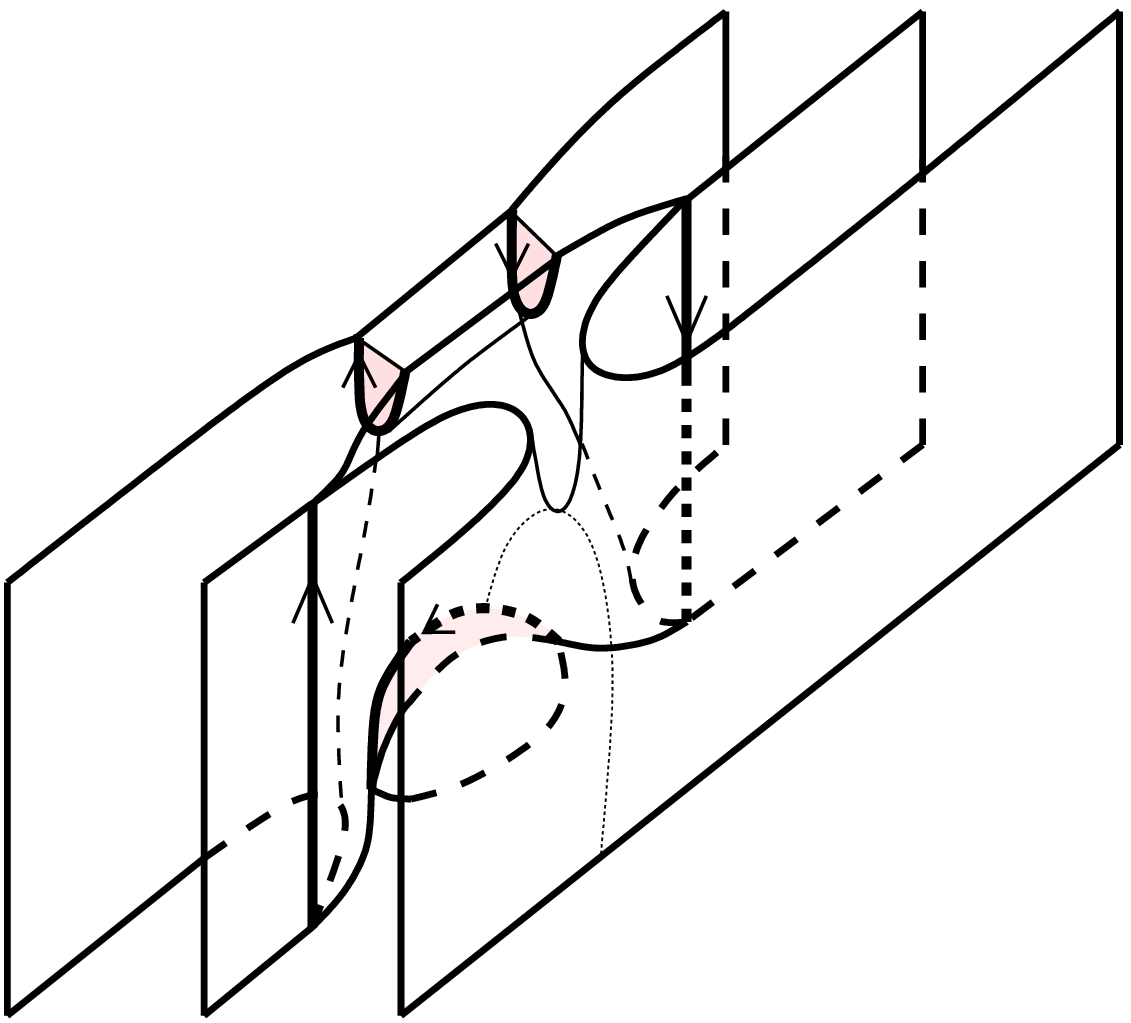}
\vspace*{2ex}
\end{equation*}
which in turn is isotopy equivalent to the foam obtained by putting
\begin{equation*}
\figins{-34}{1.2}{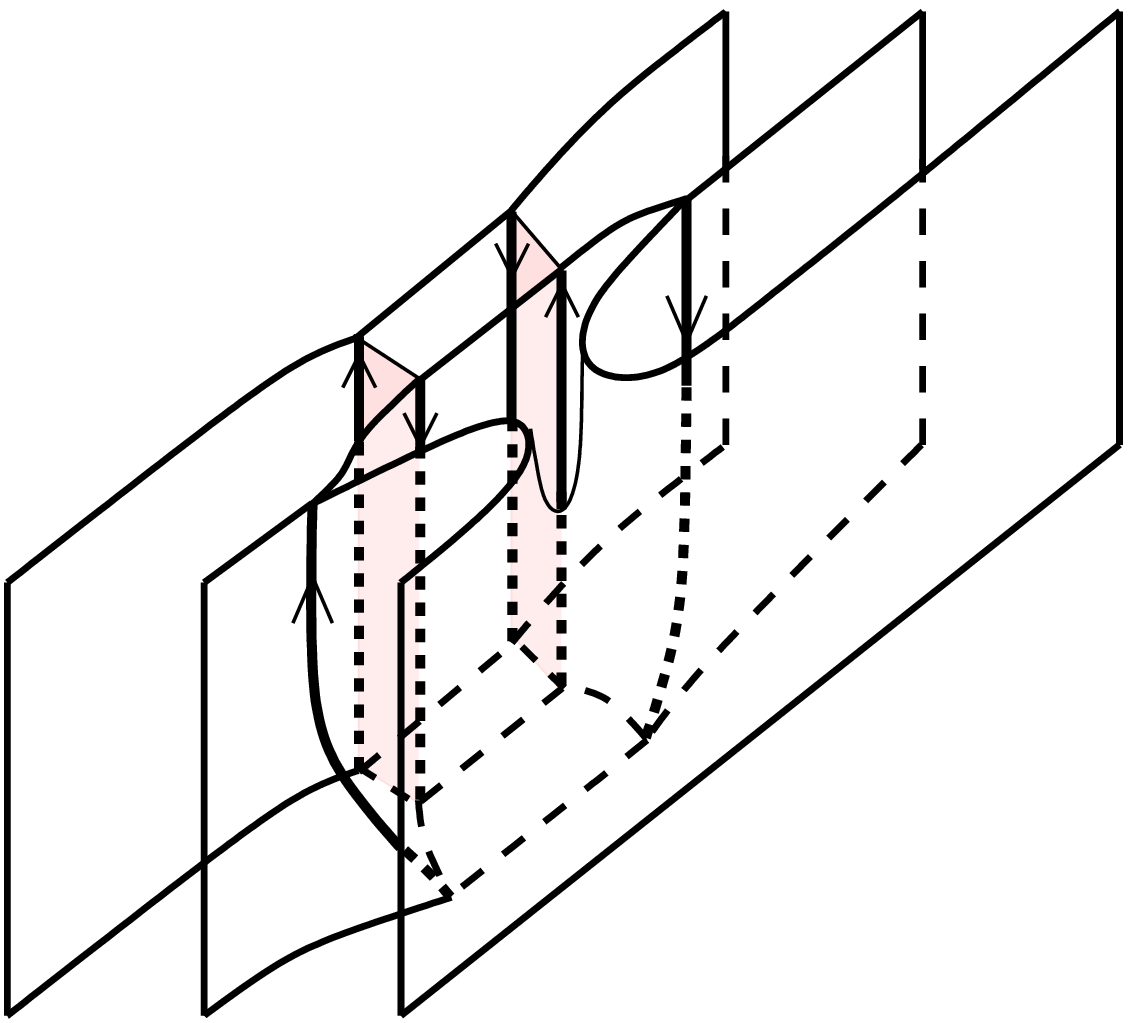}\
\mspace{12mu}\text{ on top of }\mspace{12mu}-\
\figins{-34}{1.2}{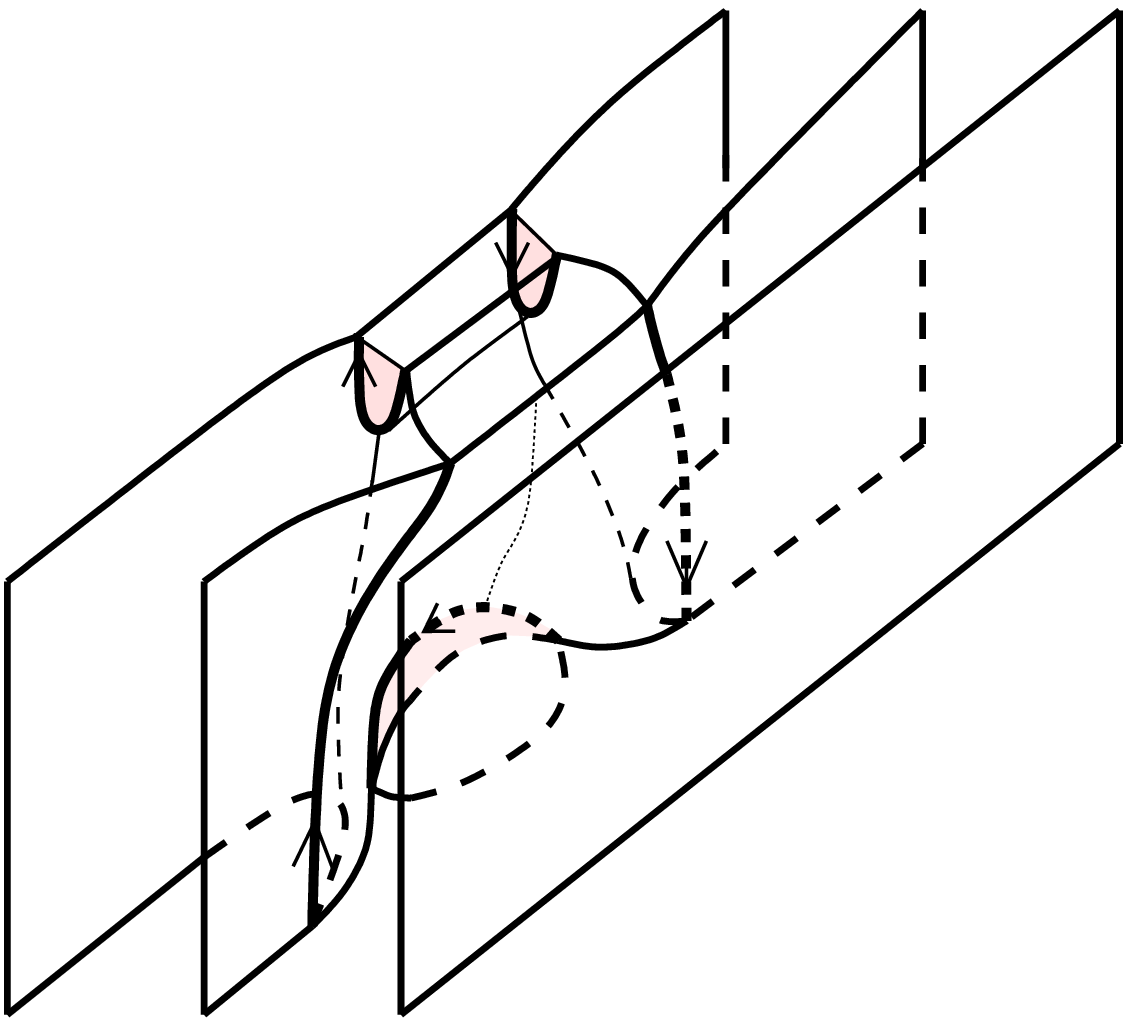}\ .
\vspace*{2ex}
\end{equation*}

The common boundary of these two foams contains two squares. Putting~\eqref{eq:sqr} on the square on the right glued 
with the identity foam everywhere else gives two terms, one isotopic to $\Fmv\bigl(\figins{-5}{0.2}{mergedots}\bigr)$ 
and the other isotopic to $\Fmv\bigl(\figins{-5}{0.2}{capcupdot}\bigr)$.

We now prove relation~\eqref{eq:reid3}. We have
\begin{equation*}
\labellist
\tiny\hair 2pt
\pinlabel $j$ at 265 20  \pinlabel $j+1$ at 280 -25  \pinlabel $j+2$ at 320 -62
\endlabellist
\Fmv
\left(\
\figins{-26}{0.8}{id-r3}\
\right)
\cong\
\figins{-36}{1.4}{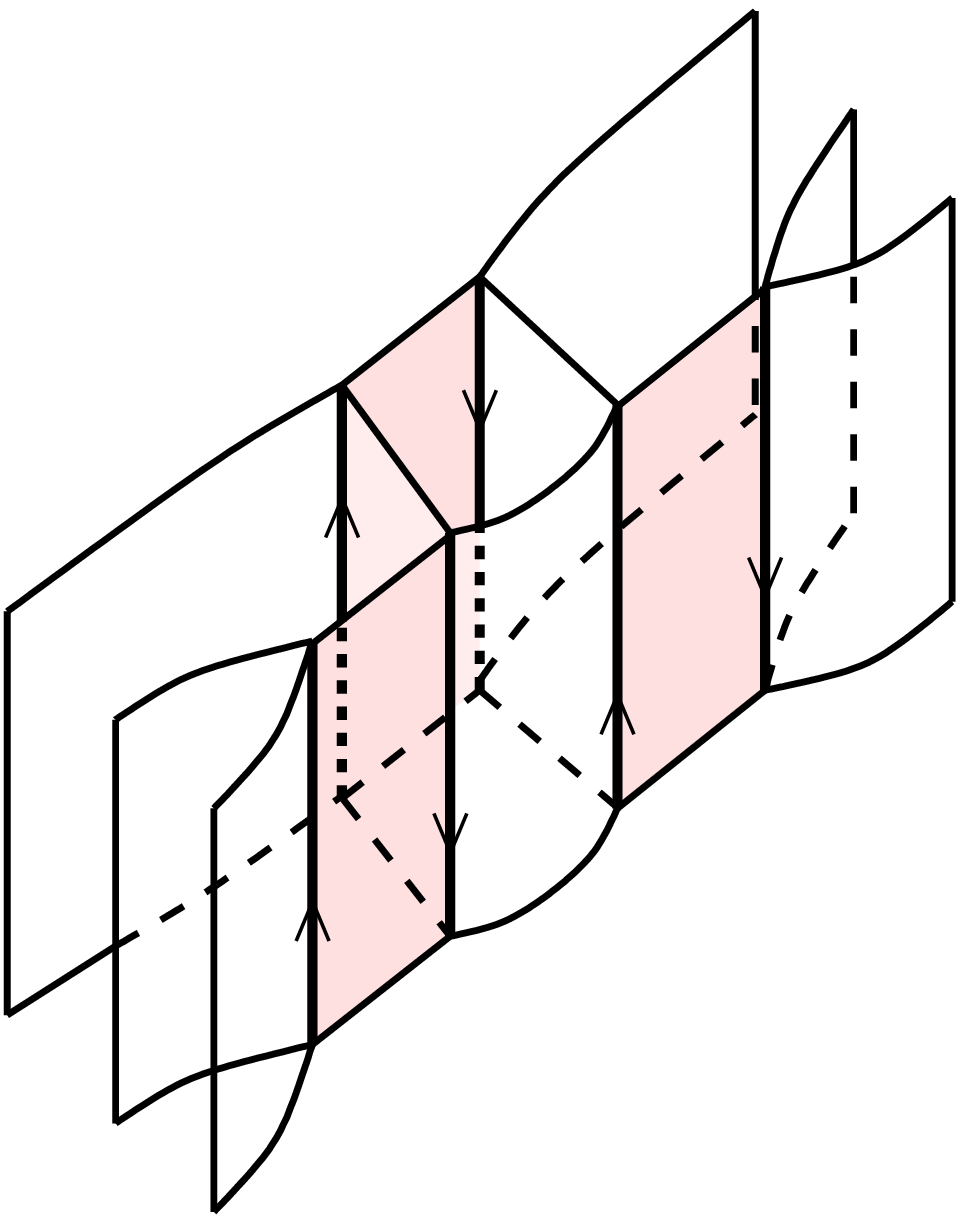}\ .
\end{equation*}

Applying~\eqref{eq:sqr} to the middle square we obtain two terms. One is isotopic to 
$-\Fmv\bigl(\figins{-6}{0.22}{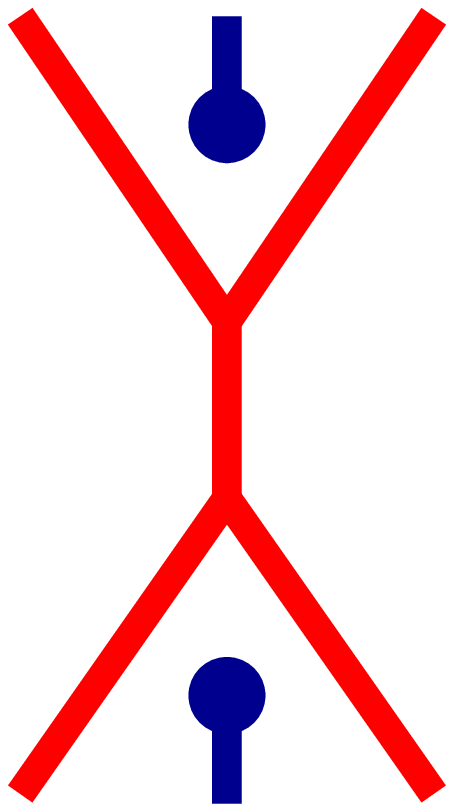}\bigr)$ 
and the other gives 
$\Fmv\bigl(\figins{-6}{0.22}{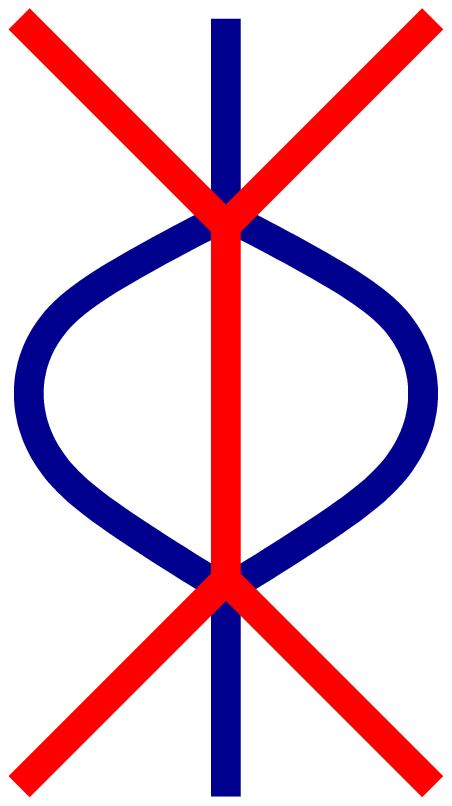}\bigr)$ after using the~\eqref{eq:bamboo} relation.

We now prove relation~\eqref{eq:dumbsq} in the form
\begin{equation*}
\figins{-26}{0.75}{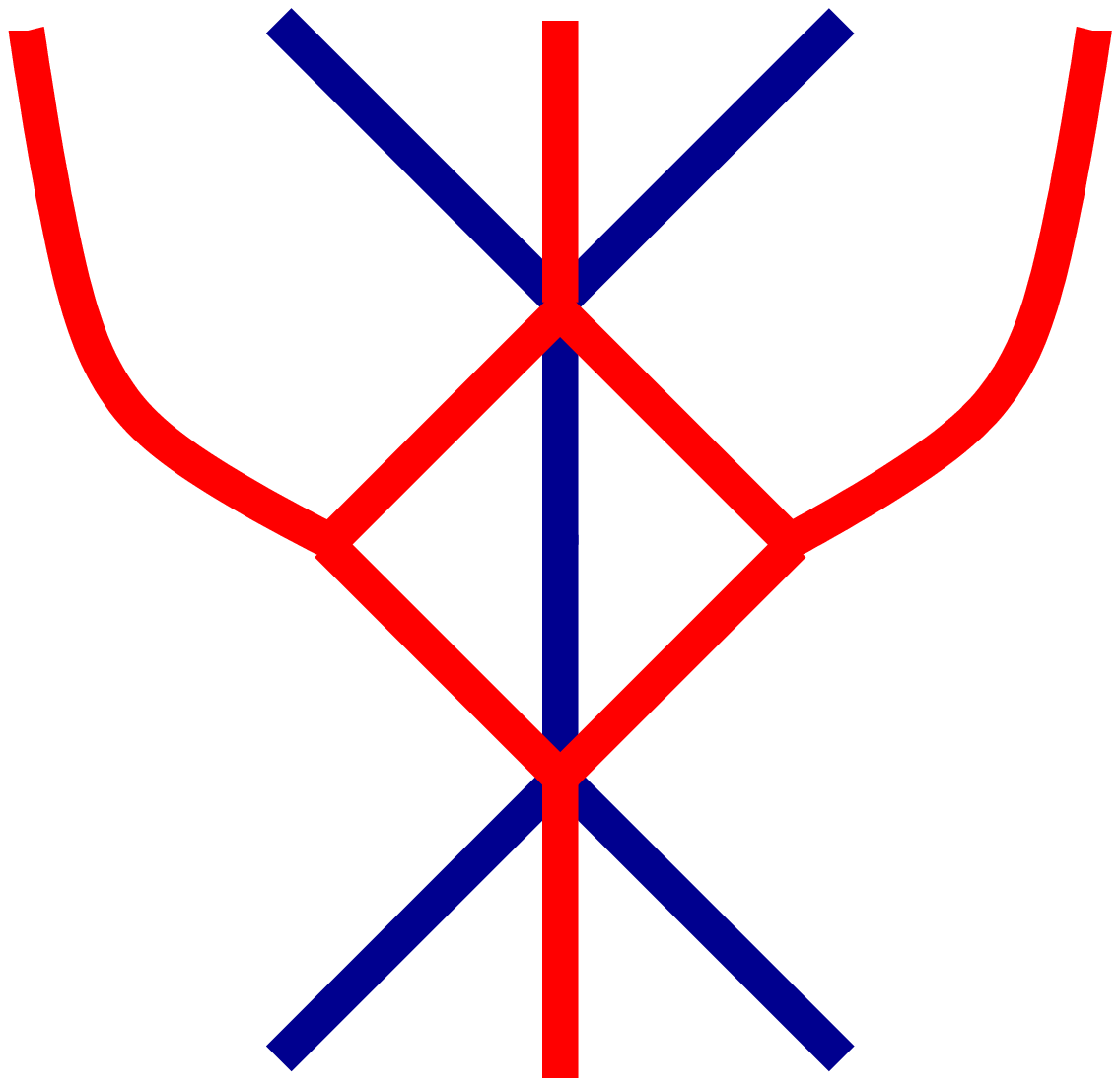}\
=\
\figins{-26}{0.75}{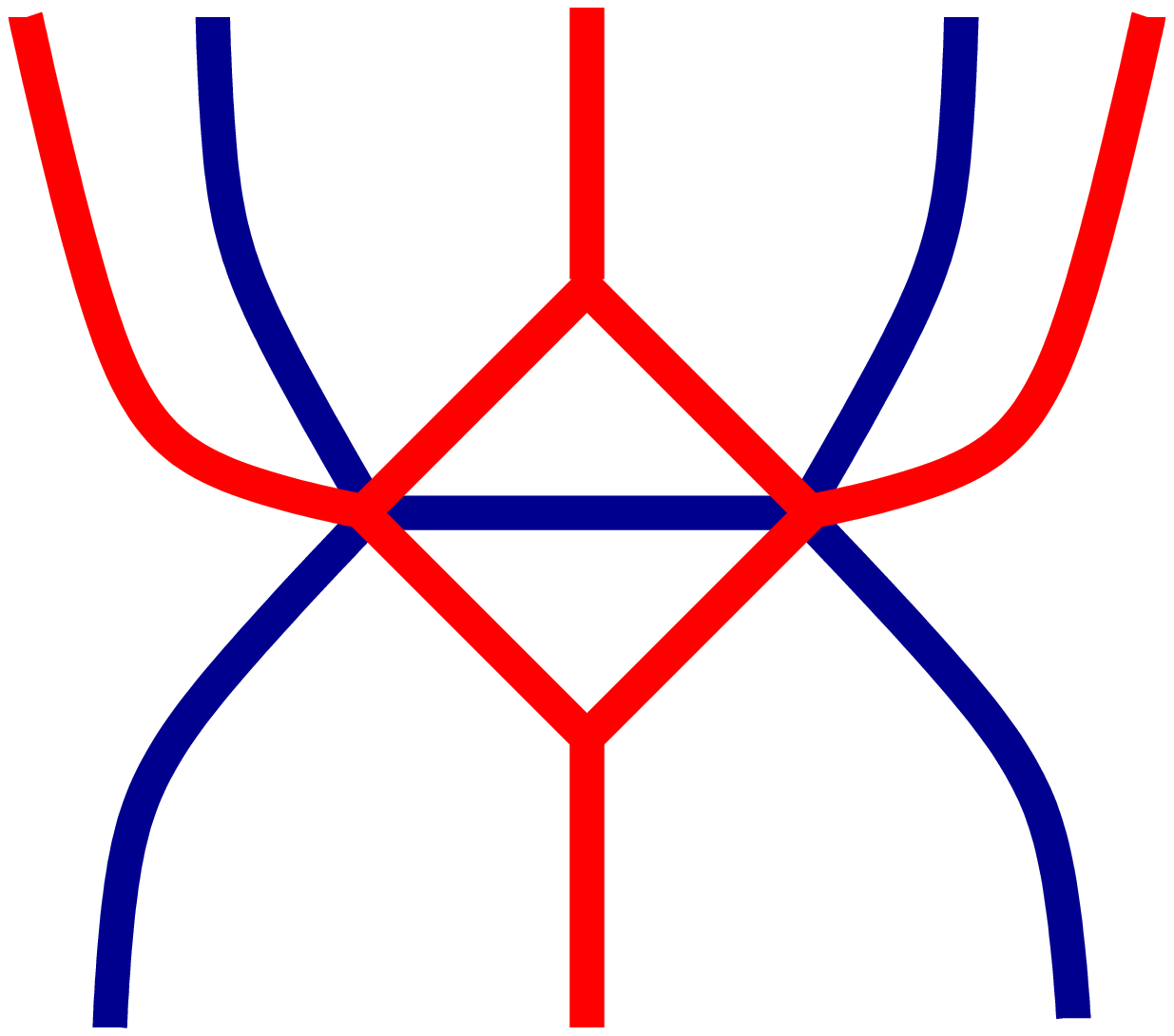}\ .
\end{equation*}

The images of the l.h.s. and r.h.s. under $\Fmv$ are isotopic to
\begin{equation*}
\labellist
\tiny\hair 2pt
\pinlabel $j$   at -20  50 
\pinlabel $j+1$ at  -5  12  
\pinlabel $j+2$ at  40 -22
\endlabellist
\figins{-50}{1.6}{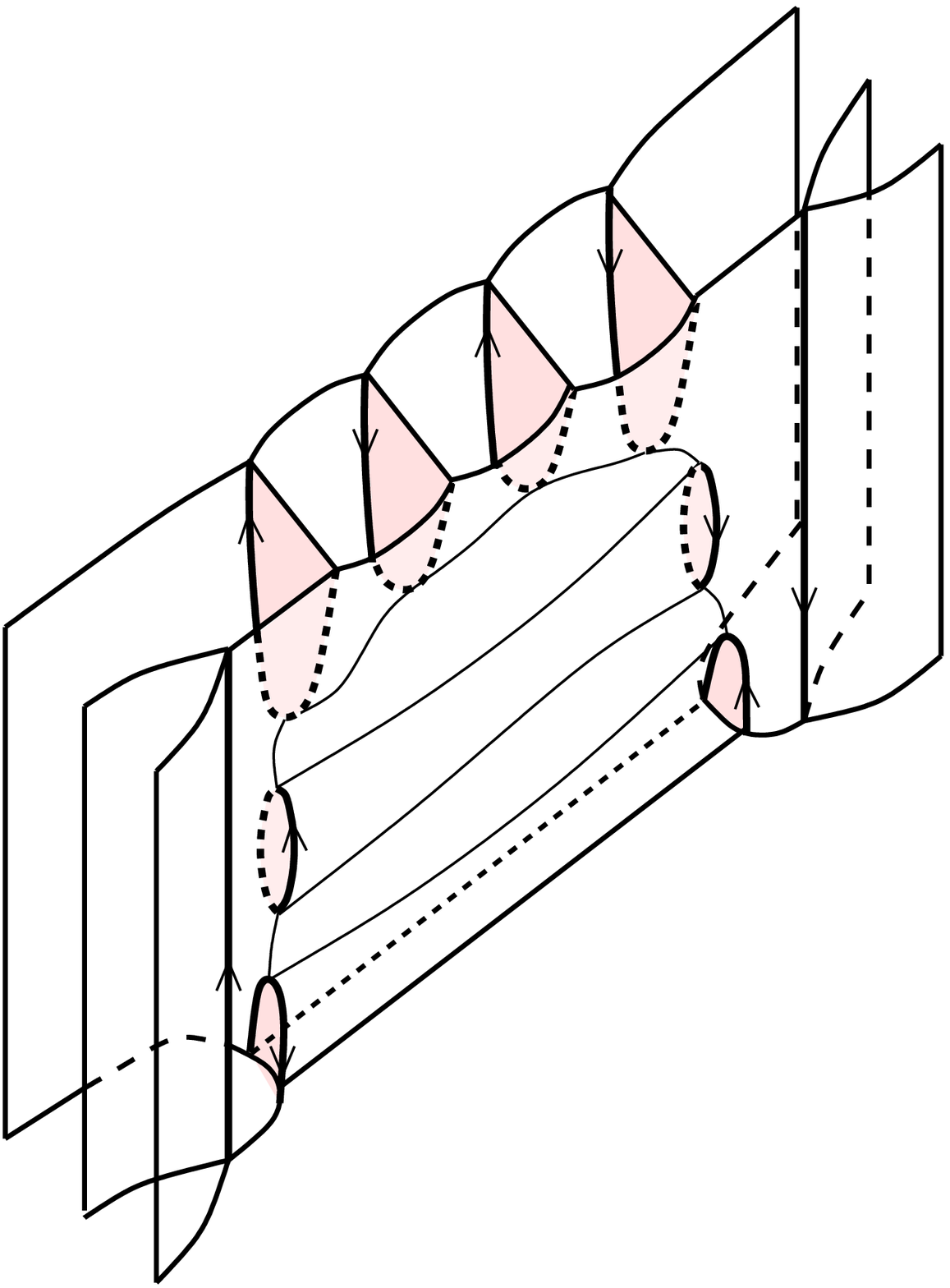}
\mspace{80mu}\text{and}\mspace{80mu}
\labellist
\tiny\hair 2pt
\pinlabel $j$   at -20  50 
\pinlabel $j+1$ at  -5  12  
\pinlabel $j+2$ at  40 -22
\endlabellist
\figins{-50}{1.6}{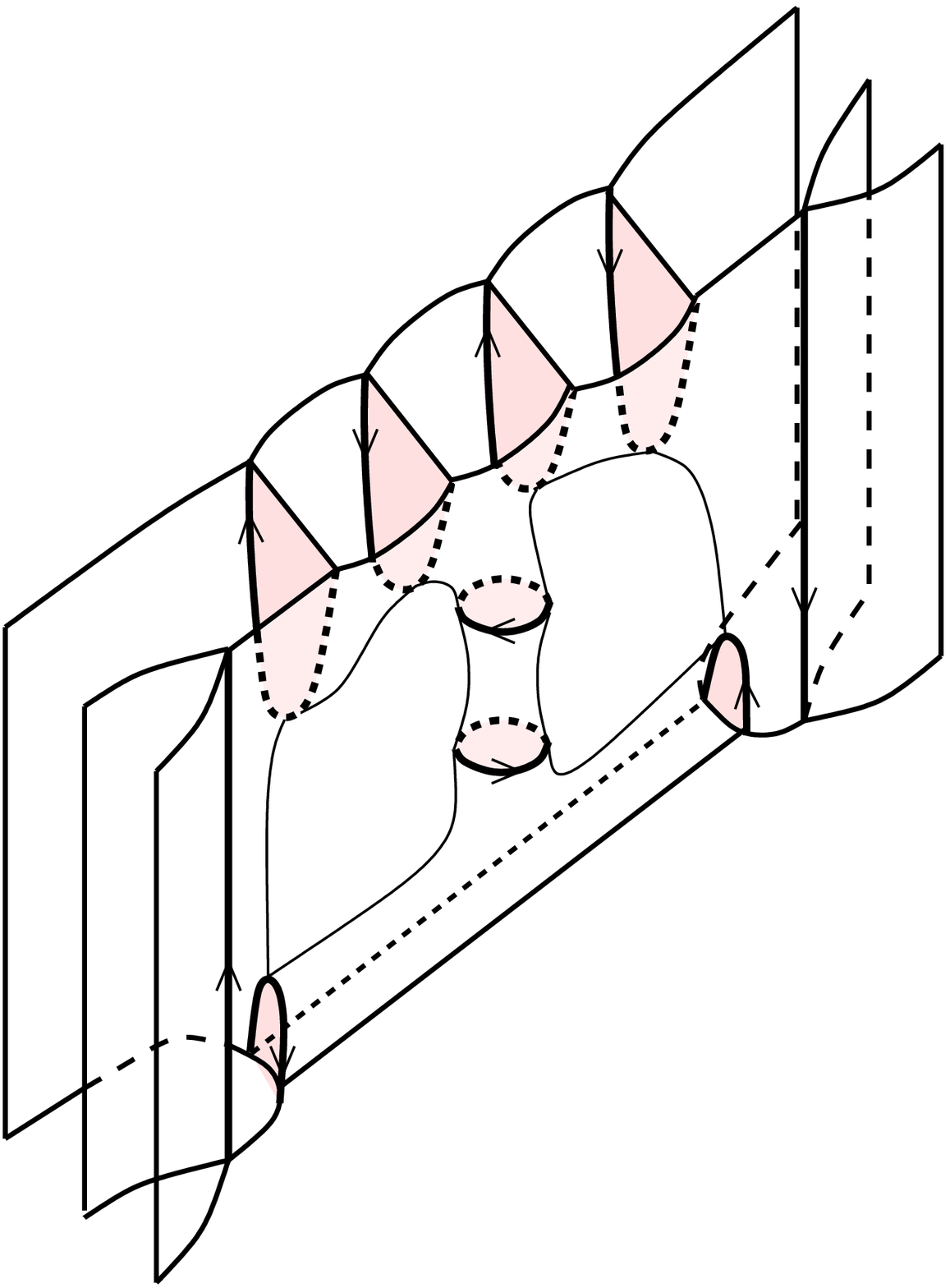}
\vspace*{2ex}
\end{equation*}
respectively, and both give the same foam after applying the~\eqref{eq:bamboo} relation.

Relation~\eqref{eq:slidenext} follows from a straightforward computation and is left to the reader.

%%%%%%%%%%%%%%%%%%%%%%%%%%%%%%%%%%%%%%%%%%
\medskip
\n\emph{Relations involving three colors:} 
Relations~\eqref{eq:slide6v} and~\eqref{eq:slide4v} follow from isotopies of the cobordisms involved.

We prove relation~\eqref{eq:dumbdumbsquare} in the form
\begin{equation}
\label{eq:dumbdumbsquare-b}
\figins{-26}{0.75}{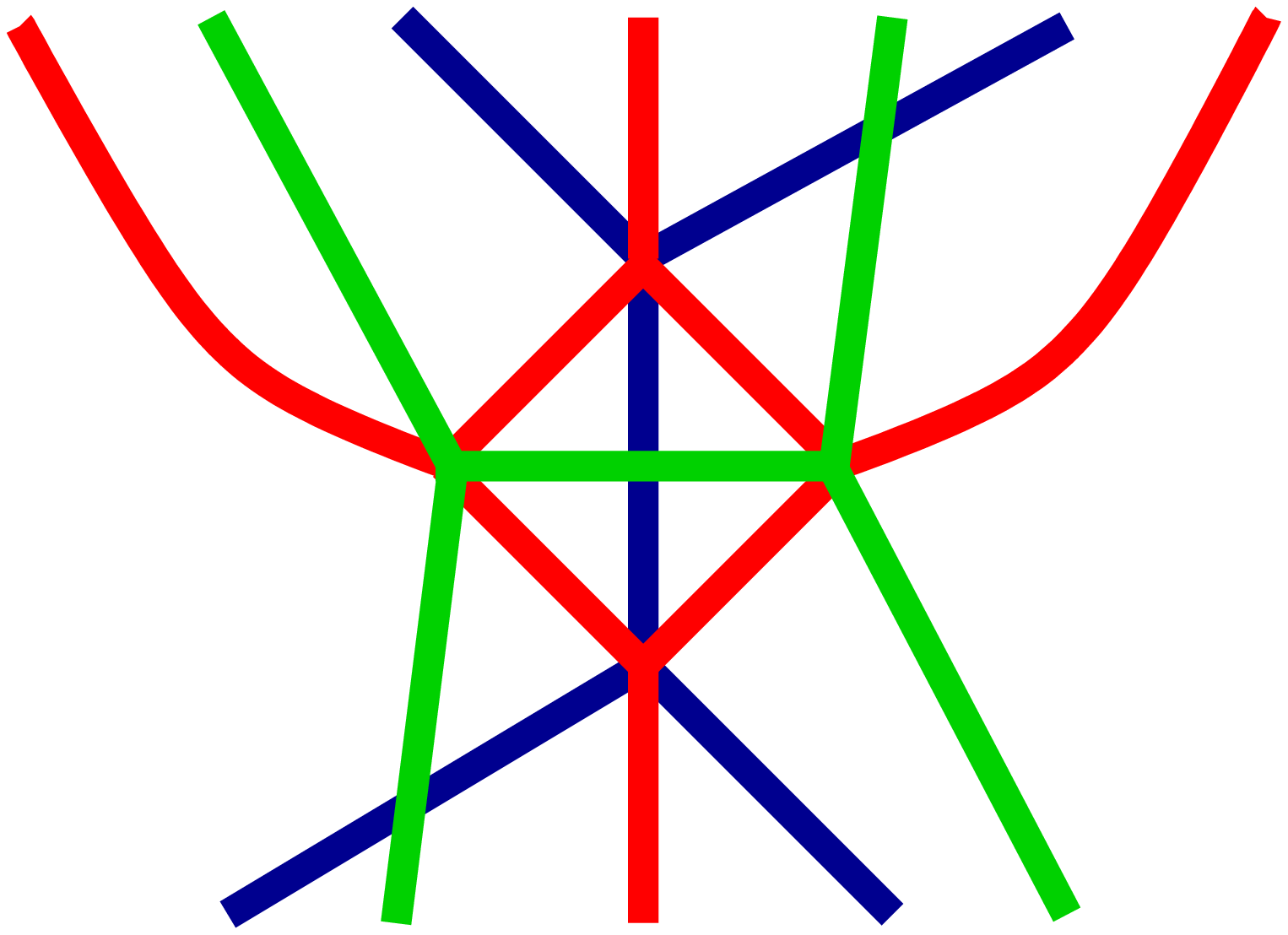}\
=\
\figins{-26}{0.75}{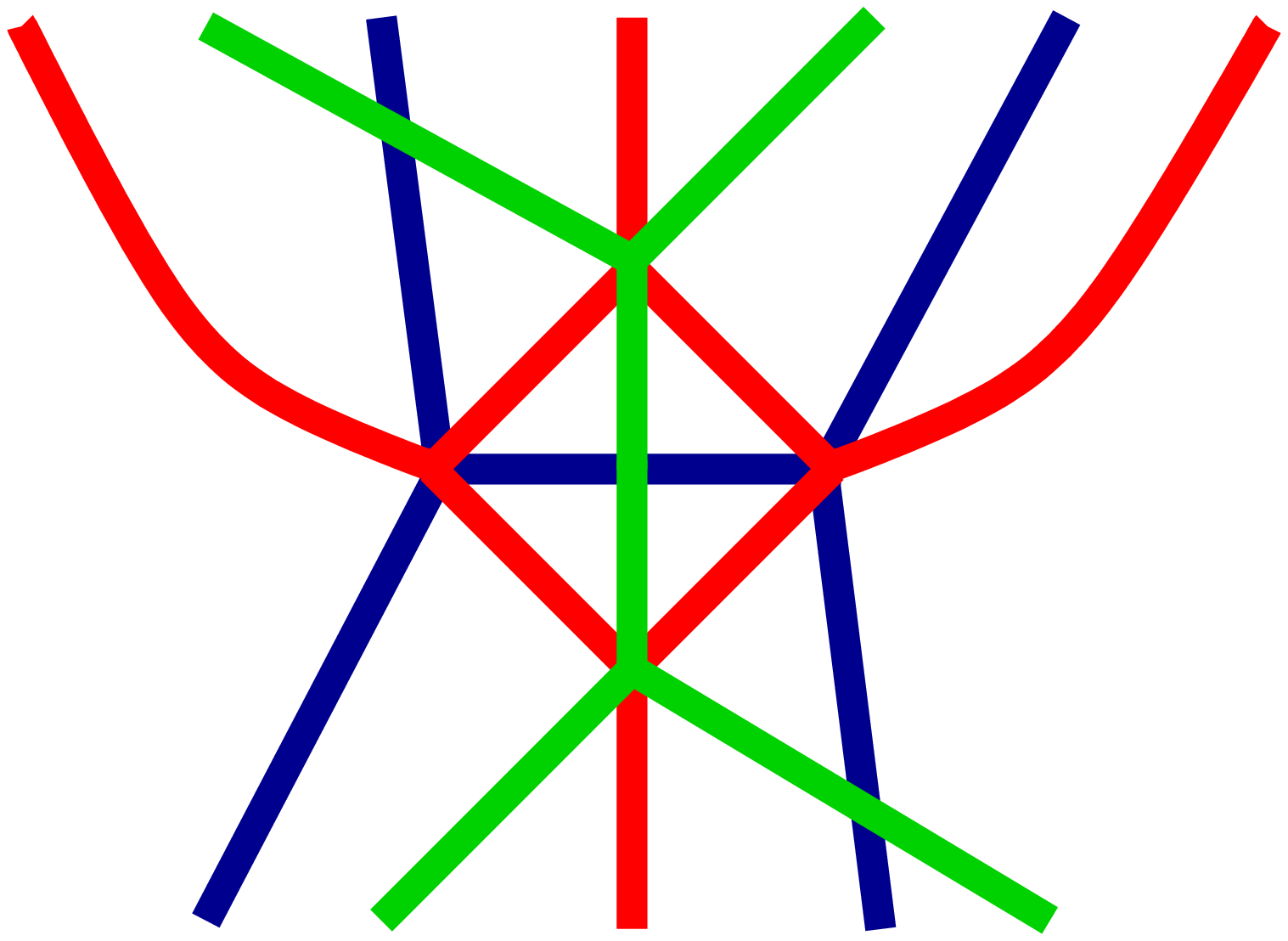}\ .
\end{equation}
We claim that $\Fmv$ sends both sides to zero.
Since the images of both sides of~\eqref{eq:dumbdumbsquare-b} can be obtained from each other 
using a symmetry relative to a vertical plane placed between the sheets labelled $j+1$  
and $j+2$ it suffices to show that one side of~\eqref{eq:dumbdumbsquare-b} is sent to zero. 
The foams involved are rather complicated and hard to visualize. To make the computations easier we 
use movies (two dimensional diagrams) for the whole foam and implicitly translate some bits to 
three dimensional foams to apply isotopy equivalences or relations from Subsection~\ref{ssec:sl3}. 
The r.h.s. corresponds to
\begin{equation*}
f_1=
\figins{-28}{0.8}{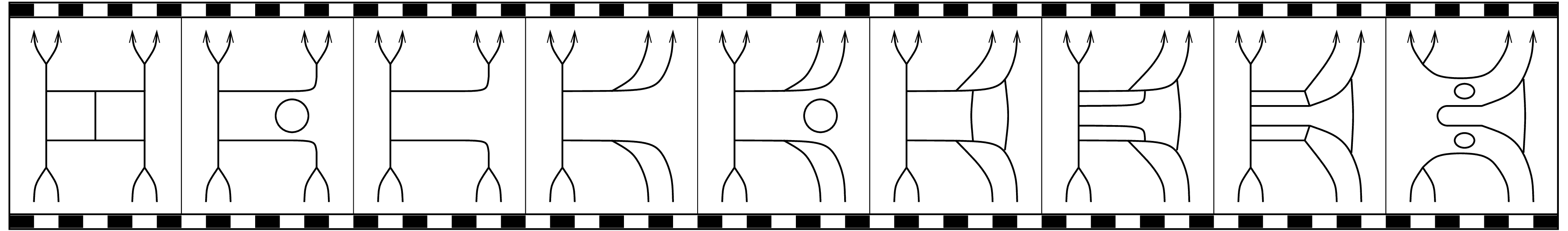}
\end{equation*}
followed by
\begin{equation*}
f_2=
\figins{-28}{0.8}{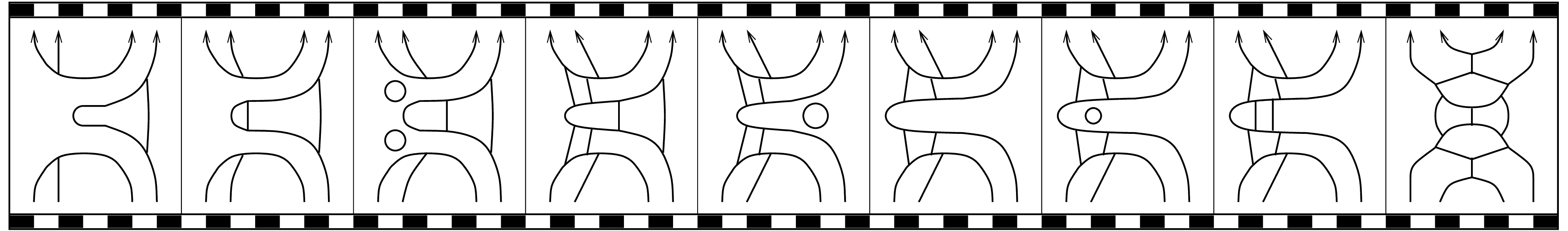}.
\end{equation*}

The foam $f_2$ is isotopic to
\begin{equation*}
%f_3=
\figins{-28}{0.8}{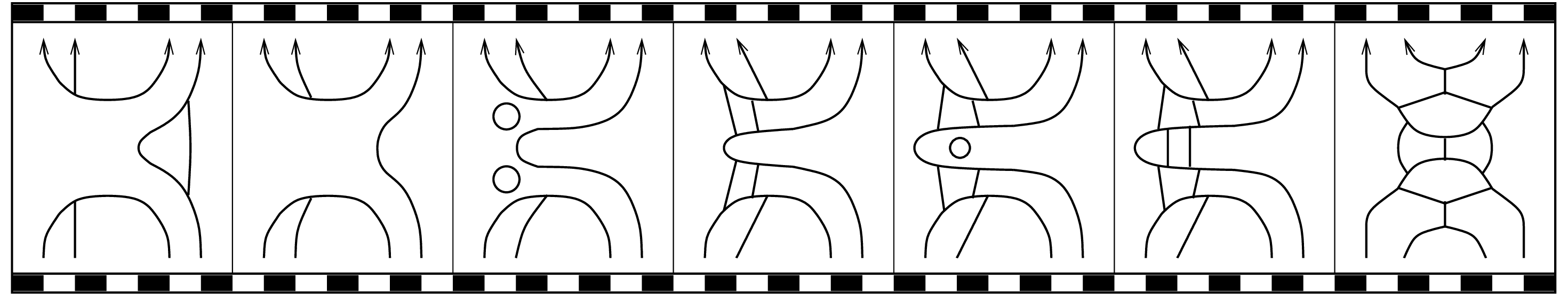}.
\end{equation*}

Using this we can also see that the foams corresponding with
\begin{equation*}
\figins{-28}{0.8}{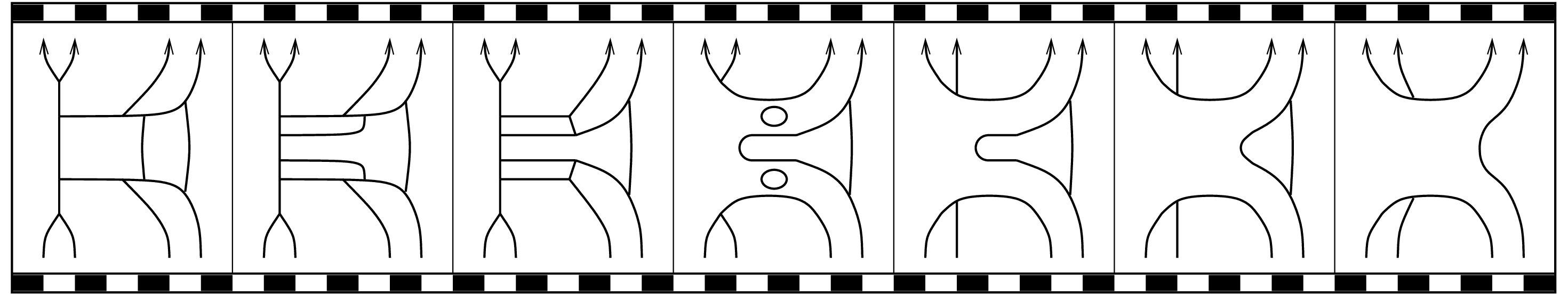}
\end{equation*}
and
\begin{equation*}
\figins{-28}{0.8}{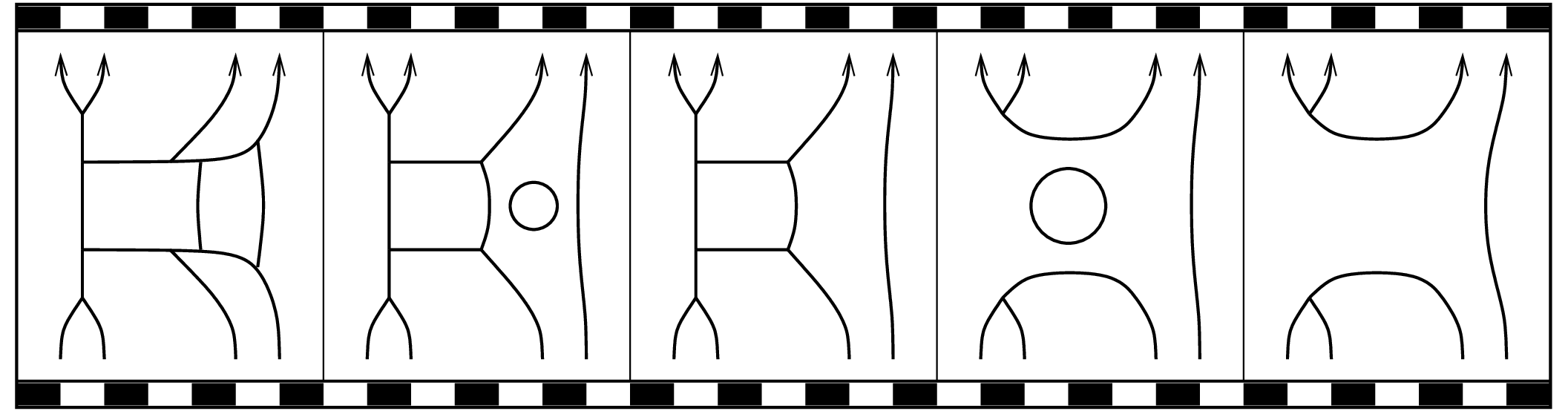}\ 
\end{equation*}
are isotopic.
We see that the foam we have contains  
\begin{equation*}
\figins{-28}{0.8}{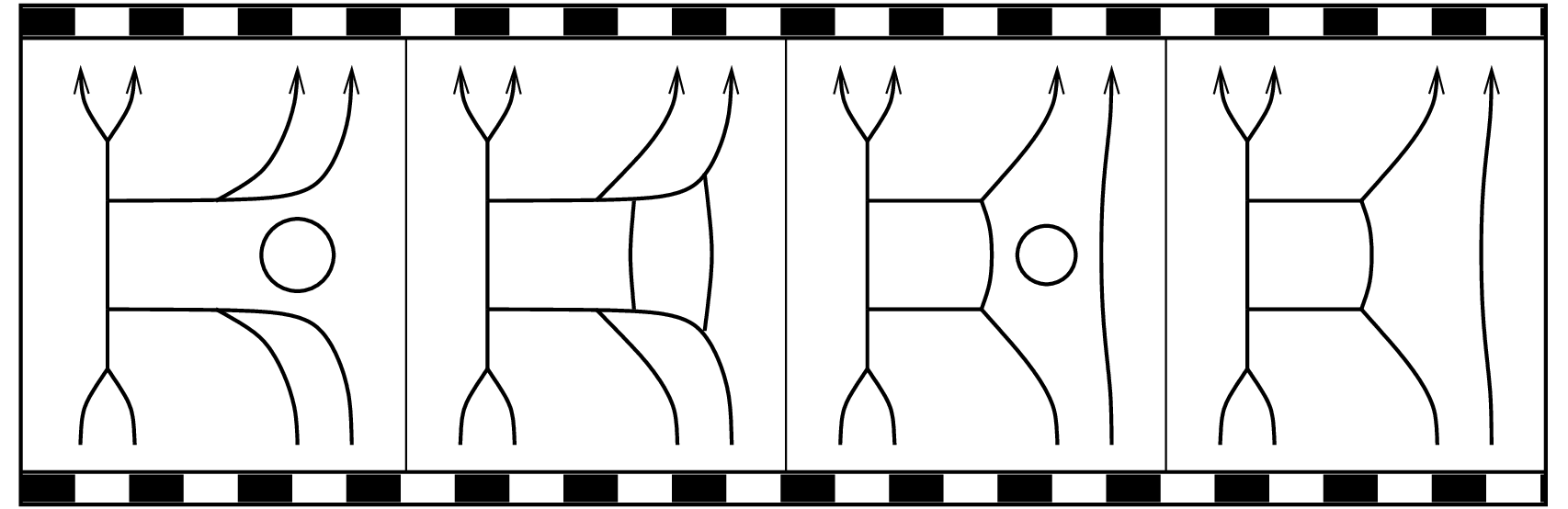}\ ,
\end{equation*}
which corresponds to a foam containing
%\begin{equation*}
$\figins{-8}{0.3}{bbubble00}$, 
%\end{equation*}
which is zero by the~\eqref{eq:bubble} relation.
\end{proof}

%%%%%%%%%%%%%%%%%%%%%%%%%%%%%%%%%%%%%%%%%%%%%%%%%%%%%%%

\vspace*{1cm}

\noindent {\bf Acknowledgements} 
The author would like to thank Mikhail Khovanov and Ben Elias for valuable comments on a 
previous version of this paper. The author also thank Ben Elias for sharing 
reference~\cite{Elias}. 
%The author was supported by the Funda\c {c}\~{a}o para a Ci\^{e}ncia e a Tecnologia (ISR/IST plurianual funding) through the programme ``Programa Operacional Ci\^{e}ncia, Tecnologia, Inova\-\c{c}\~{a}o'' (POCTI) and the POS Conhecimento programme, cofinanced by the European Community fund FEDER.
The author was financially supported by FCT (Portugal), post-doc grant number SFRH/BPD/46299/2008.

%%%%%%%%%%%%%%%%%%%%%%%%%%%%%%%%%%%%%%%%
%%%                                  %%%
%%%            Bibliography          %%%
%%%                                  %%%
%%%%%%%%%%%%%%%%%%%%%%%%%%%%%%%%%%%%%%%%

%%%%%%%%%%%%%%%%%%% end of paper %%%%%%%%%%%%%%%%%%%%%%%%%
\end{document}